\documentclass[11pt]{amsart} 

\usepackage{amssymb,amsmath,amscd,amsxtra}
\usepackage[margin=4cm]{geometry}
\usepackage[all]{xy}

\newtheorem{theorem}{Theorem}[section]
\newtheorem{lemma}[theorem]{Lemma}
\newtheorem{claim}[theorem]{Claim}
\newtheorem{proposition}[theorem]{Proposition}  
\newtheorem{corollary}[theorem]{Corollary} 

\theoremstyle{definition}
\newtheorem{remark}[theorem]{Remark}

\newtheorem{definition}[theorem]{Definition}

\newcommand{\beq}{\begin{equation}}
\newcommand{\eeq}{\end{equation}}
\newcommand{\beqa}{\begin{eqnarray}}
\newcommand{\eeqa}{\end{eqnarray}}
\newcommand{\beaa}{\begin{eqnarray*}}
\newcommand{\ben}{\begin{eqnarray*}}
\newcommand{\eaa}{\end{eqnarray*}}
\newcommand{\een}{\end{eqnarray*}}

\newcommand{\CC}{\mathbb{C}}
\newcommand{\C}{\mathcal{C}}
\newcommand{\PP}{\mathbb{P}}
\newcommand{\A}{\mathcal{A}}
\newcommand{\D}{\mathcal{D}}
\renewcommand{\O}{\mathcal{O}}
\newcommand{\RR}{\mathbb{R}}
\newcommand{\T}{\mathcal{T}}
\newcommand{\U}{\mathcal{U}}
\newcommand{\ZZ}{\mathbb{Z}}

\def\diag{\mathop{\rm diag} \nolimits}

\begin{document} 

\title[Painleve property for semi-simple Frobenius manifolds]
{\textbf{Lectures on Painleve property for semi-simple Frobenius manifolds }}
\author{Todor Milanov}
\address{Kavli IPMU (WPI), UTIAS, The University of Tokyo, Kashiwa, Chiba 277-8583, Japan}
\email{todor.milanov@ipmu.jp} 

\begin{abstract}
These notes are based on a sequence of five lectures given to graduate
students.  The main goal is to prove the so-called Painleve property
for semi-simple Frobenius manifolds. 
\end{abstract}
\maketitle

\setcounter{tocdepth}{2}
\tableofcontents

\section{Introduction}

The Painleve property for semi-simple Frobenius manifolds can be
stated as follows. Let $U$ be a contractible open subset of the configuration
space
\ben
Z_N=\{u\in \CC^N\ :\ u_i\neq u_j \mbox{ for } i\neq j\}.
\een 
As we will see later on, in order to define a semi-simple Frobenius
structure on $U$ we have to choose a holomorphic 1-form
$\sum_{i=1}^N\eta_i(u)du_i$ on $U$ such that 
\begin{enumerate}
\item[(i)] $\eta_i(u)\neq 0$ for all $i$ and for all $u\in U$.
\item[(ii)] The 1-form is closed, i.e., $\eta_{ab}:=\partial
  \eta_a/\partial u_b$ is symmetric in $a$ and $b$. 
\item[(iii)] The functions $\eta_i(u)$ ($1\leq i\leq N$) satisfy the
  following system of PDEs
\ben
(\partial_{u_1}+\cdots +\partial_{u_N})\eta_i  & = 0  & (1\leq i\leq N)\\
(D+u_1\partial_{u_1}+\cdots +u_N\partial_{u_N})\eta_i & =0  & (1\leq i\leq N)\\
\frac{\partial \eta_{ij}  }{ \partial u_k} - \frac{1}{2} \Big(
\frac{\eta_{ij}\eta_{kj} }{\eta_j} + 
\frac{\eta_{jk}\eta_{ik} }{\eta_k}+
\frac{\eta_{ki}\eta_{ji} }{\eta_i} 
\Big) & =0 & (k\neq i\neq j\neq k),
\een
where $D$ is some constant and $\partial_{u_a}:=\partial/\partial u_a$.
\end{enumerate}
The geometric interpretation of the above conditions is the following. Let us
define a diagonal bi-linear pairing on the tangent bundle $TU$ via
\ben
(\partial_{u_i},\partial_{u_j})=\delta_{ij}\eta_j(u).
\een
Condition (i) is equivalent to requiring that the metric is
non-degenerate. Condition (ii) is equivalent to requiring that the
translation vector field $e=\partial_{u_1}+\cdots +\partial_{u_N}$ is
flat with respect to the Levi--Civita connection $\nabla^{\rm L.C.}$
of the pairing. Finally, the system of PDEs in (iii) is equivalent
to requiring that the metric is translation invariant, conformal
invariant, and flat. The main goal of these lectures is to prove the
following theorem.
\begin{theorem}\label{t1}
Suppose that $U$ is lifted to the universal cover $T$ of $Z_N$. Then
every 1-form $\sum_{i=1}^N \eta_i(u)du_i$ satisfying the above
conditions (i)--(iii) extends to a global meromorphic 1-form on
$T$. 
\end{theorem}
Recall that a system of differential equations defined in some domain
$T$ is said to have the Painleve property on $T$ if the so-called movable
singularities (singularities depending on the initial conditions) are
at most poles. That is why the statement of Theorem \ref{t1} is
usually referred to as the Painleve property for semi-simple Frobenius
manifolds. 

Theorem \ref{t1} is a corollary of the theory of isomonodromic
deformations. In fact, there are two possible ways to obtain a proof.
They correspond to the fact that a semi-simple
Frobenius manifold can be obtained as a solution to two different
Riemann--Hilbert (RH) problems (see \cite{Du}). The first RH problem consists of finding a connection on
$\PP^1$ with one Fuchsian singularity at, say $\infty$, and one
irregular singular point at, say $0$. The second one is a RH problem
for a Fuchsian connection. The two connections are related via a
formal Laplace transform, so in principle one could switch between the
two languages. Using the first RH problem, the proof of Theorem
\ref{t1} follows from the results of T. Miwa in \cite{Miw}, while the
second RH problem reduces the proof to the results of Malgrange in
\cite{Mal}.  Therefore, if one wants to understand the statement of Theorem
\ref{t1}, then the main task is to understand either the results of Miwa
or the results of  Malgrange. 

In these lectures we are going to pursue the second approach based on
Malgrange results. The most difficult part in \cite{Mal} is a certain
theorem about existence of a meromorphic trivialization of a family of
vector bundles on $\PP^1$ (see Theorem \ref{thm:Malgrange-1} in these
notes). On the other hand, A. Bolibruch found an elementary proof,
so combining the work of Bolibruch and Malgrange 
we can obtain an argument that requires only basic knowledge of
complex geometry (e.g. chapter 0 in \cite{GH}) and ordinary
differential equations (e.g. \cite{Ar}). 

Summarizing, our lectures include Bolibruch's proof of the
Birkhoff--Grothendieck theorem with parameters, Malgrange's proof of
the Painleve property for the Schlesinger equations, and finally we
check that Theorem \ref{t1} is a corollary of the 
Painleve property for the Schlesinger equations. The last part can be
found also in \cite{Du,Man}.  We have also included a very interesting
theorem of Manin (see \cite{Man}) classifying the initial conditions
for the Schlesinger equations such that the corresponding solution of
the Schlesinger equations comes from a semi-simple Frobenius
manifold. 

We have to admit that many of the arguments in these lectures although elementary
are a bit cumbersome. They could be made more elegant if one is
willing to use some more advanced (but standard) techniques in
complex geometry. Fortunately, this goal is achieved by C. Sabbah. For
more details we refer to the excellent book \cite{Sab}.

\section{Levelt's theory for Fuchsian connections}

The main goal of this lecture is to prove the existence of weak Levelt solutions for Fuchsian systems. We follow closely \cite{Bol}. 
\subsection{Fuchsian systems}
Let $D=\{\lambda\ :\ |\lambda|<r_0\}$ be the open disk of radius $r_0$ and $B_0(\lambda)\in \mathfrak{gl}(\CC^p)$ be a $p\times p$-matrix whose entries depend holomorphically on $\lambda \in D$. We will be interested in the system of ODEs defined by
\ben
\frac{\partial y}{\partial \lambda} (\lambda)=B(\lambda)y(\lambda),\quad B(\lambda):=B_0(\lambda)/\lambda.
\een
Systems of this type are said to be {\em Fuchsian in a neighborhood of $0$}. 

Let us fix a small sector $S$ in $D$ containing the open interval $(0,r_0)$, e.g., 
\ben
S=\{\lambda\in D-\{0\} \ :\ 
-\epsilon < \operatorname{Arg}(\lambda) <\epsilon\}
\een
where $0<\epsilon<2\pi$ is fixed arbitrary. Furthermore, let us fix a reference point $\lambda_0\in (0,r_0)\subset S$ and denote by $X$ the space of holomorphic functions $y: S\to \CC^p$ that solve the above system. The general theory of ODEs implies that $X$ is a finite dimensional vector space of dimension $p$. More precisely
\ben
X\cong \CC^p,\quad y\mapsto y(\lambda_0).
\een
Since the coefficients of the linear system are holomorphic in $D-\{0\}$, every solution $y\in X$ can be extended analytically along any path in $D-\{0\}$. In particular, we have a linear map 
\ben
\sigma:X\to X
\een
corresponding to analytic continuation along a loop based at $\lambda_0$ that goes once around $\lambda=0$ in counter-clockwise direction.

\subsection{Fuchsian singularities are regular}
The following result is well known in the theory of
ODEs. Nevertheless, we give our own proof as well. For a different
argument, which is shorter but yields a slightly weaker result, see
\cite{Bol},  Theorem 4.1 and Lemma 4.1. 
 
\begin{proposition}\label{reg-sing}
Every solution $y\in X$ has the form
\ben
y(\lambda) = \sum_\rho \sum_{k=0}^{p-1} y_{\rho,k}(\lambda) \lambda^\rho (\log\lambda)^k,
\een
where the first sum is over all eigenvalues $\rho$ of $B_0(0)$ and $y_{\rho,k}(\lambda)$ are $\CC^p$-valued functions analytic for all $\lambda\in D$.
\end{proposition}
\proof
We will prove that the system has a fundamental matrix whose columns have the above form. Using a constant gauge transformation $y(\lambda)\mapsto Cy(\lambda)$ we can reduce the general case to the case when $B_0(0)$ is in Jordan normal form. Moreover, we may assume that the Jordan blocks are ordered in such a way that $B_0(0)=R+N^{(0)}$, where $R$ and $N^{(0)}$ have the following properties. Both
\ben
R=\operatorname{diag}(R_1,\dots,R_s)
\quad \mbox{and}\quad 
N^{(0)}=\operatorname{diag}(N^{(0)}_1,\dots,N^{(0)}_s)
\een
are block-diagonal. The block $R_i=\rho_i \operatorname{I_i}$, where $I_i$ is an identity matrix of size the multiplicity of $\rho_i$ as an eigenvalue of $B_0(0)$ and 
\ben
\operatorname{Re}(\rho_1)>\dots >\operatorname{Re}(\rho_s).
\een
The block $N_i^{(0)}$ ($1\leq i\leq s$) is an upper-triangular nilpotent matrix whose size is the same as the the size of $I_i$. Note that the commutator $[R,N^{(0)}]=0$.
\begin{claim}\label{claim:1}
There exists a formal solution 
\ben
Y(\lambda) = U(\lambda)\lambda^R \lambda^N,
\een
where 
$$
U(\lambda)=1+U_1\lambda+U_2 \lambda^2+\cdots,\quad 
U_k\in \mathfrak{gl}(\CC^p)
$$ 
and $N$ is upper-triangular nilpotent matrix of the form
$$
N=N^{(0)}+N^{(1)}+\cdots,\quad [R,N^{(k)}]=kN^{(k)}.
$$
\end{claim}
\proof
Put $B_0(\lambda)=B_{0,0}+B_{0,1}\lambda+\cdots$, substitute $Y(\lambda)$ in the differential equation, and compare the coefficients in front of the powers of $\lambda^k$. For $k=0$ we get
$R+N^{(0)}=B_0(0)=B_{0,0}$, which is true by definition. For $k>0$ we get
\ben
kU_k+[U_k,R+N^{(0)}]+N^{(k)} = 
B_{0,k} +\sum_{i=1}^{k-1} \Big(
B_{0,k-i}U_i-U_i N^{(k-i)} \Big). 
\een 
The linear operator 
\ben
\operatorname{ad}_R:\mathfrak{gl}(\CC^p)\to \mathfrak{gl}(\CC^p), 
\quad
x\mapsto [R,x]
\een
is diagonalizable, i.e., we have a decomposition 
\ben
\mathfrak{gl}(\CC^p) = \bigoplus_{a\in {\rm spec}(R)} 
\mathfrak{gl}_a(\CC^p),
\een
where ${\rm spec}(R)$ denotes the set of eigenvalues of $\operatorname{ad}_R$ and for $a\in {\rm spec}(R)$
\ben
\mathfrak{gl}_a(\CC^p) = \{x\ :\ [R,x]=ax\}
\een
is the corresponding eigen-subspace. Let us denote by 
$$\pi_a: \mathfrak{gl}(\CC^p)\to \mathfrak{gl}_a(\CC^p)$$
the projection map defined via the above decomposition. 
Let us assume that we have determined $U_1,\dots,U_{k-1}$ and $N^{(1)},\dots,N^{(k-1)}$. Then $U_k=\sum_{a\in {\rm spec}(R)} \pi_a(U_k)$ and $N^{(k)}\in \mathfrak{gl}_k(\CC^p)$ are defined by projecting via $\pi_a$ the above recursion relation and solving for $\pi_a(U_k)$ and $\pi_a(N^{(k)})$. There are two cases. First, if $a=k$, then we set $\pi_k(U_k)=0$. Note that since $N^{(0)}$ commutes with $R$, we have $\pi_k([U_k,N^{(0)}]) =[\pi_k(U_k),N^{(0)}]=0$ and $\pi_k(N^{(k)})=N^{(k)}$. Therefore, we can uniquely solve for $N^{(k)}$. The second case is if $a\neq k$, then $\pi_a(N^{(k)})=0$ and
\ben
\pi_a(kU_k + [U_k,R+N^{(0)}]) = (k-a-\operatorname{ad}_{N^{(0)}} )\pi_a(U_k).
\een
Since $N^{(0)}$ is nilpotent, the linear operator $\operatorname{ad}_{N^{(0)}}$ is also nilpotent. Therefore the linear operator  $k-a-\operatorname{ad}_{N^{(0)}}$ is invertible, so we can uniquely solve for $\pi_a(U_k)$. 
\qed

\medskip
It remains to prove that the formal series $U(\lambda)$ is convergent. 
Note that $U(\lambda)$ satisfies the following differential equation
\beq\label{deq-U}
(\lambda\partial_\lambda - \operatorname{ad}_{R+N^{(0)}}) U =
(U\alpha(\lambda) +  \beta(\lambda) U ),
\eeq
where
\ben
\alpha(\lambda) = -\sum_{i=1}^\infty N^{(i)} \lambda^i
\quad
\beta(\lambda) = \sum_{i=1}^\infty B_{0,i}\lambda^i.
\een
Let us fix an integer $k>0$, such that the set ${\rm spec}(R)$ does not contain any integers $\ell>k$. Note that $N^{(\ell)}=0$ for all $\ell>k$, so $\alpha(\lambda)$ is polynomial in $\lambda$. Let us write the formal series in the form
\ben
U(\lambda) = U_{\leq k}(\lambda) + \lambda^k V(\lambda),\quad
U_{\leq k}(\lambda) = 1+\sum_{i=1}^k U_i\lambda^i,
\een
where $V(\lambda) = \sum_{j=1}^\infty U_{j+k} \lambda^j$. Then $V(\lambda)$ satisfies the following differential equation 
\beq\label{deq-V}
(\lambda\partial_\lambda +k- \operatorname{ad}_{R+N^{(0)}}) V = 
V \alpha(\lambda)  + \beta(\lambda) V + \gamma(\lambda),
\eeq
where 
\ben
\gamma(\lambda) = \lambda^{-k} \Big(
U_{\leq k}(\lambda) \alpha(\lambda) + \beta(\lambda) U_{\leq
  k}(\lambda) -
(\lambda\partial_\lambda - \operatorname{ad}_{R+N^{(0)}})
U_{\leq k}(\lambda)\Big)
\een
By definition $U_{\leq k}(\lambda)$ satisfies the differential equation \eqref{deq-U} up to terms of order $O(\lambda^{k+1})$. Therefore, $\gamma(\lambda)$ is analytic at $\lambda=0$ and $\gamma(0)=0$.
It is enough to prove that the differential equation \eqref{deq-V} has a solution $V_{hol}(\lambda)$ analytic at $\lambda=0$. Indeed, the linear operator $k-\operatorname{ad}_{R+N^{(0)}}$ is invertible, so after substituting the Taylor series of $V_{hol}(\lambda)$ in the differential equation we get that the Taylor series must coincide with the formal series $V(\lambda)$.

In order to construct a holomorphic solution, we use the standard idea to identify $V_{hol}$ with the fixed point of a certain integral operator. Let us fix a closed disk $D_r=\{\lambda\ :\ |\lambda|\leq r \}$ with radius $r<r_0$. Let us define a sequence of holomorphic $\mathfrak{gl}(\CC^p)$-valued functions 
$$
V_n: D_r \to \mathfrak{gl}(\CC^p),\quad n=0,1,2,\dots
$$
as follows. Put $V_0(\lambda)=0$ and let $V_{n+1}(\lambda)$ be such that 
\ben
(\lambda\partial_\lambda +k-\operatorname{ad}_{R+N^{(0)}}) V_{n+1} = V_n \alpha(\lambda) + \beta(\lambda) V_n +\gamma(\lambda). 
\een
Note that 
\ben
V_{n+1}(\lambda) = \int_0^1 t^{k-ad_{R+N^{(0)}}} 
\Big(
V_n(t\lambda) \alpha(t\lambda) + \beta(t\lambda) V_n(t\lambda) +\gamma(t\lambda)
\Big)\frac{dt}{t}.
\een
The convergence of the integral follows from the fact that if we choose $k$ sufficiently large the real part of the eigenvalues of $k-ad_{R+N^{(0)}}$ will be positive. Therefore $V_{n+1}(\lambda)$ is an analytic function for all $\lambda\in D_r$.

In order to prove that the sequence $V_n$ is convergent we introduce the following norm. Let $|\ |:\mathfrak{gl}(\CC^p)\to \mathbb{R}_{\geq 0}$ be the standard matrix norm
\ben
|A| = \operatorname{sup}_{v\neq 0} \frac{|Av|}{|v|},
\een 
where $|v|=\sqrt{|v_1|^2+\cdots +|v_p|^2}$ is the standard Euclidean norm of $v\in \CC^p$. If $A:D_r\to \mathfrak{gl}(\CC^p)$ is holomorphic, then we define 
\ben
|\!|A|\!|_r = \sum_{i=0}^\infty |A_i| r^i,
\een
where $A(\lambda) = \sum_{i=0}^\infty A_i \lambda^i$ is the Taylor
series expansion. Let $B_r$ be the space of those holomorphic maps
$A$ for which $|\!|A|\!|_r<\infty$. It is known (see \cite{GR})
that $B_r$ is a Banach algebra. Using the Cauchy 
inequality it is easy to prove that if $A(\lambda)$ is holomorphic for
all $\lambda\in D$ then $A\in B_r$. 
\begin{claim}
Suppose $k>|\operatorname{ad}_{R+N^{(0)}}|$. Then the map 
\ben
F: B_r\to B_r,\quad
F(A)(\lambda):=\int_0^1 t^{k-ad_{R+N^{(0)}}} A(t\lambda)\frac{dt}{t}
\een
is a bounded linear operator of norm less or equal to $ 1$, i.e., 
$
|\!|F(A)|\!|_r \leq  |\!|A|\!|_r.
$ 
\end{claim}
\proof
Put $A(\lambda)=\sum_{i=0}^\infty A_i\lambda^i$. Then the coefficient in front of $\lambda^i$ in $F(A)$ is
\ben
F(A)_i = \int_0^1 t^{k+i-1-{\rm ad}_{R+N^{(0)}}} A_i dt.
\een
Using that 
\ben
|t^{k+i-1-{\rm ad}_{R+N^{(0)}}}| = 
t^{k+i-1}|t^{-{\rm ad}_{R+N^{(0)}}}|\leq 
t^{k+i-1}t^{-|{\rm ad}_{R+N^{(0)}}|},\quad 0\leq t \leq 1
\een
we get 
\ben
|F(A)_i|\leq \frac{|A_i|}{k+i-|{\rm ad}_{R+N^{(0)}}|} \leq |A_i|.
\qed
\een
Note that $V_{n+1}=F(V_n\alpha + \beta V_n +\gamma)$. Therefore
\ben
|\!|V_{n+1}-V_n|\!|_r \leq 
(|\!|\alpha|\!|_r +|\!|\beta |\!|_r) \,  |\!|V_n-V_{n-1}|\!|_r.
\een
Since $\alpha(0)=\beta(0)=0$ we can always choose $r$ so small that
$|\!|\alpha|\!|_r +|\!|\beta |\!|_r<1$. Then the above
inequality shows that $\{V_n\}$ is a Cauchy sequence in $B_r$, so
the limit $V_{hol}=\lim_{n\to \infty} V_n $ exists and it gives a
solution to the differential equation \eqref{deq-V}.

Finally, note that the series $U(\lambda)$ must be analytic for all
$\lambda\in D$, because the fundamental matrix
$Y(\lambda)=U(\lambda)\lambda^R\lambda^N$ extends analytically along
any path inside $D-\{0\}$. 
\qed 
\begin{corollary}\label{cor:mon-nilp}
If $B_0(0)$ is nilpotent, then the matrix of the monodromy of the
Fuchsian system with respect to a basis of $X$ given by the columns of
the fundamental matrix $Y(\lambda)$ satisfying the initial condition
$Y(\lambda_0)=1$ is $e^{2\pi\sqrt{-1} B_0(0)}$.  
\end{corollary}

\subsection{Levelt evaluations}
Let us denote by $\O[S]$ the space of holomorphic maps $y:S\to \CC^p$,
such that 
\ben
\lim_{\substack{\lambda\to 0\\ \lambda\in S}}
\frac{y(\lambda)}{|\lambda|^m}=0
\een
for some integer $m$. Such functions are also sometimes said to be of
{\em moderate growth} at $\lambda=0.$  The key to Levelt's theory is the
map
\ben
\varphi:\O[S]\to \mathbb{Z}\cup \{\infty\}
\een
defined by 
\ben
\varphi(y) := \operatorname{max}\Big\{ m\in \mathbb{Z}\ |\
\lim_{\substack{\lambda\to 0\\ \lambda\in S}}
\frac{y(\lambda)}{|\lambda|^\ell}=0\mbox{ for all } \ell<m\Big\}
\een
for all $y\in \O[S]\setminus{\{0\}}$ and $\varphi(0)=\infty$. Note
that according to Proposition \ref{reg-sing} the space of solutions
$X\subset \O[S]$. 

\begin{lemma}\label{lemma:eval} The map $\varphi$ satisfies the
  following properties.
\begin{enumerate}
\item[a)] If $y_1,y_2\in \O[S]$, then $\varphi(y_1+y_2) \geq
\operatorname{min}(\varphi(y_1),\varphi(y_2)$.  If
$\varphi(y_1)\neq \varphi(y_2)$, then the equality in the above
inequality holds.
\item[b)] If $c\in \CC\setminus{\{0\}}$, then $\varphi(cy)=\varphi(y)$ for
all $y\in \O[S]$.
\end{enumerate} 
\end{lemma}
The proof is an elementary consequence from the definitions, so it
will be omitted. 

\begin{lemma}\label{cf-vanish}
Let $\theta_i$, $1\leq i\leq n$ be real numbers such that
$\theta_i\neq \theta_j$ for all $i\neq j$. Suppose that 
\ben
f(x)=\sum_{i=1}^n a_i(x)e^{\sqrt{-1}\theta_i x},\quad 
a_i\in \CC[x]
\een 
and that there is a real number $\epsilon>0$, such that
\ben
\lim_{x\to +\infty} f(x) e^{\ell x}=0,\quad \forall \ell<\epsilon.
\een
Then $a_i=0$ for all $i$. 
\end{lemma}
\proof
It is enough to prove the lemma in the case when the polynomials are
constants. Indeed, let $m$ be the maximal degree
among the degrees of $a_i$, i.e., $a_i(x)=\sum_{\mu=0}^m a_{i,\mu}
x^\mu$ and $a_{i,m}\neq 0$ for at least one $i$. There exists a constant $C$, s.t., 
\ben
\Big|\sum_{i=1}^n a_{i,m} e^{\sqrt{-1}\theta_i x} \Big| x^m e^{\lambda
  x} \leq C 
\Big|\sum_{i=1}^n a_i(x) e^{\sqrt{-1}\theta_i x} 
\Big| e^{\lambda x},\quad \forall x\geq 0.
\een
Therefore we must have 
\ben
\lim_{x\to +\infty} \Big(\sum_{i=1}^n a_{i,m} e^{\sqrt{-1}\theta_i x} \Big) e^{\lambda
  x} =0 \quad \forall \lambda<\epsilon.
\een
Therefore, if we knew that the lemma holds for constant polynomials,
then we would get $a_{i,m}=0$ for all $i$ -- contradiction with the
definition of $m$.

Let us assume that $a_i\in \CC$ are constants. 
Using induction on $m$ it is easy to prove that if $\lambda_m<\cdots
<\lambda_1<\epsilon$ is any sequence of real numbers, then 
\ben
\lim_{x\to +\infty} \Big(\sum_{i=1}^n \frac{a_i e^{\sqrt{-1}\theta_i
    x }}{(\sqrt{-1}\theta_i+\lambda_1)\cdots
  (\sqrt{-1}\theta_i+\lambda_m)}\Big) e^{\lambda x} = 0,\quad \forall \lambda<\lambda_m.
\een
Indeed, the starting point of the induction is $m=0$ and the statement
is true by definition. Suppose the statement is true for  $m$ and that
$\lambda_{m+1}<\lambda_m$. Let us pick $\lambda'$ in the open interval
$(\lambda_{m+1},\lambda_m)$. Using the inductive assumption we get 
\ben
\Big|\sum_{i=1}^n \frac{a_i e^{(\sqrt{-1}\theta_i+\lambda_{m+1})
    y}}{(\sqrt{-1}\theta_i+\lambda_1)\cdots
  (\sqrt{-1}\theta_i+\lambda_m)}\Big|\leq C'\, e^{(\lambda_{m+1}-\lambda') y},\quad
\forall y\geq 0
\een
for some constant $C'$ depending on the choice of $\lambda'$. Integrating the function inside the absolute
value on the LHS for $y$ from $0$ to $x$ and using the above
inequality to estimate the absolute value of the integral, we get the
following inequality  
\ben
\Big|\sum_{i=1}^n \frac{a_i (e^{(\sqrt{-1}\theta_i+\lambda_{m+1})
    x} -1)}{(\sqrt{-1}\theta_i+\lambda_1)\cdots
  (\sqrt{-1}\theta_i+\lambda_m)
  (\sqrt{-1}\theta_i+\lambda_{m+1})}\Big|\leq C'\, 
\frac{e^{(\lambda_{m+1}-\lambda') x}-1}{\lambda_{m+1}-\lambda'}.
\een
If $\lambda<\lambda_{m+1}$ is any given number, then we multiply the above inequality by
$e^{(\lambda-\lambda_{m+1})x}$, and let $x\to +\infty$. 

To complete the proof of the lemma we proceed as follows. Let us
choose a sequence of $n$ numbers $0<\lambda_n<\cdots<\lambda_1<
\epsilon$ and define the matrix $C$ with entries
\ben
C_{im}:= \frac{1 }{ 
(\sqrt{-1}\theta_i+\lambda_1)\cdots
(\sqrt{-1}\theta_i+\lambda_m)},\quad 1\leq i,m\leq n.
\een
Note that for $\lambda_1=\cdots=\lambda_n$ the determinant of $C$ turns
into a Wandermond determinant, which is not 0 according to the
assumption $\theta_i\neq \theta_j$ for $i\neq j$. Therefore choosing
$\lambda_1$ sufficiently close to $\lambda_n$ we may guarantee that $C$ is
invertible. On the other hand if we define
\ben
g_m(x) = \sum_{i=1}^n a_i e^{\sqrt{-1}\theta_i x }C_{im},\quad 1\leq
m\leq n,
\een
then according to the above fact $\lim_{x\to +\infty}
g_m(x)=0$. However, since $C$ is invertible, we can solve the above
equations and express each $a_i e^{\sqrt{-1}\theta_i x }$ as a linear
combination of $g_m(x)$ with constant coefficients. Therefore $\lim_{x\to
  +\infty} a_i e^{\sqrt{-1}\theta_i x }=0$. This however is possible
only if $a_i=0$.
\qed

\begin{proposition}\label{mon-inv}
If $y\in X$, then $\varphi(\sigma y)=\varphi(y)$. 
\end{proposition}
\proof
Recalling Proposition \ref{reg-sing} we write the solution as
\ben
y(\lambda)=\sum_{i=1}^n y_i(\lambda)
\lambda^{\rho_i},\quad y_i(\lambda)=\sum_{k=0}^{p-1}
y_{i,k}(\lambda)\,(\log \lambda)^k
\een
where $y_{i,k}(\lambda)$ are analytic at $\lambda=0$. We may further
assume that  $\operatorname{Re}(\rho_1)\leq \cdots \leq
\operatorname{Re}(\rho_n)$. Let us write the solution as
\ben
y(\lambda)=\lambda^{\rho_1} \Big(
f(\lambda)+\sum_{j} y_j(\lambda)\lambda^{\rho_j-\rho_1} \Big),
\een
where the sum is over all $j$, s.t., that
$\operatorname{Re}(\rho_j)>\operatorname{Re}(\rho_1)$
and 
\ben
f(\lambda)=\sum_i y_i(\lambda)\lambda^{\rho_i-\rho_1},
\een
where the sum is over all $i$, s.t., 
$\operatorname{Re}(\rho_i) = \operatorname{Re}(\rho_1)$. Let us assume
that $y_{1,k}(0)\neq 0$ for at least one $k$. Otherwise, we can
replace $\rho_1$ with an exponent with a larger real part. Note that $\varphi(y)=\lfloor
\operatorname{Re}(\rho_1)\rfloor$, where $\lfloor x\rfloor$ is the
largest integer that does not exceed $x$. Indeed, be definition we
have that if $\ell < \lfloor
\operatorname{Re}(\rho_1) \rfloor$, then $\lim
y(\lambda)/|\lambda|^\ell =0$, so $\varphi(y)\geq
\lfloor \operatorname{Re}(\rho_1) \rfloor$. If the inequality is strict then
we can find $\epsilon >0$, such that 
\ben
\lfloor \operatorname{Re}(\rho_1) \rfloor\leq
\operatorname{Re}(\rho_1) <\operatorname{Re}(\rho_1)+2\epsilon < \varphi(y)
\een
and $\epsilon <\operatorname{Re}(\rho_j)- \operatorname{Re}(\rho_1)$
for all $j$ for which  $\operatorname{Re}(\rho_j)\neq
\operatorname{Re}(\rho_1).$ We have
\ben
\frac{y(\lambda)}{|\lambda|^{\operatorname{Re}(\rho_1)+\epsilon +\ell}}
=
\frac{\lambda^{\rho_1}}{|\lambda|^{\operatorname{Re}(\rho_1)+\epsilon}}\, 
\Big(
\frac{f(\lambda)}{|\lambda|^\ell} + \sum_j y_j(\lambda) 
\frac{\lambda^{\rho_j-\rho_1} }{
|\lambda|^\ell}
\Big).
\een
If $\ell<\epsilon$, then the LHS has limit 0 as $\lambda\to 0$, while
the limit of  the first factor on the RHS is $\infty$ and the limit of
the sum over $j$ is $0$. Therefore, we must have 
\ben
\lim_{\lambda\to 0} \frac{f(\lambda)}{|\lambda|^\ell} = 0,\quad
\forall \ell<\epsilon.
\een
If we put $\lambda=e^{-x}$, $x\in \RR_{>0}$ and let $x\to +\infty$ we
get that
\ben
\sum_i \sum_{k=0}^{p-1} y_{i,k}(0) x^k e^{\sqrt{-1}\theta_i x},\quad
\sqrt{-1}\theta_i = \rho_i-\rho_1
\een 
satisfies the condition of Lemma \ref{cf-vanish},
so it must be $0$, which contradicts the choice of $\rho_1$. 

Note also that we have 
\ben
\sigma y(\lambda) = \sum_{i=1}^n  \sum_{k=0}^{p-1} y_{i,k}(\lambda)
e^{2\pi\sqrt{-1} \rho_i} \lambda^{\rho_i} (\log \lambda+2\pi\sqrt{-1} )^{k}.
\een
Therefore, choosing $k$ to be the largest integer such that
$y_{1,k}(0)\neq 0$ we get
\ben
\varphi(y) = 
\varphi(y_{1,k}(\lambda) \lambda^{\rho_1} (\log \lambda)^{k}) =
\varphi(y_{1,k}(\lambda) e^{2\pi\sqrt{-1} \rho_1}\lambda^{\rho_1} 
(\log \lambda+ 2\pi\sqrt{-1})^{k}) \leq \varphi(\sigma y),
\een
where in the last equality we used Lemma \ref{lemma:eval}, Part a). Similarly $\varphi(y)\leq
\varphi(\sigma^{-1}y)$ for all $y\in X$.  Finally we get
\ben
\varphi(y)\leq \varphi(\sigma y)\leq \varphi(\sigma^{-1}(\sigma y)) =
\varphi(y).
\qed
\een

\subsection{Weak Levelt solutions}

The eigenvalues of $\sigma$ can be written uniquely as 
\ben
e^{2\pi\sqrt{-1} \rho_i},\quad 
0\leq \operatorname{Re}(\rho_i)<1,\quad 
1\leq i\leq s. 
\een
Let 
\ben
X=X_1\oplus \cdots \oplus X_s,\quad 
X_i:=\{y\in X\ :\ (\sigma-e^{2\pi\sqrt{-1} \rho_i})^ny=0
\mbox{ for all } n\gg 0\}
\een
be the decomposition of $X$ into generalized eigensubspaces.

Using Lemma \ref{lemma:eval} we get that $\varphi(X)$ is a finite
set. Let us define the set
\ben
\{\infty,k_i^{1},\dots,k_i^{m_i}\} :=\varphi(X_i),\quad 
1\leq i\leq s,
\een
where in addition we assume that $k_i^1>\cdots >k_i^{m_i}$. Put
\ben
X_i^\ell = \{y\in X\ |\ \varphi(y)\geq k_i^\ell\},\quad 
1\leq i \leq s,\quad 1\leq \ell \leq m_i.
\een
According to Lemma \ref{lemma:eval} the sets $X_i^\ell$ are vector
subspaces of $X_i$, so we have a strictly increasing filtration (in
particular we see that $\varphi$ could take only finitely many values
on $X_i$)
\ben
X_i^1\subset X_i^2\subset \cdots \subset X_i^{m_i} = X_i
\een
Using Proposition \ref{mon-inv} we get that the above filtration is
$\sigma$-invariant. 

A {\em weak Levelt} solution $Y(\lambda)$ is by definition a
fundamental matrix whose columns are splited into $s$ groups
\ben
Y(\lambda)=[Y_1(\lambda)\, \cdots \, Y_s(\lambda)],
\een 
where the columns in $Y_i(\lambda)$ represent a basis of $X_i$ with
the following property. We can split $Y_i(\lambda)$ into $m_i$ groups
\ben
Y_i(\lambda)=[Y_{i,1}(\lambda)\, \cdots \, Y_{i,m_i}(\lambda)]
\een
such that 
\begin{enumerate}
\item[(i)]
The columns in $Y_{i,\ell}(\lambda)$ represent a basis of the
quotient subspace $X_i^\ell/X_i^{\ell-1}$. 
\item[(ii)] 
The matrix of the linear operator in $X_i^\ell/X_i^{\ell-1}$ induced
by $\sigma$ in the basis represented by the columns of
$Y_{i,\ell}(\lambda)$ is upper-triangular.
\end{enumerate}
Let $G$ be the matrix of $\sigma$ with respect to the basis of $X$
given by the columns of a weak Levelt solution $Y(\lambda)$. Note that
the matrix $G$ is block-diagonal
\ben
G=\operatorname{diag}(G_1,\dots, G_s)
\een
where each block is a square matrix of size
$\operatorname{dim}_\CC(X_i)$. Each block $G_i$ has a
natural  block-matrix form corresponding to the filtration
$X_i^1\subset \cdots \subset X_i^{m_i}$
\ben
G_i= \begin{bmatrix}
G_i^{11} & G_i^{12} & \cdots & G_i^{1m_i} \\
0 & G_i^{22} & \cdots & G_i^{2m_i} \\
\vdots & \vdots & \ddots &  \vdots \\
0 & 0 & \cdots & G_i^{m_im_i}
\end{bmatrix},
\een
where the size of the block $G_i^{ab}$ is
$
\operatorname{dim}_\CC(X^a_i/X_i^{a-1})\times 
\operatorname{dim}_\CC(X^b_i/X_i^{b-1})$. The definition of a weak
Levelt solution implies that $G_i^{ab}=0$ for
$a>b$ ($\because$ the filtration is $\sigma$-invariant) and that the
block $G_i^{\ell\ell}$ has the form of an upper-triangular 
matrix with all diagonal entries being equal to $e^{2\pi\sqrt{-1}
  \rho_i}$ ($\because$ the matrix of the linear map in
$X_i^\ell/X_i^{\ell-1}$ induced by $\sigma$ is upper-triangular).

\subsection{Levelt's theorem}
Let $Y(\lambda)$ be a weak Levelt solution. Let us write the monodromy
matrix $G=e^{2\pi\sqrt{-1} E},$ where $E$ has the same block-matrix
structure as $G$. Namely, 
\ben
E=\operatorname{diag}(E_1,\dots, E_s)
\een 
is block-diagonal and each block $E_i$ has the form
\ben
E_i= \begin{bmatrix}
E_i^{11} & E_i^{12} & \cdots & E_i^{1m_i} \\
0 & E_i^{22} & \cdots & E_i^{2m_i} \\
\vdots & \vdots & \ddots &  \vdots \\
0 & 0 & \cdots & E_i^{m_im_i}
\end{bmatrix},
\een  
where $E_i^{\ell\ell} = \rho_i I_i^\ell + N_i^{\ell\ell}$ is
upper-triangular matrix whose diagonal entries are all equal to
$\rho_i$. We have denoted by $I_i^\ell$ the identity matrix of size
$\operatorname{dim}_\CC(X^\ell_i/X_i^{\ell-1})$ while by
$N_i^{\ell\ell}$ we have denoted the strictly upper-triangular part of
$E_i^{\ell\ell}$. 

Let us define also the matrix $K$ with the same block-diagonal
structure as $G$ and $E$, i.e., 
\ben
K=\operatorname{diag}(K_1,\dots, K_s)
\een
where the block $K_i$ is given by the diagonal matrix
\ben
K_i=\operatorname{diag}(k_i^1I_i^1,\dots,k_i^{m_i}I_i^{m_i}). 
\een
The main result of this lecture can be stated as follows.
\begin{theorem}[Levelt]\label{thm:Levelt}
Suppose that $Y(\lambda)$ is a weak Levelt solution and that $K$ and
$E$ are the matrices defined as above. Then
\ben
Y(\lambda)=U(\lambda) \lambda^K\lambda^E,
\een
where $U(\lambda)$ is holomorphic  and invertible at $\lambda=0$. 
\end{theorem}
\proof
Our argument follows \cite{Bol}. Note that the analytic continuation of $Y(\lambda)$ and $\lambda^E$
around $\lambda=0$ are respectively $Y(\lambda)G$ and $\lambda^E
G$. Therefore, the holomorphic branch of 
\ben
U(\lambda):=Y(\lambda)\lambda^{-E}\lambda^{-K}
\een
defined in the sector $S\subset D-\{0\}$ extends analytically to the
entire punctured disc $D-\{0\}$. Using Proposition \ref{reg-sing} we
get that $U(\lambda)$ has at most a finite order pole at $\lambda=0$. 

Let us prove that $U(\lambda)$ is holomorphic at $\lambda=0$. Let us
denote by $r={\rm max}_{1\leq j\leq s}
\operatorname{Re}(\rho_j)$. Since $r<1$ we can find a real number
$\epsilon>0$, such that $r+2\epsilon <1$. We claim that
$\lim_{\lambda\to 0} U(\lambda) \lambda^{r+2\epsilon}=0$.  This
clearly implies that $U(\lambda)$ does not have a pole at
$\lambda=0$. To prove that the limit is 0 we write
\ben
U(\lambda) \lambda^{r+2\epsilon}= Y(\lambda)\lambda^{-K+\epsilon} \, 
\exp\Big( (r-\lambda^K E \lambda^{-K}) \, \log \lambda\Big)\, \lambda^\epsilon.
\een
Note that the first two factors on the RHS give a matrix obtained from 
$Y(\lambda)$ by multiplying each column in $Y_{i,\ell}$ by
$\lambda^{-k_i^\ell +\epsilon}$. Since the Levelt evaluation of
every column in $Y_{i,\ell}$ is at least $k_i^\ell$ we get that the
limit of $Y(\lambda)\lambda^{-K+\epsilon} $ is 0. 
Since $K$ and $E$ have the same block-diagonal structure we get that
3rd and the 4th factor give a matrix which is also block-diagonal and the $i$-th
block is 
\beq\label{factor-34}
\lambda^{\epsilon+r-\rho_i} e^{\lambda^{K_i} N_i\lambda^{-K_i} \log \lambda},
\eeq
where $N_i$ is strictly upper triangular. Since $K_i$ is diagonal with
decreasing diagonal entries the matrix $\lambda^{K_i}
N_i\lambda^{-K_i} $ is holomorphic at $\lambda=0$. Therefore the limit
of \eqref{factor-34} is 0. 

\medskip
It remains only to prove that $U(0)$ is invertible. Substituting
$Y(\lambda)=U(\lambda) \lambda^K\lambda^E$ in the differential
equation we get 
\ben
\lambda U'(\lambda) + U(\lambda) L(\lambda)= B_0(\lambda) U(\lambda),
\een
where $L(\lambda)= K+ \lambda^K E \lambda^{-K}$. As we discussed above
the matrix $\lambda^K E\lambda^{-K}$ is holomorphic at $\lambda=0$. In
fact $L(0)$ is block-diagonal and the $i$th block is 
\ben
\begin{bmatrix}
(k_i^1+\rho_i) I_i^1 + N_i^{11} & 0 & \cdots & 0 \\
0 &  (k_i^2+\rho_i) I_i^2 + N_i^{22} & \cdots & 0 \\
\vdots & \vdots & \ddots &  0 \\
0 & 0& \cdots & (k_i^{m_i}+\rho_i) I_i^{m_i} + N_i^{m_im_i}
\end{bmatrix}.
\een
Since $U(0)L(0)=B_0(0)U(0)$, we get that $L(0)$ is a linear operator
in $\operatorname{Ker}(U(0))$. If we assume that $U(0)$ is not
invertible, then $L(0)$ has a non-zero eigenvector $c\in
\operatorname{Ker}(U(0))$. Let us denote by
$y_c(\lambda)=Y(\lambda)c$. 

Let us split the vector-column $c$ in the following way
\ben
c=\begin{bmatrix}
c_1 \\
\vdots \\
c_s
\end{bmatrix},\quad
c_i =
\begin{bmatrix}
c_{i,1} \\
\vdots \\
c_{i,m_i}
\end{bmatrix},
\een
where the length of the subcolumn $c_{i,\ell}$ is the same as the
dimension of $X_i^\ell/X_i^{\ell-1}$. Since $L(0)$ is block-diagonal
and upper triangular, we get that there exists a unique pair
$(i,\ell)$ for which $c_{i,\ell}\neq 0$ and that the eigenvalue of $c$
is $\rho_i+k_i^\ell$. Note that 
\ben
y_c(\lambda) = Y_{i,\ell} c_{i,\ell}
\een
is a linear combination of elements in $X_i^\ell$ that project to a
basis in $X_i^\ell/X_i^{\ell-1}$. Therefore
$\varphi(y_c)=k_i^\ell$. 

On the other hand, let us denote by $R$ the diagonal part of $E$ and
write $E=R+N$. Note that $[R,N]=0$. Therefore
\ben
y_c(\lambda)=U(\lambda) \lambda^K\lambda^N\lambda^{-K} c\ 
\lambda^{\rho_i+k_i^\ell},
\een
where we used that $K+R$ is the diagonal part of $L(0)$.
Furthermore, using that $\lambda^K N\lambda^{-K} = L(\lambda)-K-R$ is
a holomorphic nilpotent matrix we get
\ben
U(\lambda) \lambda^K\lambda^N\lambda^{-K} =
U(\lambda) e^{(L(\lambda)-K-R)\log \lambda} .
\een
Expanding near $\lambda=0$ we get 
\ben
U(0)\Big(1 +\sum_{k=1}^m \frac{1}{k!}(L(0)-K-R)^k(\log
\lambda)^k\Big)+ 
O(\lambda (\log \lambda)^m),
\een
where $m$ is an integer such that $N^m=0$. However
$(L(0)-K-R)c=U(0)c=0$, so we get that $\varphi(y_c)\geq 1+k_i^\ell$
-- contradiction. 
\qed

\section{Vector bundles on $\mathbb{P}^1$}

Let us assume that $E\to \mathbb{P}^1\times \widetilde{\Pi}$ is a holomorphic vector bundle of rank $p$, where 
\ben
\widetilde{\Pi} = 
\{u=(u_1,\dots,u_N)\in \CC^N\ |\ |u_i-u_i^\circ|<\widetilde{\delta}_i,\ 1\leq i\leq N\},
\een
is the polydisc with center $u^\circ:=(u_1^\circ,\dots,u_N^\circ)$ and polyradius $\widetilde{\delta}=(\widetilde{\delta}_1,\dots,\widetilde{\delta}_N)$. 
The main goal of this lecture is to prove the existence of Birkhoff factorization for the transition matrix of $E$. We follow Bolibruch \cite{Bol}.  

\subsection{Transition function}

We will be interested in transition functions of $E$ of the following type. Let us fix a point $b\in \CC\subset \mathbb{P}^1$, real numbers $0<r<R$, and a polydisc 
\ben
\Pi=\{ u\in \CC^N\ |\ |u_i-u_i^\circ|<\delta_i,\ 1\leq i\leq N\},
\een
where $0<\delta_i<\widetilde{\delta}_i$ for all $i$. The discs
\ben
D_b=\{\lambda \in \CC\ | \ |\lambda-b|<R\},\quad
D_\infty = \{\lambda\in \mathbb{P}^1 |\ |\lambda-b|>r\}
\een
give an open cover of $\mathbb{P}^1$. The open subsets $D_\nu\times \Pi$, $\nu=b,\infty$, are Stein and contractible, so according to the Grauert--Oka principle $E|_{D_\nu\times \Pi}$ is trivial. Let us define raw vectors
\ben
e_\nu=(e_{\nu,1},\dots,e_{\nu,p}),\quad e_{\nu,i}\in 
\Gamma(D_\nu\times \Pi,E),
\een
such that $\{e_{\nu,i}\}_{i=1}^p$ is a trivializing frame for $E|_{D_\nu\times \Pi}$. On the intersection the two frames are related by a holomorphic invertible matrix
\ben
e_\infty(\lambda,u) = e_b(\lambda,u) M(\lambda,u),\quad (\lambda,u)\in D_{b\infty}\times \Pi,
\een
where $D_{b\infty}=D_b\cap D_\infty$ and 
\ben
M:D_{b\infty}\times \Pi\to \operatorname{GL}(\CC^p)
\een
is a holomorphic map. Choosing different trivialization frames $\widetilde{e}_b=e_b U$ and $\widetilde{e}_\infty=e_\infty W$, where 
\ben
U:D_b\to \operatorname{GL}(\CC^p)
\quad \mbox{and}\quad
W:D_\infty \to \operatorname{GL}(\CC^p)
\een
are holomorphic maps, yields a new transition matrix $\widetilde{e}_\infty= \widetilde{e}_b \widetilde{M}$, where
\ben
\widetilde{M} (\lambda,u)= U(\lambda,u)^{-1} M(\lambda,u) W(\lambda,u),\quad (\lambda,u)\in D_{b\infty}\times \Pi.
\een
Our main goal can be stated as follows. We would like to prove that after decreasing $\Pi$ if necessary and removing an analytic hypersurface from $\Pi$ we can always arrange that 
\ben
\widetilde{M} = \operatorname{diag}((\lambda-b)^{k_1},\dots, (\lambda-b)^{k_p}),
\een
where $k_1\geq \cdots \geq k_p$ is a decreasing sequence of integers. 

\subsection{GAGA reduction}

\begin{definition}\label{def:rat}
We say that a map $M: D_{b\infty}\times \Pi \to \mathfrak{gl}(\CC^p)$ is $\Pi$-rational if the entries of $M(\lambda,u)$ are quotients of polynomials in $\O(\Pi)[\lambda]$, where $\O(\Pi)$ is the ring of holomorphic functions on $\Pi$. 
\end{definition}
We would like to reduce the general analytic problem to an algebraic one. More precisely we would like to prove the following proposition
\begin{proposition}\label{pr:gaga}
Decreasing the size of $\Pi$ if necessary, we can find a transition matrix 
\ben
M:D_{b\infty}\times \Pi\to \operatorname{GL}(\CC^p)
\een
such that 
\begin{enumerate}
\item[(i)] $M$ is $\Pi$-rational.
\item[(ii)] The zeroes of $\operatorname{det}(M(\lambda,u))$ and the poles of $M(\lambda,u)$ for $(\lambda,u)\in \mathbb{P}^1\times \Pi$ are independent of $u$. 
\end{enumerate}
\end{proposition}

Let us introduce the following notation. If $\Pi$ is an open polydisc, then we denote by $\overline{\Pi}$ the corresponding closed polydisc. If $X\subset \mathbb{P}^1\times \overline{\Pi}$ is an open subset, then we define
\ben
H(X):=\{\phi: \overline{X} \to \mathfrak{gl}(\CC^p)\ |\ \phi \mbox{ is continuous in $\overline{X}$ and holomorphic in $X$}\},
\een
and 
\ben
H^0(X):=\{\phi\in H(X)\ |\ \phi(x) \mbox{ is invertible for all $x\in X$}\}.
\een
Recall that $H(X)$ is a Banach algebra with norm 
\ben
|\!|A|\!|=\operatorname{sup}_{(\lambda,u)\in \overline{X}} |A(\lambda,u)|,
\een
where $|\ |:\mathfrak{gl}(\CC^p)\to \RR_{\geq 0}$ is the matrix norm
\ben
|A|:=\operatorname{sup}_{v\in \mathbb{C}^p-\{0\}} 
|Av|/|v|,
\een 
where for $w\in \CC^p$ we denote by $|w|=(|w_1|^2+\cdots +|w_p|^2)^{1/2}$ the Euclidean norm of $w$. 
\begin{lemma}\label{BF-1}
There exists an $\epsilon>0$, depending on $r,$ and $R$ such that 
if $B\in H(D_{b\infty}\times \Pi)$ has norm $|\!|B|\!|<\epsilon$, then $1+B
\in H^0(D_{b\infty}\times \Pi)$  and we have a factorization
\ben
1+B=UW,\quad U\in H^0(D_b\times \Pi),\quad W\in H^0(D_\infty\times \Pi).
\een
\end{lemma}
\proof
The Laurent series expansion gives a decomposition
\ben
H(D_{b\infty}\times \Pi) = H(D_{b}\times \Pi)\bigoplus H(D_{\infty}\times
\Pi)(\lambda-b)^{-1},
\quad 
B= B_++B_-.
\een
Let $\operatorname{pr}_\pm$ be the corresponding projection maps $B\mapsto B_\pm$. We have 
\ben
\operatorname{pr}_+(B)(\lambda,u) = 
\frac{1}{2\pi\sqrt{-1}} 
\int_{|\zeta-b|=R} \frac{B(\zeta,u)d\zeta}{\zeta-\lambda} 
\een 
and 
\ben
\operatorname{pr}_-(B)(\lambda,u) = -
\frac{1}{2\pi\sqrt{-1}} 
\int_{|\zeta-b|=r} \frac{B(\zeta,u)d\zeta}{\zeta-\lambda} .
\een 
It is easy to check that $|\!|\operatorname{pr}_\pm( B)|\!| \leq C \, |\!|B|\!|$ for
some constant $C$ that depends on $r$ and $R$. Using these estimates
and choosing $\epsilon$ sufficiently small ($\epsilon <1/C$ works)  we can prove that the
series
\ben
w=\sum_{n=1}^\infty (-\operatorname{pr}_-\circ B)^n 1= -B_-+(BB_-)_--(B(BB_-)_-)_-+\cdots
\een
is convergent. Note that 
$
(Bw)_-+w+B_-=0
$
therefore,  $(1+B)(1+w) = 1+u$, with $u=B_++(Bw)_+$. Decreasing
$\epsilon$ if necessary ($\epsilon<1/(4C)$ works) we can arrange that $1+u$ and $1+w$ are
invertible, so the lemma follows with $U=1+u$ and $W=(1+w)^{-1}$.
\qed

{\em Proof of Proposition \ref{pr:gaga}.}
Let us fix positive numbers $0<r'<r''<r<R<R''<R'$ and polydiscs
$\Pi\subset \Pi''\subset \Pi' \subset \widetilde{\Pi}$  with center at
$u^\circ$. Here $\Pi'$ is chosen arbitrary, while the sizes of $\Pi''$
and $\Pi$ will be specified later on.  

Let us pick an arbitrary transition matrix
\ben
M': D_{b\infty}'\times \Pi' \to \operatorname{GL}(\CC^p),\quad 
D_{b\infty}' :=\{ r'<|\lambda-b|<R'\}.
\een
Note that $M'\in H^0(D_{b\infty}''\times \Pi'')$ where $D''_{b\infty}=\{ r''<|\lambda-b|<R''\}.$ 
The Laurent series expansion of $M'(\lambda,u_0)^{-1}$ at $\lambda=b$
is uniformly convergent for $r''\leq |\lambda-b|\leq R''$, while
$M'(\lambda,u)^{-1}$ is uniformly continuous for $(\lambda,u)\in
\overline{D}''_{b\infty}\times \overline{\Pi}''$. Therefore by
truncating the Laurent series expansion of $M'(\lambda,u_0)^{-1}$
appropriately and choosing $\Pi''$ sufficiently small, we can find a
Laurent polynomial $P\in \mathfrak{gl}(\CC^p)[(\lambda-b)^{\pm 1}]$,
s.t., 
$|\!|M'P-1|\!|_{r'',R'',\Pi''}<\epsilon$, where the norm is in the
space $H(D''_{b\infty}\times \Pi'')$. Recalling Lemma \ref{BF-1}, we
find $U_1\in H^0(D''_b\times\Pi'')$ and $W_1\in
H^0(D''_\infty\times\Pi'')$, s.t., $M'P  = U_1W_1$, i.e.,  
\ben
M'=U_1W_1P^{-1}.
\een 
Similarly, we can find $Q\in \mathfrak{gl}(\CC^p)[(\lambda-b)^{\pm 1}]$, s.t.,
$Q W_1P^{-1}  = U_2 W_2$ with $U_2\in H^0(D_b\times \Pi)$ and $W_2\in
H^0(D_\infty\times \Pi)$ for some sufficiently small polydisc
$\Pi\subset \Pi''$. Therefore we get
\ben
M'=U_1 Q^{-1} U_2 W_2.
\een
We claim that the matrix $M:=Q^{-1} U_2$ has the required
properties. Condition (ii) is easy to verify. Let us prove that
$U_2$ is $\Pi$-rational. We have  
\ben
U_2=QW_1P^{-1}W_2^{-1}.
\een 
Let $g(\lambda)\in \CC[\lambda]$
be a common denominator for the entries of $Q$ and $P^{-1}$. The
matrix $g^2U_2= (gQ)W_1(gP^{-1})W_2^{-1}$ is holomorphic for all
$\lambda\in \CC$, because $U_2$ is holomorphic in $D_b\times \Pi$
while $W_1$ and $W_2$ are holomorphic in $D_\infty\times
\Pi$. Moreover, since $W_1$ and $W_2$ are holomorphic at
$\lambda=\infty$, the matrix $g^2U_2$ has at most a   
finite order pole at $\lambda=\infty$, so it must be polynomial in
$\lambda$. \qed

\subsection{Existence and uniqueness of Birkhoff factorization}

Let $\Pi$ be a polydisc with center $u^\circ$ and $\Theta_0\subset
\Pi$ be an analytic hypersurface with finitely many irreducible
components. Since $\Pi$ is Stein and contractible, there exists a holomorphic function
$f_0\in \O(\Pi)$ such that $\Theta_0$ is the zero locus of
$f_0$. Suppose that we have a transition matrix
\ben
M:D_{b\infty}\times (\Pi-\Theta_0) \to \operatorname{GL}(\CC^p),
\een 
such that 
\begin{enumerate}
\item[(i)] $M$ is $\Pi$-rational.
\item[(ii)]  The zeroes of
$\operatorname{det}(M(\lambda,u))$ and the poles of $M(\lambda,u)$ for
$(\lambda,u)\in D_b\times (\Pi-\Theta_0)$ are independent of
$u$. 
\end{enumerate}
According to the previous section such a transition matrix exists
provided we choose $\Pi$ sufficiently small and
$\Theta_0=\emptyset$. 

Note that condition (i) implies that the points $(\lambda,u)\in D_b\times
\Pi$ for which $M(\lambda,u)$ is not holomorphic form an analytic
hypersurface $Z_\infty(M)$. More precisely $Z_\infty(M)$ is the union
of all irreducible hypersurfaces $V\subset 
D_b\times \Pi$ such that there exists an entry $m=g/f$ ($g,f\in
\O(D_b\times \Pi)$) of $M$ for which 
$
\operatorname{ord}_V (f)>\operatorname{ord}_V (g),
$ 
where $\operatorname{ord}_V (h)$ denotes the order of vanishing of the
holomorphic function $h\in \O(D_b\times \Pi)$ along $V$. 
\begin{lemma}\label{lemma:polar-div}
Every irreducible component of $Z_\infty(M)$ has either the form
$\{b'\}\times \Pi$ for some $b'\in D_b-D_{b\infty}$ or $D_b\times
\Theta'_0$, where $\Theta'_0$ is an irreducible component of
$\Theta_0$. 
\end{lemma}
\proof
Let $V$ be an irreducible component of $Z_\infty(M)$. Condition (ii) implies that 
\ben
V\cap D_b\times (\Pi-\Theta_0) = \bigcup_{i=1}^s \{b_i\}\times (\Pi-\Theta_0) ,
\een
for some $b_i \in D_b$. Since $M(\lambda,u)$ is holomorphic and
invertible for $(\lambda,u)\in D_{b\infty}\times (\Pi-\Theta_0) $ we
have $b_i\in D_b-D_{b\infty}$ and 
\ben
V\subset \Big(\bigcup_{i=1}^s \{b_i\}\times \Pi\Big)\bigcup D_b\times \Theta_0.
\een
The RHS of the above inclusion relation is an analytic hypersurface, so $V$ must be one of its
irreducible components. 
\qed

\begin{proposition}\label{prop:bf}
a) There exists an analytic hypersurface $\Theta\subset \Pi$ that
contains $\Theta_0$ and has finitely many irreducible components, such
that 
\ben
M(\lambda,u)=U(\lambda,u) (\lambda-b)^K W(\lambda,u),
\een
where
\begin{enumerate}
\item[(i)] $U$ and $W$ are $\Pi$-rational.
\item[(ii)] $U(\lambda,u)$ (resp. $W(\lambda,u)$) is holomorphic and
  invertible for all $(\lambda,u)\in D_b\times (\Pi-\Theta)$
  (resp. $D_\infty\times (\Pi-\Theta)$).
\item[(iii)] $K=\operatorname{diag}(k_1,\dots,k_p)$, where $k_1\geq
  \cdots \geq k_p$ are integers. 
 \end{enumerate}

b) If 
\ben
M(\lambda,u)=U_i(\lambda,u) (\lambda-b)^{K^{(i)} } W_i(\lambda,u),\quad i=1,2,
\een
are two factorizations satisfying the conditions in part a), then $K^{(1)}=K^{(2)}.$
\end{proposition}
\proof
a) We split the proof into two cases. 

{\em Case 1:} If $\operatorname{det}(M(\lambda,u))\neq 0$ for all
$(\lambda,u)\in (D_b-\{b\})\times (\Pi-\Theta_0)$. 
We may assume that $M(\lambda,u)=L(\lambda,u)\, (\lambda-b)^K$, where 
$L(\lambda,u)$ is holomorphic for $(\lambda,u)\in D_b\times
(\Pi-\Theta_0)$ and
$K=\operatorname{Diag}(k_1,\dots,k_p)$, where $k_1\geq k_2\geq \cdots
\geq k_p$ are integers. This can be always achieved by first
multiplying $M$ from the right by matrices of the type
$(\lambda-b_0)^{K_0} (\lambda-b)^{-K_0}$, so that we clear all the poles of
$M(\lambda,u)$ from $D_b\times (\Pi-\Theta_0)$, and finally multiply by a
constant permutation matrix  to arrange that the entries of $K$ are
decreasing. Moreover, according to Lemma \ref{lemma:polar-div} 
there exists an integer $n$, such that $L(\lambda,u)f_0(u)^n$ is
holomorphic for all $(\lambda,u)\in D_b\times \Pi$.  
 
The Taylor's series expansion of $L$ has the form
\ben
L(\lambda,u) = L_0(u)+ L_1(u) (\lambda-b)+L_2(u)(\lambda-b)^2+\cdots.
\een
Let us denote by $m_i(u)$, $1\leq i\leq p$, the columns of the matrix
$L_0(u)$. We may assume that $m_1\neq 0$, otherwise we can factor
out $(\lambda-b)$ from the first column of $L(\lambda,u)$ and increase
$k_1$ by one. We can also assume that $\operatorname{det}(L_0(u))=0$,
otherwise the matrix $L(\lambda,u)$ is invertible for all
$(\lambda,u)\in D_b\times (\Pi-\Theta)$, where $\Theta\subset \Pi$ is
the union of $\Theta_0$ and the zero locus of
$\operatorname{det}(L_0(u))$, and this is already a Birkhoff 
factorization, so there is nothing to prove.

Let us denote by $i$ the maximal index, s.t., some $i\times i$ minor
of $L_0(u)$ contained in the first $i$-columns is not identically
0 for $u\in \Pi-\Theta_0$. If there are several such minors, then we
choose one of them, write it in the form $g(u)/f_0(u)^n$ for some
$g\in \O(\Pi)$ and
denote by $\Theta\subset \Pi$ the analytic hypersurface defined by the
zero locus of $g(u)f_0(u)$. There are functions $s_1(u),\dots,s_i(u)$,
holomorphic for $u\in \Pi-\Theta$ and meromorphic along $\Theta$, s.t., 
\ben
m_{i+1}(u)=s_1(u)m_1(u)+\cdots + s_i(u)m_i(u).
\een
Let 
\ben
W(\lambda,u)=1-\sum_{a=1}^i s_a(u) (\lambda-b)^{-k_a+k_{i+1}} E_{a,i+1},
\een
where $E_{a,i+1}$ is the matrix with only one non-zero entry, which is
equal to 1 and it is in row $a$ and column $i+1$. Note that $k_a\geq
k_{i+1}$, so $W(\lambda,u)$ is holomorphic and invertible for
$(\lambda,u)\in D_\infty \times (\Pi-\Theta)$ and meromorphic along
$D_\infty \times \Theta$. It is easy to 
check that $M(\lambda,u)W(\lambda,u) =
\widetilde{L}(\lambda,u)(\lambda-b)^{\widetilde{K}}$, where
$\widetilde{K}=\operatorname{Diag}(\widetilde{k}_1,\dots,\widetilde{k}_p)$
satisfies $\widetilde{k}_j=k_j$ for $j\neq i$ and
$\widetilde{k}_{i+1}> k_{i+1}$. Multiplying if necessary $W$ by a
constant permutation matrix from the right we can arrange that
$\widetilde{k}_1\geq \cdots\geq \widetilde{k}_p$. Note that 
\ben
\operatorname{det}(L(\lambda,u)) =
\operatorname{det}(\widetilde{L}(\lambda,u))\,
(\lambda-b)^{\sum_{i=1}^p (\widetilde{k}_i-k_i)},
\een
so the order of vanishing of $\operatorname{det}(L(\lambda,u)) $ at
$\lambda=b$ decreases strictly. Repeating the above procedure 
finitely many times we will eventually get a matrix $L(\lambda,u)$,
s.t., $\operatorname{det}(L(\lambda,u))\neq 0$ at $\lambda=b$, which
as explained above would give a Birkhoff factorization.

{\em Case 2:} general case. 
Just like in Case 1, multiplying $M(\lambda,u)$ from the right by an
appropriate holomorphic invertible matrix defined for all
$(\lambda,u)\in D_\infty\times \Pi$ and by $(\lambda-b)^m{\rm Id}$ with $m\gg
0$, we may reduce the proof to the case when 
$M(\lambda,u)$ is holomorphic for all  $(\lambda,u)\in D_b\times
(\Pi-\Theta_0)$ and meromorphic along $D_b\times \Theta_0$. 
We argue by induction on the number of zeroes of $\det(M(\lambda,u))$
in $D_b\times (\Pi-\Theta_0)$. If there are no zeroes, then $M(\lambda,u)$ is
holomorphic and invertible for all $(\lambda,u)\in D_b\times
(\Pi-\Theta_0)$ and there is nothing to prove.  

Let $b_1\in D_b$ be a zero. Let us choose a small disc 
$D_1=\{|\lambda-b_1| < R_1\}$ with center $b_1$,
s.t., $D_1\subset D_b$ and $D_1$ does not contain other zeroes of
$\det(M(\lambda,u))$. Let 
us recall Case 1 for $M$ and the covering of $\PP^1$ given by the discs
$D_1$ and $D_1^{\infty}:=\{|\lambda-b_1| > r_1\}$, where
$0<r_1<R_1$. We get a Birkhoff factorization 
\ben
M(\lambda,u) =
M_1(\lambda,u)(\lambda-b_1)^{K_1} W_1(\lambda,u),
\een
 where 
\begin{enumerate}
\item[(i)] $M_1$ and $W_1$ are $\Pi$-rational.
\item[(ii)] $M_1(\lambda,u)$ (resp. $W_1(\lambda,u)$) is holomorphic
  and invertible for all $(\lambda,u) \in D_1\times (\Pi-\Theta_1))$
  (resp. $D_1^{\infty}\times(\Pi-\Theta_1)$) for some analytic
  hypersurface $\Theta_1\subset \Pi$ with $\Theta_0\subset \Theta_1$. 
\item[(iii)] $K_1$ is a diagonal matrix with decreasing integer
  entries. 
\end{enumerate}
Note that  
\beq\label{u0}
M_1(\lambda,u)=M(\lambda,u) W_1(\lambda,u)^{-1} (\lambda-b_1)^{-K_1}
\eeq
is holomorphic for $(\lambda,u)\in D_{b}\times (\Pi-\Theta_1)$ and invertible for $(\lambda,u)\in D_{b\infty}\times (\Pi-\Theta_1)$. The
zeroes of $\det(M_1(\lambda,u))$ for $(\lambda,u)\in D_b\times
(\Pi-\Theta_1)$ are first of all in $(D_b-D_1)\times (\Pi-\Theta_1)$
and then by expecting the RHS of \eqref{u0}, we get that the they are
contained in the set of zeroes of $\det(M(\lambda,u))$. Note that if
$\lambda=b_1$ is the only zero of $\det(M(\lambda,u))$ for
$(\lambda,u)\in D_b\times (\Pi-\Theta_1)$, then we are done, because
$M_1(\lambda,u)$ will be holomorphic and invertible  for
$(\lambda,u)\in D_{b}\times (\Pi-\Theta_1)$. Otherwise, let 
$b_2\in D_b$ be a 2nd zero of  $\det(M(\lambda,u))$ and let $m>0$ be
an integer such that the diagonal entries of $K_1$ are greater than $-m$. We get that the number
of zeroes of $\det\Big(M_1(\lambda,u)(\lambda-b_2)^{K_1+m}\Big)$ in
$D_b\times (\Pi-\Theta_1)$ is at least 1 less than the number of zeroes of
$\det(M(\lambda,u))$. Using the inductive assumption we get a Birkhoff
factorization 
\ben
M_1(\lambda,u) (\lambda-b_2)^{K_1+m}= U(\lambda,u) (\lambda-b)^{K} W'(\lambda,u),
\een
where 
\begin{enumerate}
\item[(i)] $U$ and $W'$ are $\Pi$-rational.
\item[(ii)] $U(\lambda,u)$ (resp. $W'(\lambda,u)$) is holomorphic
  and invertible for all $(\lambda,u) \in D_b\times (\Pi-\Theta))$
  (resp. $D_\infty\times(\Pi-\Theta)$) for some analytic
  hypersurface $\Theta\subset \Pi$ with $\Theta_1\subset \Theta$. 
\item[(iii)] $K$ is a diagonal matrix with decreasing integer
  entries. 
\end{enumerate}
Therefore,
\ben
M(\lambda,u)  = 
U(\lambda,u) (\lambda-b)^{K-m} W'(\lambda,u) 
\Big(\frac{\lambda-b_1}{\lambda-b_2}\Big)^{K_1}
\Big(\frac{\lambda-b}{\lambda-b_2}\Big)^{m}
W_1(\lambda,u)
\een
which provides a Birkhoff factorization for all $u\in \Pi-\Theta$.

b) 
Put $K^{(i)}=\operatorname{Diag}(k^{(i)}_1,\dots,k^{(i)}_p)$. We argue
by induction on $i$ that $k_i^{(1)}=k_i^{(2)}$ for all $i$. Assume that 
$k_a^{(1)}= k_a^{(2)}$ for $a=1,2,\dots,i-1$ and $k_i^{(1)} > k_i^{(2)}$.
Comparing the two Birkhoff factorization, we get 
\ben
(U_2^{-1}U_1)_{a\ell} = (\lambda-b)^{k_a^{(2)}-k_\ell^{(2)} } (W_2 W_1^{-1} )_{a\ell} ,
\een
where $A_{a\ell}$ denotes the $(a,\ell)$-entry of the matrix $A$. The
LHS is analytic for $\lambda\in D_b$. If $k_a^{(2)}<k_\ell^{(1)}$,
then the RHS is analytic in $D_\infty$ and vanishes for
$\lambda=\infty$, so by Liouiville's theorem both sides must
vanish. We get that 
$(U_2^{-1}U_1)_{a\ell}=0$ for $1\leq \ell \leq i$ and $a\geq i$,
because according to our assumptions
\ben
k_a^{(2)}\leq k_i^{(2)} < k_i^{(1)} \leq k_\ell^{(1)}.
\een
The first $i$-columns of $U_2^{-1}U_1$ have non-zero entries only in
the first $(i-1)$ places, therefore they must be linearly
dependent. This however contradicts the fact that $U_2^{-1}U_1$ is
invertible for $(\lambda,u)\in D_b\times V$. Similarly, the assumption
$k_i^{(1)} < k_i^{(2)}$ would contradict the invertibility of
$W_2W_1^{-1}$, so $k_i^{(1)}=k_i^{(2)}$. \qed

\section{Painleve property for the Schlesinger equations}
In this lecture we prove two theorems of Malgrange, which will be used
later on to prove the Painleve property for Frobenius manifolds. Our
arguments follow closely Bolibruch \cite{Bol} for Theorem
\ref{thm:Malgrange-1} and Malgrange \cite{Mal} for Theorem
\ref{thm:Malgrange-2}.

\subsection{Vector bundles on $\mathbb{P}^1$}
Using the results of Lecture 2 we will prove the following theorem of
Malgrange (see \cite{Mal}).   
\begin{proposition}[Malgrange]\label{thm:Malgrange-1}
Suppose that $T$ is a connected smooth analytic variety and $E$ is a vector bundle on
$\PP^1\times T$, s.t., $E_{\PP^1\times \{t_0\}}$ and $E|_{\{b_0\}\times
  T}$ are trivial for some $(b_0,t_0)\in \PP^1\times T$.  Then 

a) The subset
\ben
\Theta=\{t\in T\ :\ E_{\PP^1\times \{t\}} \mbox{ is not trivial } \}
\een
is either empty or it is an analytic hypersurface of $T$.

b) $E|_{\PP^1\times (T-\Theta)}$ is trivial and meromorphic 
along $\PP^1\times \Theta$.  
\end{proposition}

Let us clarify the meaning of being meromorphic in Proposition
\ref{thm:Malgrange-1}. It means that  we can find a trivializing frame
$\{e_i\}_{i=1}^p$ for $E|_{\PP^1\times (T-\Theta)}$, s.t., if
$\{e^U_i\}_{i=1}^p$ is a local frame for $E$ in a neighborhood $U$ of
some point on $\PP^1\times \Theta$, then the transition function
between the two frames is a $p\times p$ matrix whose entries are
meromorphic functions on $U$ with poles along $U\cap (\PP^1\times
\Theta)$. 

\proof
a)
We argue by induction on the dimension of $T$. If $T$ is
0-dimensional, then there is nothing to prove. 
Let us define the set 
\ben
N=\{t\in T\ :\ E|_{\PP^1\times \{t\}}\ \mbox{ is trivial}\}.
\een
\begin{claim}\label{claim:open}
If $t'\in N$, then there exists an open neighborhood $V$ of $t'$ in $T$ such that $E|_{\PP^1\times V}$ is trivial. 
\end{claim}
\proof
Let $V$ be an open polydisc neighborhood of $t'$. We
can find trivializations of $E|_{D_\nu\times V}$, s.t., the transition
function $M(\lambda,t')=1$. Indeed, using that $E_{\PP^1\times \{t'\}}
$ is trivial we get that $M(\lambda,t')= U'(\lambda)W'(\lambda)$,
where $U'(\lambda)$ (resp. $W'(\lambda)$) is holomorphic and
invertible for $\lambda \in D_b$ (resp. $D_\infty$). Changing the
trivialization frames of $E|_{D_b\times V}$ and $E|_{D_\infty\times
  V}$ via $U'$ and $W'$ we can transform the the transition matrix
into $U'(\lambda)^{-1} M(\lambda,t) W'(\lambda)^{-1}$, which turns
into 1 at $t=t'$.   

Let us assume now that the transition matrix is such that $M(\lambda,t')=1$.  Decreasing $V$ if necessary, we can make 
$M(\lambda,t)$ sufficiently close to $M(\lambda,t')$. Recalling Lemma
\ref{BF-1}, we get a Birkhoff factorization
$M(\lambda,t)=U(\lambda,t)W(\lambda,t)$, which implies that
$E|_{\PP^1\times V}$ is trivial.  
\qed

The above claim shows that $N$ is an open subset. 
\begin{claim}\label{claim:trivial}
The vector bundle $E|_{\PP^1\times N}$ is trivial.
\end{claim}
\proof
Let $\Sigma$ be the set of open subsets $V\subset T$ such that
$E|_{\PP^1\times V}$ is trivial. By definition $t_0\in N$, so
according to Claim \ref{claim:open} the set $\Sigma$ is
non-empty. Using the inclusion of open subsets we can define a partial
ordering on $\Sigma$. Clearly every increasing chain $V_1\subset
V_2\subset \cdots$ in $\Sigma$ is bounded by $\cup_i V_i\in
\Sigma$. Therefore, recalling the Zorn's lemma the set $\Sigma$ has a
maximal element, say $V$. If $V\neq N$, then let $t'\in N$ be a
boundary point of $V$. According to Claim \ref{claim:open} we can find
an open neighborhood $V'\in \Sigma$ that contains $t'$. Let $e'$ and
$e$  be raw vectors whose entries give trivializations of respectively
$E|_{\PP^1\times V'}$ and $E|_{\PP^1\times V }$. Then $e'=e U$, where
$U:\PP^1\times (V'\cap V) \to \operatorname{GL}(\CC^p)$ is a
transition matrix. Since the entries of $U(\lambda,u)$ are holomorphic
for all $\lambda\in \PP^1$, by Liouiville's theorem they must be
constants independent of $\lambda$, i.e., $U(\lambda,u)=U(b_0,u)$. On
the other hand by definition $E|_{\{b_0\}\times T}$ is a trivial
bundle, so we can factorize   $U(b_0,u)=A(u) A'(u)^{-1}$. Therefore $e
A(u)=e'A'(u)$ for $u\in V\cap V'$, so we get that $E|_{\PP^1\times
  (V\cup V')}$ is trivial. Since $V$ is maximal we get $V'\subseteq
V$, which however contradicts the fact that $t'\in V'$ is a boundary
point of $V$.  
\qed

If $N=T$, then we are done. Let us assume that $N\neq T$. We have to show that
$T-N$ is an analytic subvariety of codimension 1. Let $u^0\in
T$ be a boundary point of $N$ and $\Pi$ be a polydisc with  center
$u^0$.  Let $M(\lambda,u)$ be the transition matrix for some
trivializations $E|_{D_\nu\times \Pi}$, $\nu=b,\infty$. According to
Proposition \ref{pr:gaga} we may assume that $M$ is $\Pi$-rational. Decreasing $\Pi$ if
necessary, we get a Birkhoff factorization (see Proposition
\ref{prop:bf}, part a)) $M(\lambda,u) = U(\lambda,u) (\lambda-b)^K
W(\lambda,u)$, where $U(\lambda,u)$ (resp. $W(\lambda,u)$) is
holomorphic and invertible for $(\lambda,u)\in D_b\times (\Pi-\Theta)$
(resp. $D_\infty \times (\Pi-\Theta)$).  On the other hand, if $V\subset
(\Pi-\Theta)\cap N\subset N$ is an open subset, then $E_{\PP^1\times V}$ is
trivial. Therefore,  the transition function 
$M(\lambda,u)=U'(\lambda,u)W'(\lambda,u).$ Comparing the two Birkhoff
factorizations of $M$ and recalling Proposition \ref{prop:bf}, part b) we get that
$K=0$, which implies that $E|_{\PP^1\times (\Pi-\Theta)}$ is
trivial, i.e., $\Pi-N\subset \Theta$. 

If $\Pi-N\neq \Theta$, then $N\cap \Theta\neq \emptyset$. Using the
inductive assumption we get that  $E_{\PP^1\times \Theta}$ is trivial
on the complement of some analytic hypersurface $\Theta_0\subset
\Theta$ (note that $\Theta_0$ contains the singular locus of
$\Theta$) and therefore $\Pi-\Theta_0\subset N$. Using the Hartogue's
extension theorem and that $\Theta_0$
is of complex codimension 2, it is easy to prove that $E|_{\PP^1\times
  \Pi}$ is also trivial. Indeed, since the vector bundle $E$ is trivial on $D_\nu\times \Pi$,
$\nu=b,\infty$, we can choose frames $\{e_{b,i}\}_{i=1}^p$ and
$\{e_{\infty,i}\}_{i=1}^p$ . Let $M:D_{b\infty}\times \Pi\to \operatorname{GL}(\CC^p)$ be
the transition matrix, i.e., $e_\infty = e_b\, M$. Since $E$ is
trivial on $\PP^1\times(\Pi-\Theta_0)$, we can choose a frame of $E$
on $\PP^1\times(\Pi-\Theta_0)$. Therefore, we have a Birkhoff
factorization 
\ben
M(\lambda,u) = U(\lambda,u)W(\lambda,u),
\een
where $U(\lambda,u)$ (resp. $W(\lambda,u)$) is holomorphic and
invertible in $D_b\times (\Pi-\Theta_0)$ (resp. $D_\infty \times (\Pi-\Theta_0)$).
On the other hand, since $\Theta_0$ has complex codimension 2 in $\Pi$,
the Hartogues extension theorem implies that $U(\lambda,u)$
(resp. $W(\lambda,u)$) extends analytically for all $(\lambda,u)\in D_b\times
\Pi$ (resp. $D_\infty\times \Pi$). Moreover, the zero locus of
$\operatorname{det}(U(\lambda,u))$ in $D_b\times \Pi$ is contained in
the codimension 2 analytic subset $D_b\times \Theta_0$. However, the
zero locus of an analytic function is  either empty or codimension
1. Therefore, $U(\lambda,u)$ is invertible for $(\lambda,u)\in
D_b\times \Pi$. Similar argument implies that $W(\lambda,u)$ is
invertible for $(\lambda,u)\in D_\infty\times \Pi$. Hence the trivialization of $E$
over $\PP^1\times (\Pi-\Theta_0)$ extends analytically across
$\PP^1\times \Theta_0$, i.e., $\Pi\subset N$. This however contradicts
the fact that the center of $\Pi$ is a boundary 
point of $N$. Therefore, $\Pi-N=\Theta$ is an analytic hypersurface.

b) 
Let $e=(e_1,\dots,e_p)$, $e_i\in \Gamma(\PP^1\times (T-\Theta),E)$ be
a trivializing frame. Every other trivializing frame has the from $e
C$, where $C:T-\Theta \to \operatorname{GL}(\CC^p)$. On the other hand,
since $E_{\{b_0\}\times T}$ is trivial we get that we can always
choose $C$ in such a way that $eC$ extends to a trivializing frame of
$E_{\{b_0\}\times T}$.  Therefore there exists a frame $e$ such that
$e|_{\{b_0\}\times (T-\Theta)}$ extends to a trivializing frame of
$E|_{\{b_0\}\times T}$. We claim that such a frame $e$ is meromorphic. 

To prove this, let us pick a point $u^\circ\in \Theta$, a polydisc
$\Pi$ with center $u^\circ$, and trivializations $e^\Pi_b$ and
$e^\Pi_\infty$ of respectively $E|_{D_b\times \Pi}$ and
$E|_{D_\infty\times \Pi}$ such that the transition function
$M:D_{b\infty}\to \operatorname{GL}(\CC^p)$ is $\Pi$-rational and we
have a Birkhoff factorization 
\ben
M(\lambda,u)=U(\lambda,u) (\lambda-b)^K W(\lambda,u),
\een
such that $U$ (resp. $W$) is $\Pi$-rational, holomorphic, and invertible for $(\lambda,u)\in
D_b\times (\Pi-\Theta')$ (resp. $ D_\infty\times (\Pi-\Theta')$) where
$\Theta'\subset \Pi$ is an analytic hypersurface. Such
choices are possible according to Propositions \ref{pr:gaga} and
\ref{prop:bf} provided  we choose $\Pi$ sufficiently small. As we have
proved in a), $K=0$ and $\Theta'=\Pi-N=\Theta$. Let
$\widetilde{e}^{\, \Pi}_b=e^\Pi_bU$ and
$\widetilde{e}^{\, \Pi}_\infty=e^\Pi_\infty W^{-1}$ be trivializing frames of respectively
$E|_{D_b\times (\Pi-\Theta)}$ and $E|_{D_\infty\times
  (\Pi-\Theta)}$. Note that these frames agree on the intersection
$D_{b\infty}\times (\Pi-\Theta)$, so we get a trivializing frame
$\widetilde{e}^{\, \Pi}$ of $E|_{\PP^1\times (\Pi-\Theta)}$. Therefore,
  there exists a transition matrix $C:\Pi-\Theta \to
  \operatorname{GL}(\CC^p)$, such that $\widetilde{e}^{\, \Pi} = e C$. Let
  $f\in \O(\Pi)$ be the function whose zero locus defines the
  hypersurface $\Theta\cap \Pi$. Since both
  $U(\lambda,u)$ and $W(\lambda,u)^{-1}$ are meromorphic along respectively 
  $D_b\times \Theta$ and $D_\infty\times \Theta$, there exists an
  integer $n>0$ such that $f(u)^n U(\lambda,u)$ (resp. $f(u)^n
  W(\lambda,u)^{-1}$) extends analytically  for all $(\lambda,u)\in
  D_b\times \Pi$ (resp. $D_\infty\times \Pi$). Therefore
  $\widetilde{e}^{\, \Pi} f^n= e C f^n$ is a raw vector whose entries are
  holomorphic section of $E$ on  $\PP^1\times \Pi$. Restricting to
  $\{b_0\}\times \Pi$ and recalling that $e|_{\{b_0\}\times \Pi}$ is a
  trivializing frame, we get that $H(u):=C(u) f(u)^n$ is holomorphic
  for all $u\in \Pi$, i.e., $C$ is meromorphic along
  $\Theta$. Therefore 
\ben
e^{\Pi}_b = e C(u) U(\lambda,u)^{-1},\quad 
e^{\Pi}_\infty = e C(u) W(\lambda,u)
\een
and we get that the transition matrices $C(u)U(\lambda,u)^{-1}$ and
$C(u)W(\lambda,u)$ are meromorphic along respectively $D_b\times
\Theta$ and $D_\infty\times \Theta$. 
\qed

\subsection{The Schlesinger equations}\label{sec:schl_def}

Let $\nabla^\circ$ be a Fuchsian connection on the trivial vector bundle $\PP^1\times \CC^p$. Written in coordinates
\ben
\nabla^\circ = d-A^\circ(\lambda) d\lambda,
\een
where 
\ben
A^\circ(\lambda) = \frac{A^\circ_1}{\lambda-u^\circ_1} +\cdots + \frac{A^\circ_N}{\lambda-u^\circ_N},
\een
where $A_i^\circ$ are $p\times p$ matrices and $u^\circ_i$ are the finite
poles of $\nabla^\circ$. Let us also assume
that  $\sum_{i=1}^N A_i^\circ\neq 0$, so that the connection has a
Fuchsian singularity at $\lambda=\infty$. 

The Schlesinger equations are the following system of differential
equations
\begin{align}
\notag
&\frac{\partial A_i}{\partial u_j}  = \frac{[A_j,A_i]}{u_j-u_i}, \quad 1\leq
i\neq j\leq N,\\
\notag
& \sum_{j=1}^N \frac{\partial A_i}{\partial u_j}  = 0 , \quad 1\leq i\leq N ,\\
\notag
& A_i(u^\circ)=A_i^\circ ,\quad 1\leq i\leq N ,
\end{align}
where $u^\circ=(u_1^\circ,\dots,u_N^\circ)$. Here 
\ben
A_i(u_1,\dots,u_N)\in \mathfrak{gl}(\CC^p),\quad 1\leq i\leq N,
\een
is a set of matrix-valued functions that should be thought as
deformations of the coefficients of the Fuchsian connection
$\nabla^\circ$. It is easy to check that the Schlesinger equations are
compatible (integrable). Therefore the solution exists for all
$u=(u_1,\dots,u_N)$ sufficiently close to $u^\circ$.

The main goal of this lecture is to prove that the local solutions
extend to global meromorphic functions. More precisely, let 
\ben
Z_N=\{u\in (\PP^1)^{N+1}\ :\ u_i\neq u_j \mbox{ for $i\neq j$  and $u_{N+1}=\infty$}\} 
\een
be the configuration space of $N$ points in $\CC$. Every point $u\in
Z_N$ corresponds to a punctured sphere 
\ben
\PP^1-\{u_1,\dots,u_{N+1}\} = \CC-\{u_1,\dots,u_{N}\}.
\een
Let us denote by $T$ the universal cover of $Z_N$. The point
$u^\circ\in Z_N$ will be fixed as a base point and we identify $T$ as
the set of pairs $(u,[c])$ such that $u\in Z_N$ and $[c]$ is the
homotopy class of a path $c$ in $Z_N$ from $u^\circ$ to $u$. A small
neighborhood of $u^\circ$ in $Z_N$ has a natural lift to a small
neighborhood of $t^\circ:=(u^\circ,[1])\in T$, where $[1]$ is the trivial path
from $u^\circ$ to $u^\circ$. In particular, every solution of the
Schlesinger equations (with the specified initial condition) is
defined and analytic in a neighborhood of $t^\circ\in T$. The main
goal of this lecture is to prove the following theorem due to
Malgrange \cite{Mal}. 
\begin{theorem}[Malgrange]\label{thm:Malgrange-2}
If $\{A_i(u)\}_{i=1}^N$ is a solution to the Schlesinger equations,
then each $A_i$ extends to a meromorphic function on $T$.
\end{theorem}
\subsection{Malgrange's vector bundle E}\label{sec:Mal-bun}
The Fuchsian connection $\nabla^\circ$ determines a monodromy
representation 
\ben
\mu: \pi_1(\CC-\{u^\circ_1,\dots,u^\circ_{N}\},b^\circ)\to \operatorname{GL}(\CC^p),
\een
where $b^\circ$ is a reference point. The representation is defined
as follows. Let $Y^\circ(\lambda)$ be the fundamental solution of
$\nabla^\circ$ defined in a neighborhood of $b^\circ$ such that
$Y^\circ(b^\circ)=1$. If $\gamma$ is a closed path in
$\CC-\{u^\circ_1,\dots,u^\circ_{N}\}$ then the analytic continuation
of $Y^\circ(\lambda)$ along $\gamma$ has the form
$Y^\circ(\lambda)\mu(\gamma)$. 

Let us denote by $D_i\subset \PP^1\times T$, $1\leq i\leq N+1$, the
hypersurface consisting of points $(\lambda,u,[c])$ such that
$\lambda=u_i$ ($u_{N+1}:=\infty$). Let $\C:=\PP^1\times
T-\cup_{i=1}^{N+1} D_i$. The projection map 
\ben
\pi:\C\to T,\quad (\lambda,u,[c])\mapsto (u,[c]) 
\een
is a smooth fibration with fiber diffeomorphic to
$\pi^{-1}(t^\circ)=\CC-\{u^\circ_1,\dots,u^\circ_{N}\}$. Since $\pi_k(T)=\pi_k(Z_N)=0$
for $k>1$ and $\pi_1(T)=\{1\}$, we get that $T$ is a contractible
space. Using the long exact sequence of homotopy groups we get that
the natural inclusion 
\ben
(\CC-\{u^\circ_1,\dots,u^\circ_{N}\}, b^\circ) \to (\C,(b^\circ,t^\circ))
\een  
induces an isomorphism between the fundamental groups. Therefore the
monodromy representation $\mu^\circ$ of $\nabla^\circ$ induces a
representation 
\beq\label{mon-repr}
\mu: \pi_1 (\C,(b^\circ,t^\circ))\to \operatorname{GL}(\CC^p).
\eeq
There exists a unique vector bundle $E\to \C$ of rank $p$ equipped with a flat
connection $\nabla$ such that the monodromy representation of $\nabla$
is equivalent to the given representation \eqref{mon-repr}. We will refer to
$E\to \C$ as the {\em Malgrange's vector bundle}.  The equivalence
between the monodromy representation of $\nabla$ and \eqref{mon-repr}
means that there exists a raw vector $f^\circ=(f_1^\circ,\dots,f_p^\circ)$ whose entries form a basis of
the fiber $E_{b^\circ,t^\circ}$ such that the parallel transport with respect
to $\nabla$ along a closed loop $\gamma$ based at $(b^\circ,t^\circ)$ transforms
$f^\circ$ into $f^\circ\, \mu(\gamma)$. 

For the reader's convenience let us recall the construction of $E$. We choose a covering of
$\C$ by open balls $\{B_i\}_{i\in \C}$ that have contractible connected
intersections. This can be achieved by choosing a Riemannian metric on
$\C$ and letting $B_i$ be the ball with center $i\in \C$ of radius
$r_i$, where $r_i$ is the injectivity radius of $\C$ at the point
$i$. It is known that if $x', x''\in B_i$, then there exists a unique
geodesic in $\C$ from $x'$ to $x''$ whose length is the distance
between $x'$ and $x''$. Moreover, such a geodesic is entirely in $B_i$.   
If $B_i\cap B_j\neq\emptyset$, then we choose a
smooth path $\gamma_{ij}$ in $B_i\cup B_j$ 
between the centers of $B_i$ and $B_j$. Let us also fix $B_0$ to be
the ball with center the base point $(b^\circ,t^\circ)$. Let us also
fix a path $\gamma_i$ from $B_0$ to $B_i$ consisting of paths $\gamma_{ab}$.
Then we define $E|_{B_i}:=B_i\times \CC^p$ and let
$e^i=(e^i_1,\dots,e^i_p)$ be the trivializing frame corresponding to
the standard basis of $\CC^p$. On the overlaps $B_i\cap B_j\neq
\emptyset$ the bundles are glued via  
\ben
e^j=e^i \, g_{ij},\quad g_{ij} = \mu(\gamma_i^{-1} \circ \gamma_{ji} \circ \gamma_j),
\een  
where $\mu$ is the given monodromy representation \eqref{mon-repr}. Since $g_{ij}$ are
constants, the standard flat  
connections given by the de Rham differential on $B_i\times \CC^p$
glue together, so the bundle $E$ is naturally equipped with a flat
connection.

\subsection{Extension of $E$}\label{sec:ext}
Recall that $Y^{\circ}(\lambda)$ is the fundamental solution of
$\nabla^\circ$ defined in a neighborhood of a fixed reference point
$\lambda=b^\circ$. $Y^\circ(\lambda)$ is uniquely determined by
requiring that it satisfies the initial condition
$Y^\circ(b^\circ)=1$. For every singular point $u_i^\circ$ ($1\leq i\leq N$) of
$\nabla^\circ$ let us fix a sector with vertex at $u_i^\circ$ of the
following form
\ben
\{ \lambda\in \CC\ :\ 0<|\lambda-u_i^\circ|<R_i^\circ,\ 
-\epsilon < \operatorname{Arg}(\lambda-u_i^\circ)<\epsilon\},
\een
where $R^\circ_i$ is sufficiently small so that the disc with center
$u_i^\circ$ and radius $R_i^\circ$ does not contain other singular points
$u_j^\circ$ and $0<\epsilon<2\pi$. Let us fix a path $\gamma^\circ_i$ ($1\leq i\leq
N+1$) in $\CC-\{u_1^\circ,\dots,u_N^\circ\}$  from
$b^\circ$ to a point $u_i^\circ+\lambda_i^\circ$ in the above
sector, e.g., $\lambda_i^\circ:= R_i^\circ/2$. Let us 
extend analytically $Y^\circ(\lambda)$ along $\gamma^\circ_i$. We get
an analytic  solution of $\nabla^\circ$ defined in the above sector.
Finally, let us choose an invertible 
matrix $S_i^\circ\in \operatorname{GL}(\CC^p)$, such that
$Y^\circ(\lambda)S_i^\circ$ is a weak Levelt solution for the Fuchsian
singularity of $\nabla^\circ$ at $\lambda=u_i^\circ$. We have 
\beq\label{w-Levelt}
Y^\circ(\lambda) S_i^\circ = U^\circ_i(\lambda) (\lambda-u^\circ_i)^{K_i} (\lambda-u^\circ_i)^{E_i},
\eeq
where the matrix 
\ben
E_i=\diag(E_i^1,\dots,E_i^{p_i})
\een 
is block diagonal with each block corresponding to an eigenvalue of
$E_i$, the block $E_i^j=\rho_i^j I + N_i^j$, where $N_i^j$ is an
upper-triangular nilpotent matrix and the eigenvalue $\rho_i^j$ satisfies
\ben
0\leq \operatorname{Re}(\rho_i^j)<1,
\een  
$K_i=\operatorname{diag}(K_i^1,\dots,K_i^{p_i})$ has the same
block diagonal structure as $E_i$ with each block $K_i^j$ being a
diagonal matrix with decreasing integer entries, and $U^\circ_i(\lambda)$ is
holomorphically invertible in a neighborhood of 
$\lambda=u^\circ_i$.

It is convenient to extend our notation for the singular points
of $\nabla^\circ$ in order to include also the singularity at
$\lambda=u^\circ_{N+1}=\infty$. The above statements
remain the same except that we have to replace everywhere
$\lambda-u_i^\circ$ with $\lambda^{-1}$. In particular, the
fundamental solution takes the forms 
\beq\label{wlb-infty}
Y^\circ(\lambda) S^\circ_{N+1} = U^\circ_{N+1}(\lambda) \lambda^{-K_{N+1}} \lambda^{-E_{N+1}}.
\eeq

\medskip

The vector bundle $E$ can be extended across the divisors $D_i$
($1\leq i\leq N+1$) as follows. Let us take a tubular neighborhood 
\ben
T_i=\{ (\lambda,u,[c])\ :\ |\lambda-u_i|<R_i(u)\}
\subset \PP^1\times T,
\een
where $R_i:Z_N \to \RR_{>0}$ is a smooth function satisfying 
\ben
R_i(u)< |u_j-u_i|,\quad \mbox{for all } \quad 1\leq i\neq j \leq N,
\een
and 
\ben
R_{N+1}(u)>|u_j|,\quad \mbox{for all}\quad 1\leq j\leq N.
\een
Using parallel transport with respect to the flat connection $\nabla$
we construct a multivalued flat frame $f=(f_1,\dots,f_p)$ of $E$ whose
value at a point $(\lambda,t)\in \C$  
\ben
f(\lambda,t)=(f_1(\lambda,t),\dots,f_p(\lambda,t)), \quad
f_i(\lambda,t) \in E_{\lambda,t}
\een
depends on the choice of a reference path in $\C$ from $(b^\circ,t^\circ)$ to
$(\lambda,t)$: the component $f_i(\lambda,t)$ is obtained from $f_i^\circ\in
E_{b^\circ,t^\circ}$ via a parallel transport along the reference path. Let us
trivialize $E|_{T_i-D_i}$ via the frame
\beq\label{Di-frame}
f(\lambda,t) S_i^\circ(\lambda-u_i)^{-E_i} (\lambda-u_i)^{-K_i},\quad
(\lambda,t)\in T_i-D_i,
\eeq
where $t=(u,[c])\in T$ and the path specifying the value of
$f(\lambda,t)$ is chosen as follows. We identify
$\CC-\{u_1^\circ,\dots,u_N^\circ\}$
with the fiber  $\C_{t^\circ}:=\pi^{-1}(t^\circ)$. Note that the path
$\gamma_i^\circ\subset \CC_{t^\circ}$ connects the reference point
$(b^\circ,t^\circ)$ with the point $(u_i^\circ + R_i^\circ/2,
t^\circ)\in T_i$ (provided we define $R_i^\circ:=R_i(u^\circ)$). The
path that we would like to select consists of two pieces the path
$\gamma_i^\circ$ and any path in $T_i-D_i$ connecting the end point of
$\gamma_i^\circ$ and $(\lambda,t)$.  The analytic continuation of
$f(\lambda,t^\circ)$ and $Y^\circ(\lambda)$ along a closed loop
around $\lambda=u_i^\circ$ are respectively  $f(\lambda,t^\circ)M_i$
and $Y^\circ(\lambda) M_i$, where $M_i$ is such that $M_iS_i^\circ= S_i^\circ e^{2\pi\sqrt{-1}
  E_i}$. The monodromy of $f(\lambda,t)S^\circ_i$ around $D_i$
cancels out the monodromy of  $(\lambda-u_i)^{-E_i}
(\lambda-u_i)^{-K_i}$ around $D_i$. Hence the frame \eqref{Di-frame}
provides a holomorphic trivialization of $E|_{T_i-D_i}$. 
We extend $E$ across $D_i$ in the obvious way: on the overlap of
$T_i$ and $T_i-D_i$ we identify the standard
frame of $T_i\times \mathbb{C}^p$ with the frame \eqref{Di-frame} of
$E|_{U_i-D_i}$.   

\subsection{Proof of Theorem \ref{thm:Malgrange-2}}

We are going to construct a multivalued analytic function
$Y(\lambda,t)$ with values in $\operatorname{GL}(\CC^p)$  defined for
all $(\lambda,t)\in \C$ such that 
\begin{enumerate}
\item[(1)] $Y(\lambda,t^\circ)=Y^\circ(\lambda)$. 
\item[(2)] The 1-form
$\omega:=dY(\lambda,t) Y(\lambda,t)^{-1}$ is a meromorphic 1-form on
$\PP^1\times T$ of the form
\ben
\sum_{i=1}^N \frac{A_i(t)}{\lambda-u_i} \, (d\lambda-du_i),
\een
where $A_i$ is a $\mathfrak{gl}(\CC^p)$-valued meromorphic function on
$T$ and $u_i:T\to \CC$ is the $i$th component of the projection map
$T\to Z_N$. 
\end{enumerate}
If we manage to do this then Theorem \ref{thm:Malgrange-2} follows
immediately. Indeed, the 1st condition implies that
$A_i(t^\circ)=A_i^\circ$. While the fact that $A_i(t)$ satisfy the
Schlesinger equations follows from the fact that $\omega$ is a 
1-form satisfying
\ben
d\omega +\omega\wedge \omega = d(dY\, Y^{-1}) + dY Y^{-1} \wedge dY Y^{-1} = 0.
\een 

The matrix-valued function $Y(\lambda,t)$ is constructed by comparing
two trivializing frames of $E$. The first one is the multivalued flat
frame 
\ben
f(\lambda,t)=(f_1(\lambda,t),\dots,f_p(\lambda,t)),\quad
f_i(\lambda,t)\in E_{\lambda,t},
\een
defined by the parallel transport with respect to $\nabla$ with
initial value $f(b^\circ,t^\circ):=f^\circ.$ Recall that $f^\circ$ is
the frame of $E_{b^\circ,t^\circ}$  that we fixed so that
the monodromy representation of $\nabla$ coincides with the monodromy
representation \eqref{mon-repr}.

The 2nd frame will be constructed by using Theorem
\ref{thm:Malgrange-1}, which guarantees the existence of a meromorphic
trivialization of $E$. Let us check that the conditions of Theorem
\ref{thm:Malgrange-1} are satisfied. By definition,
$D_{N+1}=\{\infty\}\times T$ and $E|_{D_{N+1}}$ is trivial. 
\begin{claim}
The restriction $E|_{\PP^1\times {t^\circ}}$ is trivial. 
\end{claim} 
\proof
We will prove that $f(\lambda,t^\circ) Y^\circ(\lambda)^{-1}$ is a
trivializing frame. By definition the monodromy of the frame
$f(\lambda,t^\circ)$ and the monodromy of the matrix
$Y^\circ(\lambda)^{-1}$ cancel each other. Therefore the above frame
provides a trivialization of $E|_{\PP^1\times {t^\circ}}$ on
$\CC-\{u_1^\circ,\dots,u_N^\circ\}$. Let us check that the
trivialization extends analytically in a neighborhood of
$\lambda=u_i^\circ$ for all $1\leq i\leq N+1$. Let us assume that
$1\leq i\leq N$. The case $i=N+1$ is the same but one has to use
slightly different notation. By definition the trivializing frame of
$E|_{\PP^1\times {t^\circ}}$ in a neighborhood of $\lambda=u_i^\circ$
is given by 
\ben
f(\lambda,t^\circ) S_i^\circ (\lambda-u_i^\circ)^{-E_i} (\lambda-u_i^\circ)^{-K_i} .
\een
However, recalling the definition of $S_i^\circ$ we get that the above
frame coincides with 
\ben
f(\lambda,t^\circ) Y^\circ(\lambda)^{-1} U_i^\circ(\lambda).
\een
According to Levelt's theorem $U_i^\circ(\lambda)$ is holomorphically
invertible at $\lambda=u_i^\circ$. Therefore the frame
$f(\lambda,t^\circ) Y^\circ(\lambda)^{-1}$ extends holomorphically and
it remains a frame at the point $\lambda=u_i^\circ$. 
\qed

According to Theorem \ref{thm:Malgrange-1}, there exists an analytic
hypersurface $\Theta\subset T$, such that $E|_{\PP^1\times (T-\Theta)}$ is a
trivial vector bundle. Let 
\ben
\widetilde{e}=(\widetilde{e}_1,\dots,\widetilde{e}_p),\quad
\widetilde{e}_i\in \Gamma(\PP^1\times (T-\Theta),E)
\een
be a trivializing frame. We may further  assume that
$\widetilde{e}(\lambda,t^\circ)=f(\lambda,t^\circ)Y^\circ(\lambda)^{-1}$. The
frame that we need in order to define $Y(\lambda,u)$ is slightly
different. The necessary modification is constructed as follows. In
the tubular neighborhood $T_{N+1}$ we have
\ben
f(\lambda,t) S_{N+1}^\circ \lambda^{E_{N+1}} \lambda^{K_{N+1}} =
\widetilde{e}(\lambda,t) \widetilde{U}(\lambda,t),\quad
\forall (\lambda,t)\in T_{N+1}-T_{N+1}\cap (\PP^1\times \Theta),
\een
where $\widetilde{U}(\lambda,t)$ is holomorphic and invertible for all
$(\lambda,t)\in T_{N+1}-T_{N+1}\cap (\PP^1\times \Theta)$ and
meromorphic along $T_{N+1}\cap (\PP^1\times \Theta)$. The Taylor
series expansion at $\lambda=\infty$ yields
\ben
\widetilde{U}(\lambda,t) = \widetilde{U}_0(t)+\widetilde{U}_1(t)
\lambda^{-1} +\widetilde{U}_2(t)
\lambda^{-2}+\cdots,
\een
where $\widetilde{U}_0(t)$ is holomorphic and invertible for all $t\in
T-\Theta$ and meromorphic along $\Theta$. The frame that we need is
\ben
e(\lambda,t)=\widetilde{e}(\lambda,t) 
\widetilde{U}_0(t)^{-1}\widetilde{U}_0(t^\circ).
\een
Note that the above frame is holomorphic for all $t\in T-\Theta$ and
meromorphic along $\Theta$. 

Let us define $Y(\lambda,t)\in \operatorname{GL}(\CC^p)$ as the
transition matrix 
\ben
f(\lambda,t) = e(\lambda,t) Y(\lambda,t),\quad (\lambda,t)\in
\C-(\C \cap(\PP^1\times \Theta)). 
\een
Note that at $t=t^\circ$ we have
$Y(\lambda,t^\circ)=Y^\circ(\lambda)$. Therefore, we need to check
that the 1-form $\omega = dY\, Y^{-1}$ has the required properties. 

To begin with, note that $\omega$ is single valued and
analytic on $\C$. Indeed, the monodromy of $Y(\lambda,t)$ is the same as the
monodromy of $f(\lambda,t)$, i.e., under the analytic continuation
along a closed loop $\gamma$ the value of $Y(\lambda,t)$ changes into
$Y(\lambda,t) \mu(\gamma)$. However, $\mu(\gamma)$ is independent of
$\lambda$ and $t$, so the value of $\omega$ remains the same. Since
being analytic is a local property and locally $Y(\lambda,t)$ is
analytic the same is true for $\omega$. 

Let us analyze the singularities of $\omega$ as a 1-form on
$\PP^1\times T$. The possible singular locus is along the following
divisors
\ben
D_i\ (1\leq i\leq N+1),\quad \PP^1\times \Theta. 
\een  
Let us fix $t\notin \Theta$ and look in a neighborhood of
$\lambda=u_i$ for $1\leq i\leq N$. We have 
\ben
f(\lambda,t)S_i^\circ (\lambda-u_i)^{-E_i}(\lambda-u_i)^{-K_i} =
e(\lambda,u) U_i(\lambda,t),
\een
where $U_i$ is holomorphic and invertible for all $(\lambda,t)\in
T_i-T_i\cap (\PP^1\times T)$ and meromorphic along  $T_i\cap
(\PP^1\times T)$. In particular, the Taylor series expansion at
$\lambda=u_i$ takes the form
\ben
U_i(\lambda,t) = U_{i,0}(t)+U_{i,1}(t) (\lambda-u_i)+\cdots,
\een
where $U_{i,0}(t)$ is holomorphic and invertible for $t\in T-\Theta$
and meromorphic along $\Theta$.  Recalling the definition of
$Y(\lambda,t)$ we get
\ben
Y(\lambda,t) S_i^\circ= U_i(\lambda,t) (\lambda-u_i)^{K_i}(\lambda-u_i)^{E_i},
\een
where the branch of $Y(\lambda,t)$ is determined by an appropriate
reference path (see Section \ref{sec:ext} and the definition of the
frame \eqref{Di-frame}).
Similarly, at $\lambda=\infty$ we get 
\ben
Y(\lambda,t) = 
\widetilde{U}_0(t^\circ) \widetilde{U}_0(t)^{-1} \widetilde{U}(t,\lambda)
\lambda^{-K_{N+1}} \lambda^{-E_{N+1}}. 
\een
Put 
\ben
A_i(t,\lambda):= -(\lambda-u_i) (\partial_{u_i} Y(\lambda,t) )Y(\lambda,t)^{-1}.
\een
If $t\notin \Theta$ is fixed then $A_i$ is an analytic matrix-valued
function on $\CC-\{u_1,\dots,u_N\}$. Near $\lambda=u_j$ with $1\leq
j\neq i\leq N$ we get 
\ben
A_i(t,\lambda) = (u_i-u_j) (\partial_{u_i} U_{j,0}(t)) U_{j,0}(t)^{-1}
+ O(\lambda-u_j),
\een
which is analytic in a neighborhood of $\lambda=u_j$.
Nera $\lambda=u_i$ we get 
\ben
A_i(t,\lambda) & = &  -(\lambda-u_i) (\partial_{u_i} U_{i}(\lambda,t))
U_{i}(\lambda, t)^{-1} + \\ 
&& U_i(\lambda,t) 
\Big(K_i + (\lambda-u_i)^{K_i} E_i (\lambda-u_i)^{-K_i}\Big)
U_i(\lambda,t)^{-1}.
\een
Using the special form of the matrices $E_i$ and $K_i$ we get that the
above expression is analytic at $\lambda=u_i$. Finally at
$\lambda=\infty$ we have
\ben
A_i(t,\lambda) = -(\lambda-u_i) \widetilde{U}_0(t^\circ) 
\partial_{u_i}(\widetilde{U}_0(t)^{-1} \widetilde{U}(t,\lambda))
\widetilde{U}(t,\lambda)^{-1} \widetilde{U}_0(t) \widetilde{U}_0(t^\circ)^{-1},
\een
and this again is analytic at $\lambda=\infty$. According to
Liouiville's theorem $A_i(t,\lambda)$ is independent of
$\lambda$. Setting $\lambda=u_i$ we get that 
\ben
A_i(t):=A_i(t,\lambda) = U_{i,0}(t) C_i U_{i,0}(t)^{-1}, 
\een
where $C_i$ is a constant upper triangular matrix. Moreover, we get
that $A_i$ is meromorphic along $\Theta$.  

Similar argument shows that the matrix
\ben
A(\lambda,t) := (\partial_\lambda Y(\lambda,t) ) Y(\lambda,t)^{-1}
\een
is holomorphic at $\lambda=\infty$ and equal to $0$ at
$\lambda=\infty$. While at $\lambda=u_i$ we have 
\ben
A(\lambda,t) = \frac{A_i(t)}{\lambda-u_i} + \cdots,
\een
where the dots stand for terms analytic at $\lambda=u_i$. This implies
that 
\ben
A(\lambda,t)-\sum_{i=1}^N \frac{A_i(t)}{\lambda-u_i}
\een
is analytic for all $\lambda\in \PP^1$ and vanishing at
$\lambda=\infty$. Recalling again Liouiville's theorem we get that 
\ben
A(\lambda,t)=\sum_{i=1}^N \frac{A_i(t)}{\lambda-u_i}.
\een
Summarizing, we get that 
\ben
\omega=dY\, Y^{-1} = \sum_{i=1}^N \frac{A_i(t)}{\lambda-u_i} \, (d\lambda-du_i),
\een
where $A_i$ are meromorphic functions on $T$. This completes the proof
of Theorem \ref{thm:Malgrange-2}.
\qed

\subsection{Levelt solution with
  parameters}\label{sec:Levelt-def}

The proof of Theorem \ref{thm:Malgrange-2} has the following
interesting corollary. Suppose that we have a Fuchsian connection
$\nabla^\circ$ of the same form as in Section \ref{sec:schl_def}. Let 
$Y^\circ(\lambda)$ be a fundamental solution defined in a neighborhood
of a fixed reference point $b^\circ\in
\CC-\{u_1^\circ,\dots,u_N^\circ\}$. Using the same notation as in
Section \ref{sec:ext}, let us  fix reference paths
$\gamma_i^\circ$ ($1\leq i\leq N+1$) connecting $b^\circ$ with a
neighborhood of $u_i^\circ$  and invertible matrices $S_i\in
\operatorname{GL}(\CC^p)$ ($1\leq i\leq N+1$) such that
$Y^\circ(\lambda)S_i$ is a weak Levelt solution of the form
\eqref{w-Levelt}.   

The isomonodromic deformations of Schlesinger preserve the form of the
Levelt solutions. Namely, let $\{A_i(u)\}_{i=1}^N$ be the solution
to the Schlesinger equations satisfying the initial condition
$A_i(u^\circ)=A_i^\circ$ and defined for all $u$ sufficiently close to
$u^\circ$. Then the system 
\begin{align}
\notag
\partial_\lambda Y(\lambda,u) & =  \Big(\sum_{i=1}^N
\frac{A_i(u)}{\lambda-u_i} \Big) \, Y(\lambda,u) \\
\notag
\partial_{u_i} Y(\lambda,u) & =  -
\frac{A_i(u)}{\lambda-u_i}  \, Y(\lambda,u) ,\quad 1\leq i\leq N,
\end{align}
satisfying the initial condition $Y(\lambda,u^\circ)=Y^\circ(\lambda)$ has a unique
solution. Moreover, for $\lambda$ close to $u_i$ ($1\leq i\leq N$) we have 
\ben
Y(\lambda,u) S_i= U_i(\lambda,u) (\lambda-u_i)^{K_i} (\lambda-u_i)^{E_i}
\een 
and for $\lambda$ close to $\infty$ we have 
\ben
Y(\lambda,u)S_{N+1}=U_{N+1}(\lambda,u) \lambda^{-K_{N+1}} \lambda^{-E_{N+1}},
\een
where the matrices $S_i, K_i$, and $E_i$ ($1\leq i\leq N+1$) are
independent of the deformation parameters $u$.

\section{Tau-function of the Schlesinger equation}

Recall that for a given Fuchsian connection $\nabla^0$ we have Malgrange's vector
bundle $E$ on $\PP^1\times T$. According to Theorem
\ref{thm:Malgrange-1}, $E$ is trivial in the
complement of $\PP^1\times \Theta$, where $\Theta\subset T$ is the subset
of all points $t$ such that $E|_{\PP^1\times \{t\}}$ is trivial. 
The main goal of this lecture is to present a simple algorithm due to
Bolibruch \cite{Bol} that allows us to compute the equation defining
$\Theta$ in terms of the solution of the corresponding Schlesinger
equations.

\subsection{Tau-function}
The notion of {\em tau-function of an isomonodromic deformation} was
introduced in the work of M. Jimbo, T. Miwa, and K. Ueno \cite{JMU}.
The key to the construction of the tau-function in our settings is the following
1-form
\beq\label{tau-form}
\omega= \frac{1}{2}\sum_{i=1}^N\sum_{j:j\neq i}
\frac{\operatorname{Tr}(A_i(u)A_j(u))}{u_i-u_j}\, (du_i-du_j),
\eeq
where $\{A_i(u)\}_{i=1}^N$ is a solution to the Schlesinger equations
satisfying given initial condition $A_i(u^\circ)=A_i^\circ$. According
to Theorem \ref{thm:Malgrange-2} $\omega$ is a meromorphic 1-form on
$T$ with poles along the divisor $\Theta$. We are going to prove the
following lemma.
\begin{lemma}\label{le:theta-loc}
Suppose that  $t^*\in \Theta$ and $\tau_i^*\in
\O_{T,t^*}$ ($1\leq i\leq s$) are the holomorphic germs whose zero
loci define the irreducible components of the germ
of $\Theta$ at $t^*$. Then there are integers $r_i$ ($1\leq i\leq s$)
such that the 1-form 
$
\omega-\sum_{i=1}^s r_i\, \frac{d\tau_i^*}{\tau_i^*}
$
is holomorphic at $t^*$.
\end{lemma}
Following Bolibruch we prove this lemma by giving an algorithm that
produces a set of meromorphic functions whose zero loci contain $\Theta$.
This lemma implies the following theorem.
\begin{theorem}\label{thm:weak-Miwa}
There exists a meromorphic function $\tau$ on $T$ such that
$\omega=d\log \tau.$  
\end{theorem}
\proof
Let us define 
\ben
\tau(t):=\exp\Big(\int_{t^\circ}^t \omega\Big),
\een
where the integral is along a path in $T-\Theta$. Note that this is a
single valued holomorphic function on $T-\Theta$ because according to
Lemma \ref{le:theta-loc} the periods of $\omega$ along closed loops
around $\Theta$ are integer multiples of $2\pi\sqrt{-1}$. 

Let us prove that $\tau$ is meromorphic at $t^*\in \Theta$.  Let
us denote by $\tau_i^*\in \O_{T,t^*}$ ($1\leq i\leq s$) the functions
that define the irreducible components of $\Theta$ at $t^*$. Let us take
$U^*\subset T$ to be an open neighborhood 
of $t^*$ such that $\tau_i^*$ and $\omega-\sum_{i=1}^sr_i \frac{d\tau_i^*}{\tau_i^*}$ can be
represented respectively by  holomorphic functions on $U^*$ and a
holomorphic 1-form on $U^*$. Finally, let us pick a point $b\in
U^*-\Theta$ and let $t\in U^*-\Theta$. Then we have 
\ben
\tau(t) = e^{\int_{t^\circ}^b\omega + \int_b^t
\omega}=\tau(b) e^{\int_b^t\omega}.
\een
On the other hand
\ben
\int_b^t\omega = \sum_{i=1}^s r_i(\log \tau_i^*(t)-\log \tau_i^*(b)) + \int_b^t \Big(\omega
-r_1\frac{d\tau_1^*}{\tau_1^*}-\cdots - r_s\frac{d\tau_s^*}{\tau_s^*}\Big). 
\een
The integral on the RHS extends analytically for all $t\in U^*$, because the
integrand is a holomorphic 1-form in $U^*$. Therefore up to a
holomorphically invertible function in $U^*$ the function
$\tau(t)$ equals $\tau_1^*(t)^{r_1}\cdots \tau_s^*(t)^{r_s} $. This
proves that $\tau$ is meromorphic in $U^*$.  To prove that $\tau$ is
globally meromorphic we have to recall that 
$T$ is a contractible Stein manifold, so every function which is
locally meromorphic must be globally meromorphic. 
\qed

Function $\tau\in \O_T(T-\Theta)$ having the properties in Theorem
\ref{thm:weak-Miwa} is unique up to a non-zero constant factor. Indeed, if
$\tau_1$ and $\tau_2$ are two such functions, then 
\ben
d(\tau_1(t)/\tau_2(t)) = (\tau_1(t)/\tau_2(t))\Big( \frac{d\tau_1(t)}{\tau_1(t)} -
\frac{d\tau_2(t)}{\tau_2(t)}\Big) = 0
\een 
for all $t\in T-\Theta$.  Every function $\tau(t)$ satisfying the
properties in Theorem \ref{thm:weak-Miwa} is called {\em tau-function of
  the isomonodromic deformation}. In fact, it is a theorem due to Miwa \cite{Miw}
(in the case of generic monodromy data) 
and Malgrange \cite{Mal}  (in all cases) that the tau-function is analytic
along $\Theta$. 

\begin{remark}
Bolibruch claimed in his notes \cite{Bol} that his algorithm implies
the analyticity of the tau-function. However, there seems to be a
missing justification. Namely Bolibruch's algorithm produces a sequence of
meromorphic functions, whose product according to the general theory
should be holomorphic. However, using only the
elementary approach pursued in these lectures, we could not justify
the analyticity of the product.   
\end{remark}

\begin{remark}
There is a notion of tau-function more generally for isomonodromic
deformations of connections with irregular singularities. 
Miwa proved the analyticity of the tau-function in full generality. However, in
Miwas's work there is a generality assumption about the monodromy data. Namely,
the monodromy operators are diagonalizable. 
\end{remark}

\subsection{Bolibruch's algorithm}
Let $t^*=(u^*,[c])\in \Theta$ be a generic point, where 
$u^*=(u^*_1,\dots,u^*_N)$. We will be interested only in a small
neighborhood of $t^*$ in $T$, which is isomorphic via the covering map
$T\to Z_N$ to a small neighborhood $\U$ of $u^*$ in $Z_N$. Using this local
bi-holomorphism we will sometimes write $u\in T$ for all $u\in \U$. 
Let us recall also the notation from Section \ref{sec:ext} involving
the following data:
\begin{enumerate}
\item
Fuchsian connection $\nabla^\circ$ on $\PP^1\times \{t^\circ\}$, where
$t^\circ=(u^\circ,[1])$ is a fixed reference point.
\item
We fixed a reference point $b^\circ$ in $\PP^1$, fundamental
solution $Y^\circ$ of $\nabla^\circ$ such that
$Y^\circ(b^\circ)=1$, and a system of paths from $b^\circ$
to each singular point $u_i^\circ$ that allows us to analytically
continue $Y^\circ(\lambda)$ in a neighborhood of each
$\lambda=u_i^\circ$. 
\item 
For each singular point we have chosen a constant matrix $S_i$ such
that $Y^\circ(\lambda)S_i$ is a weak Levelt solution of the type
\eqref{w-Levelt}. 
\end{enumerate}
This is the data necessary to define Malgrange's bundle $E\to
\PP^1\times T$ and hence determines the analytic hypersurface $\Theta$
as well.

The first step in Bolibruch's algorithm is to construct an auxiliary
Fuchsian system on $\PP^1\times \{t^*\}$ that has an extra singular
point. To avoid cumbersome notation let us
assume that $u_i^*\neq 0$ ($1\leq i\leq N$). Then for an extra
singular point we choose $0.$ Let us denote by
$f=(f_1,\dots,f_p) $ the multivalued flat frame of Malgrange's bundle
$E\to \C$ (see Section \ref{sec:Mal-bun}). The frame $f$ provides a
trivialization of $E|_{D_0\times \{t^*\}}$ where $D_0\subset
\PP^1$ is a small neighborhood of $\lambda=0$. According to
Proposition \ref{prop:bf} there exists a trivializing frame $\widetilde{e}$ of
$E|_{(\PP^1-\{0\})\times \{t^*\}}$, and a 
matrix $\widetilde{U}^*(\lambda,t^*)$ holomorphic and invertible for
all  $\lambda\in D_0$, s.t.,  
\ben
\widetilde{e}(\lambda,t^*)= f(\lambda,t^*)\,
\widetilde{U}^*(\lambda,t^*)^{-1} \lambda^{-K},
\quad 
\lambda\in D_0,
\een
where $K=\diag(k_1,\dots,k_p)$ with $k_1\geq \cdots\geq k_p$
integers and we have fixed also a
reference path in $\C\subset \PP^1\times T$ from
$(\lambda^\circ,t^\circ)$, which specifies the value $f(\lambda,t^*)$
. Note that at least one $k_i\neq 0$, otherwise $E_{\PP^1\times
  \{t^*\}}$ would be   
trivial, which
contradicts the definition of $\Theta$. Permuting the entries of the frame $f$ if necessary,
we can arrange that the matrix $\widetilde{U}^*(0,t^*)$ has
non-vanishing principal minors. 
\begin{remark}\label{re:theta-K}
The matrix $\widetilde{U}^*(\lambda,t^*)$ depends analytically on $t^*$ if we
allow $t^*$ to vary along a subset $\Theta_K\subset \Theta$ along
which the vector bundle $E|_{\PP^1\times \{t^*\}}\cong \O(k_1+\cdots
+k_p)$. It is known that $\Theta_K$ is a constructible subset of $\Theta$:
intersection of a closed and an open subsets.  
\end{remark}
\begin{lemma}\label{le:gamma-gauge}
a)
There exists a unique matrix $\Gamma(\lambda,t^*)$ polynomial in
$\lambda^{-1}$ such that  $\Gamma(\lambda,t^*)$ is invertible
for $\lambda\neq 0$ and 
\ben
\Gamma(\lambda,t^*)\, \lambda^{K} \, \widetilde{U}^*(\lambda,t^*) = 
U^*(\lambda,t^*) \lambda^{K},
\een
where $U^*(\lambda,t^*)$ is holomorphic and invertible
for all $\lambda\in D_0$. 

b)
The matrix 
\ben
\lambda^{-K} U^*(\lambda,t^*) \, \lambda^K 
\een
is holomorphic and invertible for all $\lambda\in D_0$. 
\end{lemma}
\proof
To avoid cumbersome notation let us redenote  $V(\lambda,t^*):=\widetilde{U}^*(\lambda,t^*)$ and
assume that $K=\operatorname{diag}(c_1I_1,\dots,c_s I_s)$, where $c_1>\cdots >c_s$ are the
eigenvalues of $K$ and $I_j$ is the identity matrix of size $m_j:=$
the multiplicity of $c_j$ in the sequence $(k_1,\dots,k_p)$. Given a
matrix $A$ of size $p\times p$, then we denote by $A^{lm}$, $1\leq
l,m\leq s$ the block
in position $(l,m)$, where the splitting of a $A$ into blocks is
according to the block-diagonal structure of $K$. We argue
by induction on the number of blocks $s$. Moreover, we are going to
prove that $\Gamma$ is block-lower triangular such that the block
$\Gamma^{lm}$ is a polynomial in $\lambda^{-1}$ of degree
$\leq c_m-c_l$ and the diagonal blocks $\Gamma^{mm}$ are identity matrices. 

For $s=1$, the statements are trivial. For $s>1$ let us write $K=K'+K''$, where 
\ben
K'=\diag((c_1-c_{s-1})I_1,\dots,(c_{s-2}-c_{s-1})I_{s-2},0\cdot
I_{s-1},0\cdot I_s)
\een
and 
\ben
K'' = \diag( c_{s-1}I_1,\dots,c_{s-1} I_{s-2}, c_{s-1}I_{s-1},c_sI_s).
\een
We have
\ben
\lambda^K V(\lambda,t^*)=\lambda^{K'} V'(\lambda,t^*) \lambda^{K''},
\een
where the matrix 
$V'(\lambda,t^*)=\lambda^{K''}V(\lambda,t^*)\lambda^{-K''}$ has the form 
\ben
V' (\lambda,t^*)= 
\begin{bmatrix}
A(\lambda,t^*) \phantom{\lambda^{-m} } & B(\lambda,t^*)\lambda^{m}\\
C(\lambda,t^*) \lambda^{-m}& D(\lambda,t^*)\phantom{\lambda^{m}}
\end{bmatrix},
\een
where $m=c_{s-1}-c_s$, $A$ is a $(p-m_s)\times (p-m_s)$ matrix invertible at
$\lambda=0$, and $B(\lambda,t^*)$, $C(\lambda,t^*)$, and $D(\lambda,t^*)$ are
holomorphic matrices, whose sizes are uniquely determined from the
sizes of $A$ and $V'$.   
There exists a matrix polynomial 
\ben
R(\lambda,t^*) = \sum_{j=0}^m R_j (t^*) \lambda^j,
\een
where $R_j$ is a matrix of size $m_s\times (p-m_s)$, s.t., 
\ben
(C(\lambda,t^*) +R(\lambda,t^*) A(\lambda,t^*)) \lambda^{-m}
\een 
is holomorphic and vanishing at $\lambda=0$. The matrices $R_j(t^*)$
are uniquely determined by
requiring that the coefficients in front of the non-positive powers of
$\lambda$ vanish
\ben
C_j(t^*) +(R_0(t^*) A_j(t^*) +\cdots + R_j(t^*) A_0(t^*)) = 0,\quad
0\leq j\leq m,
\een
where $C_j(t^*)$ and $A_j(t^*)$ are the coefficients in front of
$\lambda^j$ in the Taylor's expansion at $\lambda=0$ of respectively
$C(\lambda,t^*)$ and $A(\lambda,t^*)$.  The assumption about the principal minors of
$V(0,t^*)$ implies that $A_0(t^*)$ is invertible ($\because$ it is a
principal minor of $V(0,t^*)$), so the equations for $R_j(t^*)$
($0\leq j\leq m)$ can be
solved uniquely. Note that the matrix 
\ben
\Gamma_s(\lambda,t^*):= \lambda^{K'}
\begin{bmatrix}
I & 0 \\
R(\lambda,t^*)\lambda^{-m} & I
\end{bmatrix}\, \lambda^{-K'}
\een
is polynomial in $\lambda^{-1}$ and holomorphically invertible for
$\lambda\neq 0$. Moreover, the non-zero off diagonal blocks have the
form $\Gamma_s^{si}$ and they are polynomials in $\lambda^{-1}$ of
degree at most $c_i-c_s$.  
We have 
\ben
\Gamma_s(\lambda,t^*)\, \lambda^K V(\lambda,t^*) = 
\lambda^{K'} V''(\lambda,t^*)\lambda^{K''},
\een
where 
\ben
V''(\lambda,t^*)= 
\begin{bmatrix}
A(\lambda,t^*) & B(\lambda,t^*)\lambda^{m} \\
(C(\lambda,t^*)+R(\lambda,t^*) A(\lambda,t^*))\lambda^{-m} & D(\lambda,t^*)+R(\lambda,t^*) B(\lambda,t^*) 
\end{bmatrix}.
\een
Note that $V''(\lambda,t^*)$ is holomorphic at $\lambda=0$. We claim that $V''(0,t^*)$ has
non-vanishing principal minors. In order to prove this we recall that
an invertible matrix $M$ has non-vanishing principal 
minors if and only if it admits a $LDU$-decomposition, i.e., $M$ can be
written as the product of lower-triangular, diagonal, and
uper-triangular matrices. Put $A_0:=A(0,t^*)$, $B_0:=B(0,t^*)$, $C_0:=C(0,t^*)$,
$D_0:=D(0,t^*)$ and note that $R_0(t^*)=R(0,t^*) = -C_0A_0^{-1}$. We get
\ben
V''(0,t^*)=
\begin{bmatrix}
A_0 & 0 \\
0 & D_0-C_0A_0^{-1} B_0
\end{bmatrix}=
\begin{bmatrix}
I& 0\\
-C_0A_0^{-1} & I
\end{bmatrix}\,
V(0,t^*)\,
\begin{bmatrix}
I & -A_0^{-1}B_0 \\
0 & I
\end{bmatrix},
\een
so the matrix $V''(0,t^*)$ has a $LDU$-decomposition, because according to
our assumptions $V(0,t^*)$ has non-vanishing principal minors, which
implies that $V(0,t^*)$ has a $LDU$-decomposition.

Recalling the inductive assumption we find a matrix
$\Gamma'(\lambda,t^*)$, s.t., $\Gamma'(\lambda,t^*)$ is polynomial in
$\lambda^{-1}$, invertible for $\lambda \neq 0$, and 
\ben
\Gamma'(\lambda,t^*)\, \lambda^{K'} V''(\lambda,t^*) = U^*(\lambda,t^*)\, \lambda^{K'},
\een
where $U^*(\lambda,t^*)$ is a matrix holomorphically invertible
in a neighborhood of $\lambda=0$. We claim that the matrix
$\Gamma(\lambda,t^*):=\Gamma'(\lambda,t^*)\Gamma_s(\lambda,t^*)$ 
satisfies all the required properties. In fact, the only thing left to
check is that the
degree of the block $\Gamma^{lm}$ as a polynomial in $\lambda^{-1}$ is
at most $c_m-c_l$. We have 
\ben
\Gamma^{lm} = \sum_{k=m}^l (\Gamma')^{lk}\Gamma_s^{km}.
\een
If $l\leq s-1$, then $\Gamma_s^{km}\neq 0$ only for $k=m$, so
$\Gamma^{lm}=(\Gamma')^{l m}$. Recalling the inductive assumption for
$\Gamma'$ we get that the degree in $\lambda^{-1}$ does not exceed
\ben
c'_m-c'_l = (c_m-c_{t-1})-(c_l-c_{t-1}) = c_m-c_l.
\een
The case $l=m=s$ is trivial. If $l=s$ and $m=s-1$. Then since the
$(s-1,s-1)$-block of $\Gamma'$  with respect to the block-matrix
structure of $K'$ has the form
\ben
\begin{bmatrix}
(\Gamma')^{s-1,s-1} & (\Gamma')^{s-1,s}  \\
(\Gamma')^{s,s-1} & (\Gamma')^{s,s} 
\end{bmatrix}
\een
and by inductive assumption this should be an identity matrix, we get
that $(\Gamma')^{s,s-1}=0$ and $ (\Gamma')^{s,s}$ is an identity
matrix. Therefore $\Gamma^{lm}=\Gamma_s^{s,s-1}$ so the degree of
this matrix as polynomial in $\lambda^{-1}$ as we proved above is
$\leq c_{s-1}-c_s=c_m-c_l$. Finally, if $l=s$ and $m\leq s-2$, then
$\Gamma_s^{km}\neq 0$ either if $k=m$ or if $k=s$,
i.e., 
\beq\label{gamma_tm}
\Gamma^{sm} = (\Gamma')^{sm} + \Gamma_s^{sm} .
\eeq
Note that the $(s-1,m)$-block of $\Gamma'$ with respect to the
block-matrix structure of $K'$ has the form
\ben
\begin{bmatrix}
(\Gamma')^{s-1,m} \\
(\Gamma')^{s,m}
\end{bmatrix}.
\een 
The inductive assumption about the degree with respect to
$\lambda^{-1}$ implies that the degree of $(\Gamma')^{s,m}$ does not
exceed $c'_m-c'_{s-1} = c'_m=c_m-c_{s-1}< c_m-c_s$. Therefore, both 
terms on the RHS of \eqref{gamma_tm} are polynomials in $\lambda^{-1}$ of degree
$\leq c_m-c_s$, which completes the proof of the inductive assumption
and hence of part a) as well.

For part b), it is enough to check that 
\ben
\lambda^{-K}\Gamma(\lambda,t^*)\lambda^K
\een 
is holomorphically invertible at $\lambda=0$, but this follows
easily if we recall the block-lower triangular structure of $\Gamma$
combined with the degree estimates of the blocks $\Gamma^{lm}$. 
\qed

\medskip

Let us define $e(\lambda,t^*)=\widetilde{e}(\lambda,t^*) \,
\Gamma(\lambda,t^*)^{-1}C(t^*)^{-1}$, where $C(t^*)$ is a constant
invertible matrix. The choice of $C(t^*)$
will be specified bellow. We have the following relation  
\ben
f(\lambda,t^*) = e(\lambda,t^*)\, C(t^*)\,U^*(\lambda,t^*) \, \lambda^K,\quad
\lambda\in D_0. 
\een
Let us denote by $Y^*(\lambda)$ the multivalued function on
$\PP^1-\{0,u_1^*,\dots,u_N^*,\infty\}$ (here $u^*\in Z_N$ is the projection
of $t^*$) with matrix values defined by  
\ben
f(\lambda,t^*)= e(\lambda,t^*) \, Y^*(\lambda).
\een
In particular, if $\lambda$ is close to $0$ we have 
\ben
Y^*(\lambda) = C(t^*)U^*(\lambda) \, \lambda^K.
\een
The local forms of $Y^*(\lambda)$ near the remaining singularities are
\ben
Y^*(\lambda) = 
U_i^*(\lambda) \, (\lambda-u_i^*)^{K_i}\, (\lambda-u_i^*)^{E_i}\, S_i^{-1},
\een
if $\lambda$ is close to $u_i^*$ ($1\leq i\leq N$) and 
\ben
Y^*(\lambda) = U_{N+1}^*(\lambda,t^*) \, \lambda^{-K_{N+1}} \lambda^{-E_{N+1}}\, S_{N+1}^{-1},
\een
if $\lambda$ is close to $u_{N+1}^*:=\infty$. Here the matrices
$U^*_j(\lambda)$ are holomorphically invertible near
$\lambda=u_j^*$ for $1\leq j\leq N+1$. Changing the matrix $C(t^*)$ if
necessary we can arrange that the Taylor's series of
$U^*_{N+1}(\lambda)$ has a constant term $U^*_{N+1}(\infty)=1.$ 
Using the above expansions, it is easy to verify that
\beq\label{Schl:1-form}
\partial_\lambda Y^*(\lambda) =
\Big(
\frac{A_0^*}{\lambda}\, d\lambda + \sum_{i=1}^N
\frac{A_i^*}{\lambda-u_i^*}
\Big)
Y^*(\lambda).
\eeq
Shrinking the neighborhood $\U$ of $t^*$ if necessary we may assume
that the Schlesinger equations 
\ben\label{Schl}
&&
d\widehat{A}_i  = \sum_{\substack{j:j\neq i\\
0\leq j\leq N}}
\frac{[\widehat{A}_j,\widehat{A}_i]}{u_j-u_i} \, (du_j-du_i),
\quad 0\leq i\leq N, \\
&&
\notag
\widehat{A}_i(0,u_1^*,\dots,u_N^*)  =A_i^*
\een
have holomorphic solutions $\widehat{A}_i(u_0,u_1,\dots,u_N)$ defined
for all $u_0$ close to $0$ and for all $(u_1,\dots,u_N)\in \U$. Then we
define 
\ben
A_i(u_1,\dots,u_N):=\widehat{A}_i(0,u_1,\dots,u_N).
\een 
In other words  we have constructed an isomonodromic deformation that
keeps the singular point $0$ fixed.
The system  
\begin{align}
\label{isom-family-1}
\partial_\lambda Y(\lambda,u) & = \Big(\frac{A_0(u)}{\lambda}+\sum_{i=1}^N \frac{A_i(u)}{\lambda-u_i}\Big) Y(\lambda,u) ,\\
\label{isom-family-2}
\partial_{u_i} Y(\lambda,u)  & =-\frac{A_i(u)}{\lambda-u_i}\, Y(\lambda,u) ,\quad 
1\leq i\leq N,
\end{align}
has a unique solution $Y(\lambda,u)$ satisfying the initial condition
$Y(\lambda,u^*)=Y^*(\lambda)$, where
$u^*=(u^*_1,\dots,u^*_N)\in Z_N$ is the projection of $t^*$. Finally,
according to the remark of Section \ref{sec:Levelt-def} the
deformation $Y(\lambda,u)$ has the following 
local expansions  
\beq\label{expansion-infty}
Y(\lambda,u)=U_{N+1}(\lambda,u) \lambda^{-K_{N+1}} \lambda^{-E_{N+1}} S_{N+1}^{-1}
\eeq
for $\lambda$ near $\infty$,
\beq\label{expansion-ui}
Y(\lambda,u) = U_i(\lambda,u) (\lambda-u_i)^{K_i} (\lambda-u_i)^{E_i}
S_i^{-1},\quad 1\leq i\leq N,
\eeq
for $\lambda$ near $u_i$, and 
\beq\label{expansion-0}
Y(\lambda,u) = U(\lambda,u) \lambda^K 
\eeq
for $\lambda$ near $0$.

Let us express the coefficients of the Fuchsian
connection in terms of the coefficients of the Taylor's series
expansion of $U(\lambda,u)$. Substituting 
\ben
Y(\lambda,u) = (U_0(u)+U_1(u)\,\lambda+\cdots ) \, \lambda^K
\een
in the differential equations and
comparing the coefficients in front of $\lambda$ we get the following
relations 
\ben
A_0(u)=U_0(u)KU_0(u)^{-1},
\een
\beq\label{g:der-ui}
u_i\partial_{u_i} (U_0^{-1} U_1) = \frac{U_0^{-1} A_i U_0}{u_i},\quad 1\leq i\leq N,
\eeq
and
\ben
U_0^{-1} U_1 +[U_0^{-1} U_1,K]=-\sum_{i=1}^N  \frac{U_0^{-1} A_i U_0}{u_i}.
\een
More generally, comparing the coefficients in front of $\lambda^{k-1}$ for $k>0$
yields
\begin{align}
\notag
u_i\partial_{u_i} (U_0^{-1} U_k) & =  
(U_0^{-1} A_i U_0)\,\sum_{s=1}^k u_i^{-s} U_0^{-1} U_{k-s} , \quad
1\leq i\leq N, \\
\notag
(k-\operatorname{ad}_K)(U_0^{-1} U_k) & =   -\sum_{i=1}^N (U_0^{-1} A_i U_0)\,\sum_{s=1}^k u_i^{-s} U_0^{-1} U_{k-s} .
\end{align}

\subsection{Gauge transformations}\label{sec:gauge}

Bolibruch has introduced gauge transformations of the following type. 
Let $g(u)$ be an entry of $U_0^{-1}U_1$, and denote its position by
$(\alpha,\beta)$.  Put 
\ben
\Gamma_1(\lambda,u) := I - \frac{U_0(u)E_{\beta\alpha}U_0(u)^{-1} }{g(u)\lambda},
\een
where $E_{\beta\alpha}$ is the matrix with 1 on position
$(\beta,\alpha)$ and $0$ elsewhere. 
The gauge transformation $Y_1(\lambda,u)=\Gamma_1(\lambda,u)\,
Y(\lambda,u)$ yields an isomonodromic system of
differential equations corresponding to a Schlesinger deformation
\begin{align}
\notag
\partial_\lambda Y_1(\lambda,u) & = 
\Big(\frac{A^1_0(u)}{\lambda}+\sum_{i=1}^N
                                  \frac{A^1_i(u)}{\lambda-u_i}\Big)
                                  Y_1(\lambda,u) ,\\ 
\notag
\partial_{u_i} Y_1(\lambda,u) & =
-\frac{A^1_i(u)}{\lambda-u_i}\,
Y_1(\lambda,u) ,\quad  
1\leq i\leq N,
\end{align}
where the matrices $A^1_i$ are given by the following formulas. Put 
\ben
N_{\alpha\beta} = U_0(u)E_{\beta\alpha}U_0(u)^{-1}
\een
and note that this is a nilpotent matrix $N_{\alpha\beta}^2=0$. Then 
\ben
A^1_i(u) = 
\Big(I-\frac{N_{\alpha\beta}}{g(u) u_i}\Big)\, A_i\, 
\Big(I+\frac{N_{\alpha\beta}}{g(u) u_i}\Big),
\quad
1\leq i\leq N,
\een
and 
\ben
A_0^1(u)=A_0+\sum_{i=1}^N (A_i-A_i^1).
\een
The local form of the expansion of $Y_1(\lambda,u)$ at
the singular points $\lambda=u_i$, $1\leq i\leq N+1$ is the same as
for $Y(\lambda,u)$. Let us define $U^1(\lambda,u)$ such that
\ben
Y_1(\lambda,u) =: U^1(\lambda,u)\, \lambda^{K^1},
\een
$K^1=\diag(k^1_1,\dots,k^1_p)=K-E_{\alpha\alpha}+E_{\beta\beta}$, i.e., 
\ben
k^1_i=
\begin{cases}
k_i-1, & \mbox{if } i=\alpha,\\
k_i+1, & \mbox{if } i=\beta, \\
k_i, & \mbox{otherwise.}
\end{cases}
\een
The properties of the matrix $U^1$ are summarized in the following
lemma. 
\begin{lemma}
The local expansion at $\lambda=0$ of $U^1(\lambda,u)$ has the form
\ben
U^1(\lambda,u)=U^1_0(u)+U^1_1(u)\lambda+U^1_2(u)\lambda^2+\cdots,
\een
where $U^1_0(u)$ is holomorphic an invertible for all $u\in \U$ 
such that $g(u)\neq 0$. 
\end{lemma}
\proof
By definition
\ben
U^1(\lambda,u) = \Gamma_1(\lambda,u) U(\lambda,u) 
\lambda^{E_{\alpha\alpha}-E_{\beta\beta}}.
\een
Substituting the expansions at $\lambda=0$ we get a Laurent series that
has a pole of order at most 2. It is easy to see that the coefficients
in front of $\lambda^{-2}$ and $\lambda^{-1}$ are 0. While the
remaining coefficients are as follows. 
\ben
U^1_0(u) & = & U_0 \Big( 
\sum_{i:i\neq \alpha,\beta} E_{ii}+ U_0^{-1} U_1 E_{\beta\beta} +\\
&&
-\frac{1 }{g}\Big(E_{\beta\alpha}+ 
\sum_{i:i\neq \alpha,\beta}
(U_0^{-1}U_1)_{\alpha i} E_{\beta i}
+(U_0^{-1}U_2)_{\alpha\beta} E_{\beta\beta}\Big)
\Big)
\een
and 
\ben
U^1_k(u) & = & U_0 \Big( 
U_0^{-1} U_{k-1} E_{\alpha\alpha}+ 
\sum_{i:i\neq \alpha,\beta} U_0^{-1} U_k E_{ii}+ U_0^{-1} U_{k+1} E_{\beta\beta} 
+ \\
&&
- \frac{1 }{g}\Big((U_0^{-1}U_k)_{\alpha \alpha} E_{\beta\alpha}+ 
\sum_{i:i\neq \alpha,\beta} (U_0^{-1}U_{k+1})_{\alpha i} E_{\beta i}
+(U_0^{-1}U_{k+2})_{\alpha\beta} E_{\beta\beta}\Big)
\Big),
\een
where we have denoted by $(U_0^{-1}U_\ell)_{ab}$ the $(a,b)$-entry of
the matrix $U_0^{-1}U_\ell$. 

We have to prove that the matrix $U_0^1$ is invertible. Note that 
\ben
 U_0^1 = U_0\, g^{E_{\alpha\alpha}-E_{\beta\beta}} \, \widetilde{U}_0^1
\een
where the matrix $\widetilde{U}_0^1$ has the form
\ben
\begin{bmatrix}
1 &  &  &  &  & &  & b_1 & & &\\
& \ddots  &  &  &  & &  & \vdots  & & &\\
&  & 1  &  &  & &  & b_{\alpha-1} & & &\\
&  &  & 0 &  & &  & 1  & & & \\
&  &  &  &1  & &  & b_{\alpha+1} & & &\\
&  &  &  &  &\ddots &  & \vdots & &  &\\
&  &  &  &  & &  1& b_{\beta-1} & & &\\
a_1 &\dots  & a_{\alpha-1}  & -1 &a_{\alpha+1}  &\dots &a_{\beta-1}
&x & a_{\beta+1}&\dots   & a_p\\
 &  &  &  &  & &  & b_{\beta+1} & 1 &  &\\
&  &  &   &  & &  & \vdots & &\ddots  &\\
&  &  &  &  & &  & b_p & &  & 1
\end{bmatrix}.
\een
The entries 
\ben
a_i:=-(U_0^{-1}U_1)_{\alpha i},\quad i\neq \alpha,\beta,
\een
\ben
b_i:= (U_0^{-1}U_1)_{i\beta},\quad i\neq \alpha,\beta,
\een
and $x=-(U_0^{-1}U_2)_{\alpha\beta} + g (U_0^{-1}U_1)_{\beta\beta}.$
The inverse of the matrix $\widetilde{U}_0^1$ is straightforward to
compute. The answer is the following 
\ben
(\widetilde{U}_0^1)^{-1}=
\begin{bmatrix}
1 &  &  &-b_1  &  & &  & & & &\\
& \ddots  &  & \vdots  &  & &  &  & & &\\
&  & 1  & -b_{\alpha-1}  &  & &  & & & &\\
a_1 &\dots  & a_{\alpha-1}  & y &a_{\alpha+1}  &\dots &a_{\beta-1}
&-1 & a_{\beta+1}&\dots   & a_p\\
&  &  & -b_{\alpha+1} &1  & &  & & & &\\
&  &  & \vdots &  &\ddots &  &  & &  &\\
&  &  &-b_{\beta-1}   &  & &  1& & & &\\
&  &  & 1 &  & &  & 0  & & &\\
 &  &  & -b_{\beta+1} &  & &  & & 1 &  &\\
&  &  &  \vdots  &  & &  &  & &\ddots &  \\
&  &  & -b_p  &  & &  & & & &1 
\end{bmatrix},
\een
where 
\ben
y=-x-\sum_{i:i\neq \alpha,\beta} a_i b_i .
\een
The proof of the Lemma is complete. 
\qed

\begin{lemma}\label{le:g-value}
If the position $(\alpha,\beta)$ of $g$ satisfies
$k_\alpha-k_\beta>1$, then $g(t^*)=0$. 
\end{lemma}
\proof
By construction $U(\lambda,u^*)=C(t^*) U^*(\lambda)$. Therefore
$U_0(u^*)^{-1}U_1(u^*)= (U_0^*)^{-1} U_1^*$, where $U_i^*$ is the
coefficient in fron of $\lambda^i$ of the Taylor's series expansion of
$U^*(\lambda)$ at $\lambda=0$. Recalling Lemma \ref{le:gamma-gauge},
part b) we get that 
\ben
(U_\ell^*)_{ab}=0 \quad\mbox{ if $k_a-k_b> \ell$},
\een
where $A_{ab}$ denotes the $(a,b)$-entry of $A$. 
This implies that the non-zero entries of  $(U_0^*)^{-1}$ are in
positions $(a,b)$ such that $k_a-k_b\leq 0$. We have
\ben
g(t^*)=((U_0^*)^{-1} U_1^*)_{\alpha\beta} = \sum_{m=1}^p
((U_0^*)^{-1})_{\alpha m} (U_1^*)_{m\beta}.
\een
The only non-zero terms could be for $m$ such that $k_\alpha-k_m\leq
0$ and $k_m-k_\beta\leq 1$. However such $m$ do not exist otherwise
$k_\alpha-k_\beta\leq 1$.
\qed

Let us determine how the 1-form 
\ben
\omega^*:=\frac{1}{2}\sum_{i=0}^N\sum_{j:j\neq i}
\frac{\operatorname{tr}(A_i(u)A_j(u))}{u_i-u_j} \, (du_i-du_j)
\een
changes under the gauge transformation. Recall that we are working
only with deformations that keep $u_0$ fixed, i.e., $u_0=0.$ 

\begin{lemma}\label{le:form-gauge}
The 1-form 
\ben
\omega^1:= \frac{1}{2}\sum_{i=0}^N\sum_{j:j\neq i}
\frac{\operatorname{tr}(A^1_i(u)A^1_j(u))}{u_i-u_j} \, (du_i-du_j)
\een
satisfies 
\ben
\omega^1=\omega^* + d\log g(u).
\een
\end{lemma}
\proof
The form $\omega^1$ can be written also as 
\ben
\frac{1}{2}\, \sum_{i=1}^N du_i \, \operatorname{Res}_{\lambda=u_i}  
\operatorname{tr}\Big( \sum_{j=0}^N \frac{A^1_j(u)}{\lambda-u_j} \Big)^2
\een
The sum over $j$ is precisely $\partial_\lambda Y_1(\lambda,u)
Y_1(\lambda,u)^{-1}$. Furthermore, recalling the gauge transformation
and using some elementary properties of the trace operation we
get 
\ben
\operatorname{tr}\Big(\partial_\lambda Y_1 Y_1^{-1} \Big)^2 = 
\operatorname{tr}\Big(
(\partial_\lambda \Gamma_1 \Gamma_1^{-1} )^2 +
2 \Gamma_1^{-1}\partial_\lambda\Gamma_1 \, \partial_\lambda Y Y^{-1}
+ (\partial_\lambda Y Y^{-1})^2 \Big).
\een
Substituting the formula for $\Gamma_1=1- (U_0 E_{\beta\alpha}U_0^{-1}
)\, g^{-1}\lambda^{-1}$ we get 
\ben
\frac{1}{2}\operatorname{tr}\Big(\partial_\lambda Y_1 Y_1^{-1} \Big)^2 - 
\frac{1}{2}\operatorname{tr}\Big(\partial_\lambda Y Y^{-1} \Big)^2 = 
\frac{1}{g} \sum_{j=0}^N  
\operatorname{tr}\Big( 
U_0 E_{\beta\alpha} U_0^{-1} \frac{A_j(u)}{\lambda^2 (\lambda-u_j)} \Big).
\een
The residue of the above function at $\lambda=u_i$ is 
\ben
\operatorname{tr}\Big( 
U_0 E_{\beta\alpha} U_0^{-1} \frac{A_i}{u_i^2 } \Big) = 
\frac{(U_0^{-1}  A_i U_0 )_{\alpha\beta} }{u_i^2 } =
\partial_{u_i} g,
\een
where in the last equality we used formula \eqref{g:der-ui}.
\qed

\subsection{Proof of Lemma \ref{le:theta-loc}}

Let us split the matrix $K$ into a block diagonal form
$\diag(c_1I_1,\dots,c_sI_s)$, where $I_j$ is the identity matrix of
size equal to the multiplicity of the number $c_j$ in the sequence
$(k_1,\dots,k_p)$. If $A$ is a $p\times p$ matrix, then we split it
into blocks according to the block-diagonal structure of $K$ and
denote by $A^{lm}$ the block in position $(l,m)$. 
\begin{lemma}
If $c_l-c_m>1$, then at least one entry among the entries of the
blocks $(U_0^{-1}A_iU_0)^{lm}$, $1\leq i\leq N$, is not identically
$0$.   
\end{lemma}
\proof
Assume that this is not true, i.e., $(U_0^{-1}A_iU_0)^{lm}=0$ for all
$i=1,2,\dots,N$ and all $l$ and $m$, s.t., $c_l-c_m>1$. Let us make a
gauge transformation  
\ben
\widetilde{Y}(\lambda,u) = \lambda^{-K} U_0(u)^{-1} Y(\lambda,u).
\een
We get the following differential equation
\ben
\partial_\lambda　\widetilde{Y}(\lambda,u)=
\Big(
\sum_{i=1}^N  \lambda^{-K} \,\frac{U_0^{-1} A_i U_0}{\lambda-u_i} \, \lambda^{K} \Big)
\widetilde{Y}(\lambda,u).
\een
Our assumption implies that the above system is Fuchsian at $\lambda=0$ and that the coefficients $\widetilde{B}_0$ in front of $\lambda^{-1}$ is a nilpotent matrix. This would imply that the monodromy around 
$\lambda=0$ is $e^{2\pi\sqrt{-1} \widetilde{B}_0}=1$ (see Corollary \ref{cor:mon-nilp}), i.e., $\widetilde{B}_0=0$. Therefore, the matrix $\widetilde{Y}(\lambda,u)$ is regular at $\lambda=0$ for all $u$ sufficiently close to $u^*$. Note that $f\, Y(\lambda,u)^{-1} \, U_0(u) $ is a frame for the vector bundle $E|_{\PP^1\times \{u\}}$ on $\PP^1-\{0\}$, while  $f\, \widetilde{Y}(\lambda,u)^{-1} $ is a frame for $E|_{\PP^1\times \{u\}}$ near $\lambda=0$. The relation 
\ben
f\, Y(\lambda,u)^{-1} \, U_0(u)= (f\, \widetilde{Y}(\lambda,u)^{-1}) \, \lambda^K
\een
implies that $E_{\PP^1\times \{u\}}\cong E_{\PP^1\times \{u^*\}}$ is a non-trivial vector bundle for all $u$ in a neighborhood of $u^*$, which contradicts the fact that $\Theta$ is at most a hypersurface.
\qed
 
The above lemma implies that at least one entry in
$(U_0^{-1}U_1)^{lm}$ with $c_l-c_m>1$ is not identically $0$. Let us
choose such an entry $g(u)$ and let $(\alpha,\beta)$ be its
position. By definition $k_\alpha=c_\ell$ and $k_\beta=c_m$ so
$k_\alpha-k_\beta>1$. We apply the gauge transformation 
\ben
Y_1(\lambda,u)=\Gamma_1(\lambda,u)Y(\lambda,u),\quad 
\Gamma_1 =1-(U_0 E_{\beta\alpha} U_0^{-1}) g^{-1}\lambda^{-1}.
\een
The resulting matrix-valued function is a fundamental solution to a Fuchsian
system that has the same type of expansion at $\lambda=u_i$ for
$i=1,\dots, N+1$ while at $\lambda=0$ we have
$Y_1(\lambda,u)=U^1(\lambda,u)\lambda^{K^1}$  (see Section
\ref{sec:gauge}). 
Note that 
\ben
\operatorname{Tr}((K^1)^2) -\operatorname{Tr}((K)^2) =
2-2k_\alpha+2k_\beta \leq -2.
\een
Repeating this process we get a sequence of fundamental matrices
\ben
Y_i(\lambda,u) = U^i(\lambda,u)\, \lambda^{K^i},
\een
satisfying
$\operatorname{Tr}((K^i)^2)<\operatorname{Tr}((K^{i-1})^2)$. Therefore,
the sequence stops after finitely many steps when $K^s=0$ for some
$s$. Let us denote by $g_\ell(u)$ the non-zero entry of
$
(U_0^{\ell-1})^{-1} U_1^{\ell-1}
$
that we choose in order to construct $Y_{\ell}$. The function
$g_1:=g\in \O_{T,t^*}$ is holomorphic, but the remaining ones are
meromorphic at $t^*$, i.e., $g_\ell\in
\operatorname{Frac}(\O_{T,t^*})$ ($2\leq \ell \leq s$), where for an
integral domain $R$ we denote by 
$\operatorname{Frac}(R)$ the quotient field of $R$. Put
$g_\ell=b_\ell/h_\ell$, where $b_\ell,h_\ell\in \O_{T,t^*}$
are relatively prime (recall that $\O_{T,t^*}$ is a UFD). According to
our construction $h_1=1$ and for each $\ell>1$ there exists an integer
$m_\ell\geq 0$ such that $h_\ell$ is a divisor of
$(b_1\cdots b_{\ell-1})^{m_\ell}$. Shrinking the
neighborhood $\U$ of $u^*$ if necessary we may assume that $b_\ell$
are represented by holomorphic functions in $\U$. Put
$\tau^*(u):=b_1(u)\cdots b_s(u)$. Then since $K^s=0$ we get that  $f\,
Y_s(\lambda,u)^{-1}$  is a global trivializing frame for
$E|_{\PP^1\times \{u\}}$ for all $u$ such that $\tau^*(u)\neq 0$. In
particular, the analytic germ   
\ben
(\Theta,t^*)\subset \{ \tau^*(t)=0\}. 
\een
On the other hand if $\widetilde{Y}(\lambda,u)$ is the fundamental
solution for the Schlesinger deformation of $\nabla^\circ$ satisfying
the initial condition $Y(\lambda,t^\circ)=Y^\circ(\lambda)$, then $f\,
\widetilde{Y}(\lambda,u)^{-1}$ is also a global trivialization  of
$E|_{\PP^1\times \{u\}}$ for all $u\notin \Theta$. Therefore
$Y_s(\lambda,u) = C \widetilde{Y}(\lambda,u) $ for some constant
invertible matrix $C$. Therefore the connection 1-form of the Fuchsian
connection corresponding to $Y_s$ is conjugated with the connection
1-form corresponding to $\widetilde{Y}$ via $C$. We get that 
$\omega^s$ coincides with the 1-form \eqref{tau-form}. Recalling Lemma
\ref{le:form-gauge} we get
\ben
\omega = \omega^s = \omega^* + \frac{d(g_1\cdots g_s)}{g_1\cdots g_s}.
\een
It remains only to factorize $g_1(u)\cdots g_s(u)=\tau_1^*(u)^{r_1}\cdots
\tau_s^*(u)^{r_s} h(u)$ where $\tau_i^*=0$ are the local equations of
the irreducible components of $\Theta$ at $t^*$ and $h$ is relatively
prime to $\tau_1^*\cdots \tau_s^*$. Note that $h(u^*)\neq 0$, otherwise the form
$\omega$ will have a pole along the hypersurface $\{h=0\}$ which is
not contained in $\Theta$. Lemma \ref{le:theta-loc}
follows. 
\qed

\medskip

Since $g_1$ is holomorphic, if the algorithm stops on the first step
we would get that the isomonodromic tau-function is analytic at
$t^*$. Therefore we get the following corollary.
\begin{corollary}
If $E|_{\PP^1\times \{t^*\}}\cong \O(1)\oplus \O(-1)\oplus \O^{p-2}$
for a generic $t^*\in \Theta$, then the isomonodromic tau-function is
analytic on $T$. 
\end{corollary}

\begin{remark}
One can prove that if $t^*\in \Theta$ is a smooth point and we choose
$g_1$ appropriately, then $g_2$ is also holomorphic. The analysis of the
poles of $g_\ell$ for $\ell>2$ however becomes more difficult. 
\end{remark}
\begin{remark}
According to Lemma \ref{le:g-value} the
subset $\Theta_K\subset \Theta$ (see Remark \ref{re:theta-K}) 
is contained in the zero locus of $b_1=g_1$ (the function generated on
the first step of the algorithm). In particular, if $t^*\in \Theta$ is
generic (so that $\Theta_K$ is open in $\Theta$) then the first step of the
algorithm already determines the germ of $\Theta$ at $t^*$.  
\end{remark}

\section{Frobenius manifolds}
The main goal of this lecture is to recall the notion of a semi-simple Frobenius manifold and to prove that semi-simple Frobenius manifolds can be classified by solutions of a certain system of PDEs. The general reference for more details is \cite{Du} (see also  \cite{Man}).

\subsection{Definition}
There are several ways to introduce the notion of a Frobenius manifold. We have chosen a set of axioms convenient for our purposes. Our definition is equivalent to (Definition 1.2 in \cite{Du}). 
Let $M$ be a complex manifold and $\T_M$ denotes the sheaf of holomorphic vector fields on $M$. Let us assume that $M$ is equipped with the following structures
\begin{enumerate}
\item[(a)] Each tangent space $T_tM$, $t\in M$,  is equipped with the structure of a {\em Frobenius algebra} depending holomorphically on $t$. In other words, we have a commutative associative multiplication $\bullet_t$ and symmetric non-degenerate bi-linear  pairing $(\ ,\ )_t$  satisfying the Frobenius property
\ben
(v_1\bullet_t w,v_2)=(v_1,w\bullet_t v_2),\quad v_1,v_2,w\in T_tM.
\een 
The pointwise multiplication $\bullet_t$ defines a multiplication $\bullet$ in $\T_M$, i.e., a $\O_M$-bilinear map 
\ben
\T_M\otimes \T_M\to \T_M, \quad v_1\otimes v_2\mapsto v_1\bullet v_2.
\een 
The pairing $(\ ,\ )_t$ determine a $\O_M$-bilinear pairing 
\ben
(\ ,\ ): \T_M\otimes \T_M \to \O_M.
\een 
\item[(b)] There exists a global vector field $e\in \T_M$, called {\em unit vector field},  such that 
\ben
\nabla_v^{\rm L.C.} e=0,\quad e\bullet v=v,\quad \forall v\in \T_M,
\een
where $\nabla^{\rm L.C.}$ is the Levi--Civita connection on $\T_M$ corresponding to the bi-linear pairing $(\ ,\ )$. 
\item[(c)] There exists a global vector field $E\in \T_M$, called {\em Euler vector field}, such that 
\ben
E(v_1,v_2)-([E,v_1],v_2)-(v_1,[E,v_2]) = (2-D) (v_1,v_2),
\een
for all $v_1,v_2\in \T_M$ and  for some constant $D\in \CC$. 
\end{enumerate}
The above data allows us to define the so-called {\em structure connection} $\nabla$ on the vector bundle $\operatorname{pr}_M^*TM\to M\times \CC^*$, where 
\ben
\operatorname{pr}_M:M\times \CC^*\to M, \quad (t,z)\mapsto t
\een
is the projection map. Namely,
\ben
\nabla_v & := & \nabla^{\rm L.C.}_v -z^{-1} v\bullet,\quad v\in \T_M\\
\nabla_{\partial/\partial z} & := & \frac{\partial}{\partial z} -z^{-1} \theta + z^{-2} E\bullet,
\een
$v\bullet$ and $E\bullet$ are $\O_M$-linear maps $\T_M\to \T_M$ corresponding to the Frobenius multiplication by respectively $v$ and $E$. The $\O_M$-linear map $\theta:\T_M\to \T_M$ is defined by 
\ben
\theta(v):=\nabla^{\rm L. C.}_vE - (1-D/2) v.
\een
The operator $\theta$ is sometimes called {\em Hodge grading operator}. Let us point out that the term $(1-D/2)v$ in the definition of $\theta(v)$ is inserted so that $\theta$ becomes skew-symmetric with respect to the Frobenius pairing
\ben
(\theta(v_1),v_2)+(v_1,\theta(v_2)) = 0,\quad v_1,v_2\in \T_M. 
\een
\begin{definition}\label{def:frob}
The data $((\ ,\ ), \bullet, e, E)$ satisfying the conditions (a), (b), and (c) from above is said to be a  {\em Frobenius structure} on $M$ of {\em conformal dimension} $D$ if the structure connection $\nabla$ is flat. 
\end{definition}

\subsection{Properties}
The following proposition is a direct consequence of the definition.
\begin{proposition}\label{prop:Frob-properties}
Suppose that $(M,(\ ,\ ),\bullet, e,E)$ is a Frobenius structure. Then 

a) The Levi--Civita connection $\nabla^{\rm L.C.}$ is flat. 

b) Let $t=(t_1,\dots,t_N)$ be $\nabla^{\rm L.C.}$-flat coordinates defined on a contractible open subset $U\subset M$. There exists a holomorphic function $F\in \O_M(U)$, such that 
\ben
(\partial/\partial t_a\bullet \partial/\partial t_b,\partial/\partial t_c) = 
\frac{\partial^3F}{\partial t_a\partial t_b\partial t_c}
\een
and 
\ben
EF = (3-D)F + H, 
\een
where $H$ is a polynomial in $t_1,\dots,t_N$ of degree at most 2.

c) The  Hodge grading operator is covariantly constant: $\nabla^{\rm L.C.}_v\theta=0$. In particular, in flat coordinates $t=(t_1,\dots,t_N)$ the matrix $(\theta_{ab})_{a,b=1}^N$ of $\theta$ defined by 
\ben
\theta(\partial /\partial t_b) = \sum_{a=1}^N \theta_{ab} \partial/\partial t_b
\een
is constant. 

d) The following identity holds
\ben
[E,v\bullet w]-[E,v]\bullet w-v\bullet [E,w] =v\bullet w,\quad v,w\in \T_M.
\een
\end{proposition}
\proof 
Parts a) and b) are straightforward. We will prove c) and d) simultaneously. To begin with note that both c) and d) are $\O_M$-linear in $v$ and $w$. Therefore, we may assume $v$ and $w$ are flat with respect to $\nabla^{\rm L.C.}$. The flatness of $\nabla$ implies that 
\ben
\nabla_{z\partial_z +E}\nabla_v w-\nabla_v \nabla_{z\partial_z +E}w -\nabla_{[E,v]} w = 0.
\een
By definition 
\ben
\nabla_v=\nabla_v^{\rm L.C.} - z^{-1} v\bullet
\een
and 
\ben
\nabla_{z\partial_z +E} = z\partial_z + \nabla_E^{\rm L.C.} -\theta.
\een
Substituting these operators in the 0-curvature equation and using that $v$ and $w$ are flat we get a polynomial expression in $z^{-1}$ of degree 1 for which the coefficient in front of $z^0$ is 
\ben
\nabla^{\rm L.C.}_v\theta(w)
\een
and the coefficient in front of $z^{-1}$ is
\beq\label{zinv-coeff}
v\bullet w +[E,v]\bullet w - v\bullet \theta(w) +\theta(v\bullet w)-\nabla^{\rm L.C.}_E(v\bullet w).
\eeq
Therefore both expressions must vanish. The vanishing of $\nabla^{\rm L.C.}_v\theta(w)$ for all flat vector fields $v$ and $w$ is equivalent to the statement in c). For the 2nd expression, using the definition of $\theta$ we get
\ben
- v\bullet \theta(w) +\theta(v\bullet w) = 
v\bullet [E,w] - [E,v\bullet w]+ \nabla^{\rm L.C.}_E(v\bullet w).
\een
Substituting this identity in \eqref{zinv-coeff} we get the identity of part d).
\qed

Note that locally the Frobenius structure is completely determined by
the Euler vector field $E$ and the holomorphic function $F$. It is
possible to state the definition in terms of $F$ as well (see
\cite{Du}). This leads to the so-called WDVV equations for $F$. In
many applications the Frobenius structures arise as solutions of the
WDVV equations. However, in our lectures this point of view would not
play an important role.  

\subsection{Example: quantum cohomology} 

Let $X$ be a smooth projective variety. Recall that a stable map $(\Sigma,z_1,\dots,z_n;f)$ is a holomorphic map $f:\Sigma\to X$, where $\Sigma$ is a nodal Riemann surface, $z_i$ are marked points (pairwise distinct and non-singular), such that the automorphism group of $(\Sigma,z_1,\dots,z_n;f)$ is finite. The homology class $d=f_*[\Sigma]\in H_2(X,\ZZ)$ is called the {\em degree} of the stable map. Let us denote by $\mathcal{M}_{g,n}(X,d)$ the moduli space of stable maps $(\Sigma,z_1,\dots,z_n;f)$ such that the arithmetic genus of $\Sigma$ is $g$, the number of marked points is $n$, and the degree of $f$ is $d$. This is a proper Delign--Mumford stack equipped with a virtual fundamental cycle $[X_{g,n,d}]$ of dimension (over $\CC$) 
\ben
3g-3+n +(1-g) D + \int_{[X]} c_1(TX),
\een
where $D=\operatorname{dim}_\CC(X)$. The Gromov--Witten invariants of $X$ are defined by the following correlators
\ben
\langle \alpha_1,\dots,\alpha_n\rangle_{g,n,d} = \int_{[X_{g,n,d}]} 
\operatorname{ev}^*(\alpha_1,\dots,\alpha_n),\quad \alpha_i\in H^*(X;\CC),
\een
where 
\ben
\operatorname{ev}:\mathcal{M}_{g,n}(X,d)\to X^n,\quad 
(\Sigma,z_1,\dots,z_n;f)\mapsto (f(z_1),\dots,f(z_n))
\een
is the evaluation map. 

In order to define quantum cohomology we need also to recall the definition of the {\em Novikov ring}. The degrees of stable maps form a cone in $H_2(X,\ZZ)$ usually denoted by $\operatorname{Eff}(X).$ The Novikov ring is by definition the formal group algebra of $\operatorname{Eff}(X)$, i.e., 
\ben
\CC[Q]:=\Big\{ \sum_{d\in \operatorname{Eff}(X)} c_d Q^d\ |\ c_d\in \CC\Big\}.
\een
Let us fix a set of ample line bundles $L_1,\dots,L_r$ on $X$, such that $p_i:=c_1(L_i)$ form a $\CC$-basis of $H^{1,1}(X;\CC)$. Then the map 
\ben
Q^d\mapsto Q_1^{\langle p_1,d\rangle}\cdots Q_r^{\langle p_r,d\rangle}
\een
gives an embedding 
\ben
\CC[Q]\to \CC[\![Q_1,\dots,Q_r]\!].
\een
Let $H^{\rm ev}(X;\CC):=\oplus_{d=0}^D H^{2d}(X;\CC)$ and  let us fix a homogeneous basis $\{\phi_i\}_{i=1}^N$ of $H^{\rm ev}(X;\CC)$ such that $\phi_1=1$ and $\phi_{i+1}=p_i$ for $1\leq i\leq r$. Put $t=\sum_i t_i\phi_i$. Then  the genus-0 potential of $X$ is defined by 
\ben
F^{(0)}(Q,t) = \sum_{n=0}^\infty \sum_{d\in \operatorname{Eff}(X)} \frac{Q^d}{n!}
\langle t,\dots,t\rangle_{0,n,d}.
\een
The GW invariants satisfy the so-called divisor equation which implies that $\partial_{t_{i+1}}F^{(0)} = Q_i\partial_{Q_i} F^{(0)}$ for all $1\leq i\leq r$. Therefore, the genus-0 potential has the form
\ben
F^{(0)}(Q,t) = F(t_1,Q_1e^{t_2},\dots, Q_re^{t_{r+1}},t_{r+2},\dots, t_N).
\een
Let us fix $Q_1,\dots, Q_r$ as complex parameters (e.g. set $Q_i=1$ for all $i$). In many important examples, the formal series defining quantum cohomology is convergent on a domain 
\ben
M=\{t\in H^*(X;\CC)\ |\ \operatorname{Re}(t_{i+1})<-R,\mbox{ for $1\leq i\leq r$},\quad |t_j|<\epsilon \mbox{ for $r+1\leq j\leq N$} \},
\een
where $R>0$  and $\epsilon>0$ are real numbers.  Let us introduce also the vector fields 
\ben
e=\partial/\partial t_1,\quad E= \sum_{i=1}^N ((1-d_i) t_i+r_i)\partial_{t_i},
\een
where $d_i=\operatorname{deg}(\phi_i)/2$ and $r_i$ are the coordinates of $c_1(TX)$, i.e., $c_1(TX)=:\sum_{i=1}^N r_i \phi_i$. 
If the domain of convergence $M$ exists, then the Poincare pairing, the vector fields $e$ and $E$, and the multiplication defined in terms of $F^{(0)}$ via the formulas of Proposition \ref{prop:Frob-properties}, part b),  determine a Frobenius structure on $M$ of conformal dimension $D=\operatorname{dim}_\CC(X). $

\subsection{Semi-simple Frobenius manifolds}
\begin{definition}
A Frobenius manifold $(M,(\ ,\ ),\bullet,e,E)$ is said to be {\em semi-simple} if there are local coordinates $u=(u_1,\dots,u_N)$ defined in a neighborhood of some point on $M$ such that 
\ben
\partial/\partial u_i\bullet \partial/\partial u_i = \delta_{ij}\partial/\partial u_j,\quad 1\leq i,j\leq N.
\een
The coordinates $u_i$ are called {\em canonical coordinates}.
\end{definition}

As we will see now, canonical coordinates are unique up to parmutation and constant shifts. To avoid cumbersome notation we put $\partial_{u_i}:=\partial/\partial u_i$. 
\begin{proposition}
Let $u=(u_1,\dots,u_N)$ be canonical coordinates defined on some open subset $U\subset M$. Then 

a) The Frobenius pairing takes the form 
\ben
(\partial_{u_i},\partial_{u_j}) = \delta_{ij}\eta_j(u),\quad 1\leq i,j\leq N,
\een
where $\eta_j\in \O_M(U)$ and $\eta_j(u)\neq 0$ for all $u\in U$. 

b) The unit vector field takes the form $e=\sum_{i=1}^N \partial_{u_i}$.

c) The 1-form $\sum_{i=1}^N \eta_i(u) du_i$ is closed.

d) There are constants $c_i$ ($1\leq i\leq N$) such that 
\ben
E=\sum_{i=1}^N (u_i+c_i)\partial_{u_i}.
\een 
\end{proposition}
\proof
a) If $i\neq j$ then we have 
\ben
(\partial_{u_i},\partial_{u_j}) = (e\bullet \partial_{u_i},\partial_{u_j}) = (e,\partial_{u_i}\bullet\partial_{u_j}) =0.
\een
The fact that $\eta_i(u):=(\partial_{u_i},\partial_{u_i})\neq 0$ follows from the non-degeneracy of the Frobenius pairing. 

b) Let $e=\sum_{i=1}^N e_i(u) \partial_{u_i}$. Then 
\ben
\partial_{u_j} = \partial_{u_j} \bullet e = e_j(u) \partial_{u_j}.
\een
Therefore $e_j(u)=1$ for all $j$.

c) We have to check that $\partial_{u_j}\eta_i =\partial_{u_i}\eta_j$. On the other hand
\ben
\partial_{u_j}\eta_i = \partial_{u_j}(\partial_{u_i},e) = 
(\nabla^{\rm L.C.}_{\partial_{u_j}}\partial_{u_i},e),
\een
where we used the Leibnitz rule and the fact the $e$ is a flat vector field. It remains only to recall that the Levi--Civita connection is torsion free, so 
\ben
\nabla^{\rm L.C.}_{\partial_{u_j}}\partial_{u_i}=
\nabla^{\rm L.C.}_{\partial_{u_i}}\partial_{u_j}.
\een
d) Put $E=\sum_{i=1}^N E_i(u)\partial_{u_i}$. Let us recall Proposition \ref{prop:Frob-properties}, part d) with $v=\partial_{u_i}$ and $w=\partial_{u_j}$. For $i\neq j$ we get 
\ben
(\partial_{u_i}E_j) \partial_{u_i}+(\partial_{u_j}E_i) \partial_{u_j} =0.
\een 
Hence $\partial_{u_i}E_j=0$ for $i\neq j$. If $i=j$ then we get $\partial_{u_i}E_i=1$. Therefore $E_i(u)=u_i+c_i$ for some constant $c_i$. 
\qed

Part d) of the above proposition shows that in every canonical coordinate system up to some constant shifts the canonical coordinates coincide with the eigenvalues of the operator $E\bullet$. Therefore, up to constant shifts and permutations the canonical coordinates are uniquely determined. From now on we will work only with canonical coordinates such that 
\ben
E=\sum_{i=1}^N u_i\partial_{u_i}.
\een
The question that we would like to answer now is the following. Let us assume that $U$ is an open subset of the universal cover $T$ of $Z_N$ and $\sum_{i=1}^N\eta_i(u) du_i$ is a closed 1-form on $U$. The tangent bundle of $T$ and hence of $U$ as well is trivial, because $T$ is a contractible Stein manifold, so according to the Grauert--Oka principle every holomorphic vector bundle on $T$ is trivial. Alternatively, we can prove that $\T_T$ is a free $\O_T$-module by using that the vector fields $\partial_{u_i}$ of the configuration space $Z_N$ lift naturally to vector fields on $T$ and provide a global trivialization of $\T_T$. Using the 1-form we define a pairing
\ben
(\partial_{u_i},\partial_{u_j})=\delta_{ij} \, \eta_j(u).
\een
Let us also define multiplication
\ben
\partial_{u_i}\bullet \partial_{u_j}= \delta_{ij}\partial_{u_j}
\een
and vector fields 
\ben
e=\sum_{i=1}^N \partial_{u_i},\quad  E =\sum_{i=1}^N u_i\partial_{u_i}.
\een
The problem then is to classify all 1-forms $\sum_{i=1}^N \eta_i(u)du_i$ such that the above data determines a Frobenius structure on $U$.   The answer is given by the following theorem.
\begin{theorem}\label{thm:Frob-class}
The closed 1-form $\sum_{i=1}^N \eta_i(u)du_i$ determines a Frobenius structure on $U$ of conformal dimension $D$ if and only if the following conditions are satisfied
\begin{enumerate}
\item[(1)]  $\eta_i(u)\neq 0$ for all $i$ and for all $u\in U$.
\item[(2)] $e\eta_i(u)=0$ for all $i$.
\item[(3)] $E\eta_i(u)=-D \eta_i(u)$.
\item[(4)] For all $k\neq i\neq j\neq k$ we have 
\ben
\frac{\partial \eta_{ij}  }{ \partial u_k} = \frac{1}{2} \Big(
\frac{\eta_{ij}\eta_{kj} }{\eta_j} + 
\frac{\eta_{jk}\eta_{ik} }{\eta_k}+
\frac{\eta_{ki}\eta_{ji} }{\eta_i}
\Big),
\een
where $\eta_{ab}(u):=\partial_{u_a}\eta_b(u)$. 
\end{enumerate}
\end{theorem}
\proof
{\bf Step 1.} Determine when does the 1-form $\sum_i \eta_i(u)du_i$ defines data satisfying conditions (a), (b), and (c) in the definition of a Frobenius manifold. 

In part (a), we would like the multiplication and the pairing to give a holomorphic family of Frobenius algebras. This is clearly satisfied for any choice of the 1-form. The requirement that the pairing is non-degenerate yields that $\eta_i(u)\neq 0$ for all $i$ and for all $u\in U$. 

For condition (b), we would like to know when is $e$ a flat vector field. Let $\Gamma_{ij}^k$ be the Christoffel's symbols of the pairing $g_{ij}(u)=\delta_{ij}\eta_j$. A straightforward computation yields
\ben
\Gamma_{ij}^j =\Gamma_{ji}^j = \frac{\eta_{ij}}{2\eta_j},\quad 1\leq i,j\leq N,
\een
\ben
\Gamma_{ii}^j = -\frac{\eta_{ij}}{2\eta_j},\quad 1\leq i\neq j\leq N,
\een
and 
\ben
\Gamma_{ij}^k =0 ,\quad k\neq i\neq j\neq k.
\een
Using the above formulas we compute directly that 
\ben
\nabla^{\rm L.C.}_{\partial_{u_i}} e = \frac{e\eta_i}{2\eta_i}\, \partial_{u_i}.
\een
Therefore $e$ is a flat vector field if and only if $e\eta_i=0$ for all $i$.  

Finally, for condition (c) to hold we must have $E\eta_i=-D\eta_i$ for all $i$. Therefore, the 1-form will define a data satisfying conditions (a), (b), and (c) if and only if the functions $\eta_i(u)$ satisfy conditions (1), (2), and (3) in Theorem \ref{thm:Frob-class}. 

{\bf Step 2.} When is the Levi--Civita connection flat? 

The flatness of $\nabla^{\rm L.C.}$ is equivalent to: the expression 
\ben
2(\nabla^{\rm L.C.}_{\partial_{u_i}}\nabla^{\rm L.C.}_{\partial_{u_j}} \partial_{u_k},\partial_{u_\ell}) 
\een
is symmetric in $i$ and $j$. Using the Leibnitz rule we transform this expression into 
\beq\label{Rijkl}
\partial_{u_i} \Big(2\Gamma_{jk}^\ell \eta_\ell\Big) -\sum_{a=1}^N 2\Gamma_{jk}^a \Gamma_{i\ell}^a \eta_a.
\eeq
Let us assume first that $i,j$, and $k$ are pairwise distinct. Then we get 
\ben
\delta_{\ell j} \Big( \frac{\partial \eta_{jk}}{\partial u_i} - \frac{\eta_{ij}\eta_{kj}}{2\eta_j}\Big)  +
\delta_{\ell i} \Big( \frac{\eta_{kj} \eta_{ij}}{2\eta_j} + \frac{\eta_{jk}\eta_{ik}}{2\eta_k}\Big)+
\delta_{\ell k} \Big( \frac{\partial \eta_{jk}}{\partial u_i} - \frac{\eta_{i k}\eta_{jk}}{2\eta_k}\Big) .
\een
The last term is symmetric in $i$ and $j$, so a non-trivial condition will be obtained either if $\ell=i$ or $\ell=j$. Due to the symmetry between $i$ and $j$ we may assume that $\ell=j$. Then we get 
\ben
\frac{\partial \eta_{jk}}{\partial u_i} - \frac{\eta_{ij}\eta_{kj}}{2\eta_j} = 
\frac{\eta_{ki} \eta_{ji}}{2\eta_i} + \frac{\eta_{ik}\eta_{jk}}{2\eta_k}.
\een
This is exactly the PDE given in condition (4). 

There are 3 more cases to analyze. Indeed, since we may assume that $i\neq j$ we get that $k=i$ or $k=j$. Again exchanging the role of LHS and RHS provides a symmetry between $i$ and $j$, which allows us to assume that $k=i$. Therefore the remaining cases are: 
$(k,\ell)=(i,i),(i,j),$ or $k=i$ and $\ell\neq i,j$. The first case yields $\partial_{u_i}\eta_{ij}=\partial_{u_j}\eta_{ii}$, which is always satisfied because the 1-form $\sum_i \eta_idu_i$ is closed. The 2nd case $(k,\ell)=(i,j)$ yields
\ben
\partial_{u_i}\eta_{ij}+\partial_{u_j}\eta_{ij} = 
\frac{\eta_{ij}\eta_{ij}}{2\eta_i} + \frac{\eta_{ij}\eta_{ij}}{2\eta_j} +
\frac{\eta_{ij}\eta_{ii}}{2\eta_i} + \frac{\eta_{ij}\eta_{jj}}{2\eta_j} 
-\sum_{a:a\neq i,j} 
\frac{\eta_{ia}\eta_{ja}}{2\eta_a} .
\een
It is easy to see that this identity is a consequence of (2) and (4). 
In the last case if $k=i$ and $\ell\neq i,j$ we get 
\ben
\partial_{u_j} \eta_{i\ell} = \frac{\eta_{ij}\eta_{i\ell}}{2\eta_i} +
\frac{\eta_{i\ell}\eta_{j\ell}}{2\eta_\ell}+
\frac{\eta_{ij}\eta_{j\ell}}{2\eta_j} 
\een
which is equivalent to (4). 

{\bf Step 3.}
It remains only to verify that under the conditions (1)--(4) the structure connection $\nabla$ is flat. The argument is similar to the argument in Step 2, so it will be left as an exercise.
\qed

\section{Painleve property for semi-simple Frobenius manifolds}

\subsection{The second structure connection}

Let $U\subset Z_N$ be an open contractible neighborhood of some fixed point $u^\circ\in Z_N$. Suppose that $U$ is equipped with a semi-simple Frobenius structure $((\ ,\ ), \bullet, e, E)$. Put $H=T_{u^\circ}U$ and let us trivialize the tangent bundle 　
\beq\label{tb-triv}　
TU\cong U\times H\cong U\times \CC^N
\eeq
using the Levi--Civita connection. In other words, we fix a basis $\{\phi_a\}_{a=1}^N$ of $H$ and let $\partial_{t_a}\in \T_U$ be the flat vector field on $U$ obtained by parallel transport with respect to the Levi--Civita connection. Then the isomorphisms \eqref{tb-triv} are given by the maps
\ben
(u,v) \in TU\mapsto 
(u, v_1 \phi_1+\cdots+ v_N\phi_N) \in U\times H\mapsto 
(u,v_1,\dots,v_N)\in U\times \CC^N,
\een 
where $v\in T_uU$ and $v=: v_1 \partial_{t_1}+\cdots + v_N \partial_{t_N}$. 
The isomorphism \eqref{tb-triv} identifies the structure connection of the Frobenius structure with the flat connection on the trivial bundle 
\ben
(U\times \CC^*)\times \CC^N\to U\times \CC^*
\een
defined by 
\ben
\nabla_{\partial_{u_i} } & = & \partial_{u_i} -z^{-1} P_i(u),\quad 1\leq i\leq N ,\\
\nabla_{\partial_{z\phantom{_i}} } & = & \partial_{z} -z^{-1}\theta + z^{-2} \mathcal{E}(u),
\een
where $P_i:U\to \mathfrak{gl}(\CC^N)$ is a holomorphic map whose $(a,b)$-entry $P_{iab}(u)$ is defined by the identity 
\ben
\partial_{u_i}\bullet \partial_{t_b} = \sum_{a=1}^N P_{iab}(u) \partial_{t_a} ,
\een
$\mathcal{E}=\sum_{i=1}^N u_i P_i(u)$, and $\theta$ is a constant matrix whose $(a,b)$-entry $\theta_{ab}$ is defined by 
\ben
\theta(\partial_{t_b}) = [\partial_{t_b},E]-(1-D/2) \partial_{t_b} =:\sum_{a=1}^N \theta_{ab} \partial_{t_a}. 
\een 
In order to justify the definition of the second structure connection we make the following heuristic argument. Suppose that the structure connection has a solution 
\ben
J:U\times \CC^*\to \CC^N
\een
given by a Laplace transform 
\ben
J(u,z)=\frac{(-z)^{n-\frac{1}{2}} } {\sqrt{2\pi}} \int_\Gamma e^{\lambda/z} I^{(n)}(u,\lambda) d\lambda
\een
along an appropriate contour $\Gamma\subset \CC$ of some $\CC^N$-valued function $I^{(n)}(u,\lambda)$ holomorphic for all $(u,\lambda)\in U\times \Gamma$. Here $n\in \CC$ is an arbitrary number. Assuming that the Laplace transform works, we would get that $J(u,z)$ is a solution to the structure connection if and only if  $I^{(n)}(u,\lambda)$ is a solution to the following connection 
\ben
\nabla^{(n)}_{\partial_{u_i} } & = & \partial_{u_i}  + (\lambda-\mathcal{E})^{-1} P_i(u) (\theta-n-1/2),\quad 
1\leq i\leq N, \\
\nabla^{(n)}_{\partial_{\lambda} } & = & \partial_{\lambda}  - (\lambda-\mathcal{E})^{-1} (\theta-n-1/2). \\
\een
This is a connection on 
\ben
(U\times \CC)'\times \CC^N\to (U\times \CC)',
\een
where 
\ben
(U\times \CC)' = \{ (u,\lambda)\in U\times \CC\ |\ \operatorname{det}(\lambda-\mathcal{E})\neq 0\}. 
\een
\begin{proposition}\label{deln-flat}
The connection $\nabla^{(n)}$ is flat for all $n\in \CC$. 
\end{proposition}
The proof is left as an exercise. The connection $\nabla^{(n)}$ is called the {\em second structure connection}. 

\subsection{Proof of Theorem \ref{t1}}
\begin{lemma}\label{le:psi}
Let $\widetilde{\Psi}$ be the matrix whose $(a,i)$-entry is given by $\widetilde{\Psi}_{ai} = \partial t_a/\partial u_i$. Then 
\ben
\widetilde{\Psi}^{-1} P_i \widetilde{\Psi} = E_{ii},\quad \widetilde{\Psi}^{-1}\mathcal{E}\widetilde{\Psi}=\operatorname{diag}(u_1,\dots,u_N),
\een
where $E_{ii}$ is the matrix whose entry in position $(i,i)$ is 1 and all other entries are $0$. 
\end{lemma} 
\proof
We have 
\ben
\partial_{u_i}\bullet \partial_{t_b} = 
\partial_{u_i}\bullet \sum_{j=1}^N \frac{\partial u_j}{\partial t_b} \, \partial_{u_j} = 
\frac{\partial u_i}{\partial t_b} \, \partial_{u_i} =
\sum_{a=1}^N \frac{\partial u_i}{\partial t_b}\, \frac{\partial t_a}{\partial u_i}\, \partial_{t_a}.
\een
Therefore
\ben
P_{iab} = \frac{\partial u_i}{\partial t_b}\, \frac{\partial t_a}{\partial u_i}.
\een
Using this formula we find that the $(a,j)$-entry of $P_i\widetilde{\Psi}$ is
\ben
\sum_{b=1}^N 
\frac{\partial u_i}{\partial t_b}\, 
\frac{\partial t_a}{\partial u_i}\,
\widetilde{\Psi}_{bj} = 
\sum_{b=1}^N 
\frac{\partial u_i}{\partial t_b}\, 
\frac{\partial t_a}{\partial u_i}\,
\frac{\partial t_b}{\partial u_j} =
\delta_{ij} \frac{\partial t_a}{\partial u_i} = \delta_{ij}\widetilde{\Psi}_{aj}.
\een
The latter is precisely the $(a,j)$-entry of $\widetilde{\Psi} E_{ii}$. Therefore $P_i\widetilde{\Psi} = \widetilde{\Psi} E_{ii}$. 
\qed

\begin{lemma}\label{le:Frob-Schl}
Let $n\in \CC$ be arbitrary. Then the matrix-valued functions 
\ben
A_i^{(n)}(u):= P_i(u) (\theta-n-1/2),\quad 1\leq i\leq N,
\een
satisfy the Schlesinger equations.
\end{lemma}
\proof
We have to prove that the connection 
\ben
\widetilde{\nabla}^{(n)}_{\partial_{u_i}}  & = &  
\partial_{u_i}+ \frac{A_i^{(n)}(u)}{\lambda-u_i} ,\quad
1\leq i\leq N \\
\widetilde{\nabla}^{(n)}_{\partial_\lambda} & = &  
\partial_\lambda -\sum_{i=1}^N \frac{A_i^{(n)}(u)}{\lambda-u_i}
\een
is flat. However, using Lemma \ref{le:psi} we get
\ben
(\lambda-\mathcal{E})^{-1} P_i (\theta-n-1/2)  = \frac{A_i^{(n)}(u)}{\lambda-u_i}.
\een
Therefore $\widetilde{\nabla}^{(n)}=\nabla^{(n)}$, so it remains only to recall Proposition \ref{deln-flat}.
\qed

\medskip

The proof of Theorem \ref{t1} can be given as follows. Let us choose
$n\in \CC$ such that the operator $\theta-n-1/2$ is invertible. Then  
\ben
\eta_i(u)=(\partial_{u_i},\partial_{u_i})=(P_i(u) e, e) = (A_i^{(n)}(u) (\theta-n-1/2)^{-1} e,e).
\een 
According to Theorem \ref{thm:Malgrange-2} and Lemma
\ref{le:Frob-Schl} the functions $\eta_i(u)$ extend to meromorphic
functions on $T$.\qed 

\subsection{Special initial conditions}
In this section we are going to prove a theorem of Manin \cite{Man}
which answers the question of what kind of initial conditions for the
Schlesinger equations determine a semi-simple Frobenius
structures. Following Manin we introduce the following definition.
\begin{definition}\label{sic}
Let $H$ be a vector space equipped with a non-degenerate symmetric
bi-linear pairing $(\ ,\ )$ and a distinguished vector $e\in
H$. Suppose also that we have a set of linear operators $\theta,
\{P_i^\circ\}_{i=1}^N\in \mathfrak{gl}(H)$. The data $(H, (\ ,\ ), e,
\theta, \{P_i^\circ\}_{i=1}^N)$ is said to be a {\em special initial
  condition} if the following conditions are satisfied:
\begin{enumerate}
\item
$\theta$ is skew-symmetric: $(\theta(a),b)+(a,\theta(b))=0$ for all
$a,b\in H$.
\item $e$ is an eigenvector of $\theta$ with eigenvalue $D/2$.
\item The set $\{P_i^\circ\}_{i=1}^N$ is a complete set of
  orthogonal projectors of $H$, i.e., 
\begin{enumerate}
\item $P_i^\circ P_j^\circ = \delta_{ij} P_j^\circ$ for all
  $1\leq i,j\leq N$.
\item $P_1^\circ+\cdots + P_N^\circ =1$.
\item $(P_i^\circ(a),b)=(a,P_i^\circ(b))$ for all $1\leq i\leq N$ and for all
  $a,b\in H$.
\item $P_i^\circ e\neq 0$ for all $1\leq i\leq N$.\qed
\end{enumerate}
\end{enumerate}
\end{definition}

Suppose that $((\ ,\ ), \bullet, e, E)$ is a semi-simple Frobenius
structure on some complex manifold $M$ and that $u^\circ\in M$ is a
semi-simple point, i.e., a neighborhood of $u^\circ$ admits canonical
coordinates. Then the data 
\ben
H:=T_{u^\circ}M,\  (\ ,\ ),\ e,\ 
\theta:=\nabla^{\rm L.C.}E-(1-D/2),\  
P_i^\circ=P_i(u^\circ) , \ 1\leq i\leq N,
\een
is a special initial condition. In fact the only property that we did
not check yet is that $e$ is an eigenvector of $\theta$. However
\ben
\theta(e) = [e,E]-(1-D/2) e = e-(1-D/2)e = (D/2) e,
\een
where in the first equality we used that $e$ is flat and in the second
equality we used that $e=\sum_i \partial_{u_i}$ and $E=\sum_i
u_i\partial_{u_i}$. 
\begin{proposition}\label{prop:manin}
Given a special initial condition $(H, (\ ,\ ), e,
\theta, \{P_i^\circ\}_{i=1}^N)$ and a point $u^\circ \in Z_N$, then
there exists an open neighborhood $U\subset Z_N$ of $u^\circ$ and an
isomorphism $T_{u^\circ}U\cong H$ such
that the special initial condition is obtained from a semi-simple
Frobenius structure on $U$. 
\end{proposition}
\proof
Let $A_i^{(n)}(u),$ $1\leq i\leq N$ be solutions to the Schlesinger
equations such that 
\ben
A_i^{(n)}(u^\circ)=P_i^\circ (\theta-n-1/2). 
\een
If $n+\frac{1}{2}$ is not an eigenvalue of $\theta$, then we define 
\ben
P_i^{(n)}(u)=A_i^{(n)}(u) (\theta-n-1/2)^{-1}.
\een
\begin{lemma}\label{le:fos}
The set $\{P_i^{(n)}(u)\}_{i=1}^N$ is a complete set of orthogonal
projections for all $u$ sufficiently close to $u^\circ$.
\end{lemma}
\proof
Let us fix a basis $\{\phi_i\}_{i=1}^N$ of $H$ and identify
$\mathfrak{gl}(H)$ with the space of $p\times p$-matrices.  
Let $\A$ be the polynomial ring
\ben
\A = \CC[(u_i-u_j)^{\pm 1}\ :\ 1\leq i<j\leq N]\otimes \CC[A_1,\dots,A_N],
\een
where $A_i=(A_{iab})_{a,b=1}^N$ are matrix variables. We define
derivations $\partial_{u_1},\dots,\partial_{u_N}$ of $\A$
such that
\ben
\partial_{u_i} A_j & := & \frac{[A_i,A_j]}{u_i-u_j},\quad 1\leq i\neq
j\leq N ,\\
(\partial_{u_1}+\cdots +\partial_{u_N}) A_j & := & 0,
\een
and if $f\in \A$ depends only on $u_1,\dots,u_N$ then $\partial_{u_i}$
is defined to be the usual derivative. It is easy to check that these
differentiations pairwise commute, so $\A$ becomes a $\D$-module for
the ring $\D$ of differential operators on $Z_N$. 

Let us define $\mathcal{I}\subset \A$ to be the ideal generated by the
relations corresponding to conditions (a)--(c) in Definition
\ref{sic}. More precisely, we replace $P_i^\circ$ by
$A_i(\theta-n-1/2)^{-1}$ and take the entries of the corresponding
matrix identities as generators of $\mathcal{I}$. 
Condition (a) yields generators given by the entries of 
\ben
R_{ij}(A_1,\dots,A_N) = A_i (\theta-n-1/2)^{-1} A_j -\delta_{ij}
A_j,\quad 1\leq i,j\leq N.
\een
Condition (b) gives the entries of 
\ben
R(A_1,\dots,A_N) = A_1+\cdots +A_N -\theta+n+\frac{1}{2}.
\een 
Finally, condition (c) gives the entries of 
\ben
R_i(A_1,\dots,A_N) = A_i (\theta-n-1/2)^{-1} + (\theta+n+1/2)^{-1}
A_i^T ,\quad 1\leq i\leq N,
\een
where ${ }^T$ is the transposition operation in $\mathfrak{gl}(H)$ with
respect to the pairing $(\ ,\ )$. 

We claim that in order to prove the lemma it is enough to check that
$\mathcal{I}$ is $\D$-invariant. Indeed, condition (a) in Definition
\ref{sic} will be satisfied if 
$$R_{ij}(A_1^{(n)}(u),\dots,A_N^{(n)}(u) )=0.$$ On the other hand, the
Taylor series expansion of $R_{ij}(A_1^{(n)}(u),\dots,A_N^{(n)}(u) ) $
at $u=u^\circ$ has the form
\ben
\sum_{m_1,\dots , m_N=0}^\infty
\Big( \frac{\partial_{u_1}^{m_1} }{m_1!} \cdots
\frac{\partial_{u_N}^{m_N} }{m_N!}  
\ R_{ij}
\Big)( A_1^{(n)}(u^\circ),\dots,A_N^{(n)}(u^\circ) ) 
(u_1-u_1^\circ)^{m_1}\cdots (u_N-u_N^\circ)^{m_N}, 
\een
where we used that $A_i^{(n)}(u)$ solve the Schlesinger equations, so
the evaluation maps $A_i\mapsto A_i^{(n)}(u)$ are $\D$-equivariant. 
It remains only to notice that all Taylor's coefficients must vanish,
because $P_i^{(n)}(u^\circ)=P_i^\circ$ form a complete system of
orthogonal projections, so the evaluation
$R(A_1^{(n)}(u^\circ),\dots,A_N^{(n)}(u^\circ))=0$ for all generators
$R$ of $\mathcal{I}$ and hence for all $R\in \mathcal{I}$. 

Let us check that $\mathcal{I}$ is $\D$-invariant. We will prove only
that $\partial_{u_k} R_{ij}\in \mathcal{I}$ because the remaining
cases can be dealt in the same way. It is more convenient to prove
that 
\ben
d R_{ij}:= \sum_{k=1}^N \partial_{u_k} R_{ij}\otimes du_k \quad
\in\quad \mathcal{I}\otimes \Omega^1(Z_N), 
\een
where $\Omega^1(Z_N)$ denotes the ring of holomorphic 1-forms on
$Z_N$. By definition $dR_{ij}$ is
\ben
&&
\sum_{k: k\neq i} \frac{[A_k,A_i]}{u_k-u_i} (\theta-n-1/2)^{-1}
A_j \otimes (du_k-du_i) + \\
&& +
\sum_{k :k\neq j} A_i (\theta-n-1/2)^{-1}\frac{[A_k,A_j]}{u_k-u_j}
\otimes (du_k-du_j)  +\\
&& - 
\delta_{ij} \sum_{k :k\neq j} \frac{[A_k,A_j]}{u_k-u_j}
\otimes (du_k-du_j) 
\een
On the other hand
\ben
[A_k,A_i] (\theta-n-1/2)^{-1}
A_j  = \delta_{ij} A_k A_j - \delta_{kj} A_i A_j\quad ({\rm mod}\ \mathcal{I})
\een
and 
\ben
A_i (\theta-n-1/2)^{-1} [A_k,A_j] = \delta_{ik} A_k A_j -\delta_{ij}
A_j A_k \quad ({\rm mod}\ \mathcal{I})
\een
Therefore modulo terms in $\mathcal{I}$ the differential $dR_{ij}$
coincides with the sum of the following 2 terms
\ben
\delta_{ij} \Big( 
A_kA_j \otimes \frac{du_k-du_i}{u_k-u_i} - 
A_jA_k \otimes \frac{du_k-du_j}{u_k-u_j} -
[A_k,A_j] \otimes \frac{du_k-du_j}{u_k-u_j}\Big) 
\een
and 
\ben
-\delta_{kj} A_iA_j \otimes \frac{du_k-du_i}{u_k-u_i} +
\delta_{ik} A_kA_j \otimes \frac{du_k-du_j}{u_k-u_j} .
\een
Both terms vanish, which proves that the entries of $dR_{ij}$ are in
$\mathcal{I}\otimes \Omega^1(Z_N)$. 

This completes the proof that the set $\{P_i^{(n)}(u)\}_{i=1}^N$ satisfies conditions
(a)--(c) in Definition \ref{sic}. The last condition (d) will be satisfied for all $u$
sufficiently close to $u^\circ$, because $P_i^{(n)}(u)$ is continuous
and $P_i^{(n)}(u^\circ)e = P_i^\circ e \neq 0$. 
\qed

\begin{lemma}\label{le:nind}
If $n+\frac{1}{2}$ and $m+\frac{1}{2}$ are not eigenvalues of
$\theta$, then $P_i^{(m)}(u)=P_i^{(n)}(u)$.
\end{lemma}
\proof
According to Lemma \ref{le:fos} the matrices $P_i^{(n)}(u)$ pairwise
commute. Using that $A_i^{(n)}(u)$ satisfy the Schlesinger equations
we get
\ben
d P_i^{(n)}(u) = \sum_{j:j\neq i} \frac{du_j-du_i}{u_j-u_i} \, 
\Big( P_j^{(n)}(u) \theta P_i^{(n)}(u)-P_i^{(n)}(u) \theta P_j^{(n)}(u)\Big).
\een
Using these equations and the fact that $P_i^{(n)}(u)$ pairwise
commute we get that the matrix-valued functions 
$\widetilde{A}_i^{(n)}(u):= P_i^{(m)}(u) \Big(
\theta-n-\frac{1}{2}\Big)$ ($1\leq i\leq N$) satisfy the Schlesinger
equations. However the initial condition
$\widetilde{A}_i^{(n)}(u^\circ) = A_i^{(n)}(u^\circ)$. Therefore
$\widetilde{A}_i^{(n)}(u)= A_i^{(n)}(u).$
\qed

According to Lemma \ref{le:nind} the matrices $P_i(u):=P_i^{(n)}(u)$
are independent of $n$, while Lemma \ref{le:fos} implies that they
form a complete system of orthogonal projections.
\begin{lemma}\label{sic-Frob}
The  1-form 
\ben
\sum_{i=1}^N \eta_i(u)du_i,\quad \eta_i(u):=(P_i(u) e,e),\quad 1\leq
i\leq N,
\een
defines a Frobenius structure on every sufficiently small neighborhood
$U$ of $u^\circ$.
\end{lemma}
\proof
Let us first check that the above 1-form is closed. We have 
\ben
\eta_{ij}(u):=\partial_{u_j}\eta_i = \partial_{u_j} (P_i(u)e,e) =
\frac{2}{D-1-2n} (\partial_{u_j}A_i^{(n)}(u)e,e),
\een
where we used that $P_i(u)=A_i^{(n)}(u)\, (\theta-n-1/2)^{-1}$ and
that $\theta(e)=(D/2)e.$ We have to prove that
$\eta_{ij}(u)=\eta_{ji}(u)$. Let us assume that $i\neq j$. Since $A_i^{(n)}(u)$ ($1\leq i\leq N$)
satisfy the Schlesinger equations we get
\ben
\partial_{u_j}A_i^{(n)} = \frac{[A_j,A_i]}{u_j-u_i} = \partial_{u_i}A_j^{(n)},
\een
which implies that $\eta_{ij}=\eta_{ji}$, so the 1-form is closed. To
complete the proof we have to check that the 4 conditions of Theorem
\ref{thm:Frob-class} are satisfied. 

Note that the vectors $P_i^\circ e$ ($1\leq i\leq N$) form a basis of
$H$. Indeed, if $\sum_i \alpha_i P_i^\circ e=0$, then applying to both
sides $P_i^\circ$ we get $\alpha_i P_i^\circ e=0$. By assumption
$P_i^\circ e\neq 0$, so $\alpha_i=0$. The matrix of the form $(\ ,\ )$
is diagonal in the basis $P_i^\circ e$ with diagonal entries
$\eta_i(u^\circ).$ Therefore $\eta_i(u^\circ)\neq 0$ for all $i$
otherwise the form will be degenerate. By continuity there exists a
small neighborhood $U$ of $u^\circ$ such that $\eta_i(u)\neq 0$ for
all $i$ and for all $u\in U$. 

The second condition that we have to check is $e\eta_i=0$. This follows from
the fact that 
\ben
\sum_{j=1}^N \eta_j(u) = \Big(\sum_{j=1}^N P_j(u)e,e\Big) = (e,e)
\een
is a constant independent of $u$. 

The third condition that we have to check is $E\eta_i=-D\eta_i$. We have (see above)
\ben
\eta_i(u) = 
\frac{2}{D-1-2n} (A_i^{(n)}(u)e,e).
\een
Note that 
\ben
E A_i^{(n)}(u) = \iota_E d A_i^{(n)}(u) =
\iota_E \sum_{j:j\neq i}
\frac{du_j-du_i}{u_j-u_i}[A_j^{(n)}(u),A_i^{(n)}(u)] = [\theta,A_i^{(n)}(u)],
\een
where in the second equality we used the Schlesinger equations and
in the third one we used that 
$$
\sum_{j=1}^N A_j^{(n)}(u) = \sum_{j=1}^N P_j(u) (\theta-n-1/2) = \theta-n-1/2.
$$
Therefore
\ben
E\eta_i = \frac{2}{D-1-2n} ([\theta, A_i^{(n)}(u)]e,e).
\een
It remains only to use that $\theta(e)=(D/2) e$ and that $\theta$ is
skew-symmetric with respect to the pairing. 

Finally, the last condition that we have to check is 
\beq\label{pde-ijk}
\frac{\partial \eta_{ij}} {\partial u_k} =\frac{1}{2}
\Big(
\frac{\eta_{ik}\eta_{jk}}{\eta_k} + \frac{\eta_{ji}\eta_{ki}}{\eta_i}
+
\frac{\eta_{kj}\eta_{ij}}{\eta_j} \Big), \quad k\neq i\neq j\neq k.
\eeq
Let us explain how to express the LHS as a quadratic expression
in the functions $\eta_{ab}$. 
Recall that we have the following differential equation 
\ben
\partial_{u_j} P_i = \frac{1}{u_j-u_i} \Big( P_j\theta P_i-P_i\theta
P_j\Big). 
\een
Using the above differential equations and the fact that the operators
$P_a$ are self-adjoint and $\theta $ is skew 
symmetric with respect to $(\ ,\ )$ we get 
\beq\label{can-theta}
\eta_{ij} = (\partial_{u_j}P_i(u)e,e) = 
\frac{2}{u_i-u_j} (P_i(u)e,\theta P_j(u)e). 
\eeq
The derivative $\partial_{u_k}\eta_{ij}$ becomes 
\ben
\frac{2}{u_i-u_j} \Big( 
\frac{(P_k\theta P_i e, \theta P_j e)}{u_k-u_i}-
\frac{(P_k\theta P_je,\theta P_ie)}{u_k-u_j}+
\frac{(P_i\theta P_k e, \theta P_j e)}{u_i-u_k}  -
\frac{(P_j\theta P_k e, \theta P_i e)}{u_j-u_k}
\Big).
\een
Using the projection formula $P_i x = (x,P_ie) \frac{P_ie}{\eta_i}$ we
get 
\beq\label{trace-1}
\frac{(P_k\theta P_i e, \theta P_je)}{u_k-u_i} = 
(\theta P_ie , P_k e) (P_ke,\theta P_j e)\frac{1}{\eta_k} = 
\frac{\eta_{ik}\eta_{jk}}{4\eta_k}\,
(u_k-u_j).
\eeq
Similar formulas hold for the remaining 3 terms above, so for the
derivative  $\partial_{u_k}\eta_{ij}$ we get
\ben
\frac{2}{u_i-u_j} \Big(
\frac{\eta_{ik}\eta_{jk}}{4\eta_k}\,(u_k-u_j) - 
\frac{\eta_{ik}\eta_{jk}}{4\eta_k}\,(u_k-u_i) +
\frac{\eta_{ki}\eta_{ji}}{4\eta_i}\,(u_i-u_j) -
\frac{\eta_{kj}\eta_{ij}}{4\eta_j}\,(u_j-u_i)
\Big).
\een
The above expression is precisely the RHS of \eqref{pde-ijk}.
\qed

\medskip
The proof of the proposition can be completed as follows. Let us define the
isomorphism
\ben
T_{u^\circ}U \cong H,\quad \partial_{u_i}\mapsto P_i^\circ e, 
\een
where slightly abusing the notation we have denoted by $\partial_{u_i}$
the tangent vector in $T_{u^\circ}U$ representing the value of the
coordinate vector field $\partial_{u_i}$ at $u^\circ$. 
We claim that the special initial condition corresponding to the
Frobenius structure defined by Lemma \ref{sic-Frob} coincides with the
given special initial condition. The easiest way to see this is if we
fix the basis of $H$ to be $\phi_i=P_i^\circ e$. Then for the given special
initial condition we have: the matrix
$P_j^\circ$ is $E_{jj}$ (the matrix with 1 on place $(j,j)$ and 0
elsewhere), the matrix of the pairing $(\ ,\ )$ is diagonal with diagonal
entries  $(P_i^\circ e,e) = \eta_i(u^\circ)$, the vector $e$ has
coordinates $(1,\dots,1)$, and the matrix of $\theta$ becomes (see
formula \eqref{can-theta}) 
\ben
\theta_{ij}=(u_i^\circ-u_j^\circ)
\frac{\eta_{ij}(u^\circ)}{2\eta_i(u^\circ)},\quad 1\leq i,j\leq N.
\een
Comparing with the special initial condition corresponding to the
Frobenius structure we see that the only thing left to prove is that
the Hodge grading operator $\widetilde{\theta}|_{T_{u^\circ}U}$ coincides with
$\theta$.   Let us compute the matrix of $\widetilde{\theta}$ in
canonical coordinates. Note that $\widetilde{\theta}_{ij}=0$ for $i=j$
due to skew-symmetry. Let us assume that $i\neq j$. Then 
\ben
\widetilde{\theta}_{ij}(u)\eta_i(u) = (\partial_{u_i},\nabla^{\rm
  L.C.}_{\partial u_j} E) =\partial_{u_j}( \partial_{u_i},E) -
\sum_{k=1}^N \Gamma_{ij}^k(u) (\partial_{u_k},E),
\een
where $\Gamma_{ij}^k$ are the Christoffel's symbols of the Frobenius
pairing. Recalling the formulas for the Christoffel's symbols (see
Step 1 in the proof of Theorem \ref{thm:Frob-class}) we get  
\beq\label{theta-can}
\widetilde{\theta}_{ij}(u)\eta_i(u) =(u_i-u_j)\,
\frac{\eta_{ij}(u)}{2} \quad \Rightarrow \quad
\widetilde{\theta}_{ij}(u) =(u_i-u_j)\,
\frac{\eta_{ij}(u)}{2\eta_i(u)}. 
\eeq
Restricting to $u=u^\circ$ we get that
$\widetilde{\theta}(u^\circ)=\theta$. 
\qed

\subsection{The genus-1 potential}

We would like to finish this lecture by deriving the relation
between the genus-1 primary potential of the semi-simple Frobenius
structure and the isomonodromic tau-function. Following Givental
(see \cite{G2} and the references there in) we introduce the genus-1 potential as 
\ben
F^{(1)}(u) = \frac{1}{2}\int \sum_{i=1}^N R_1^{ii}(u) du_i - 
\frac{1}{48} \log (\eta_1(u)\cdots \eta_N(u)),
\een
where $R_1(u)$ is a matrix and $R_1^{ii}$ is the $(i,i)$-entry. The
function is called genus-1 potential, because in the case of
quantum cohomology of some manifold $X$, the above formula coincides
with the generating function of genus-1 Gromov--Witten invariants of
$X$. 

In order to define the matrix $R_1(u)$ we have to make a choice of
square root and define $\eta_i(u)^{1/2}$ for all $i$. Let $\Psi(u)$ be
the matrix with entries
\ben
\Psi_{ai}(u) = \widetilde{\Psi}_{ai} (u)\,\eta_i^{-1/2} = \frac{\partial
  t_a}{\partial u_i}\, \eta_i^{-1/2},
\een
where $t=(t_1,\dots,t_N)$ is a flat coordinate system. 
Dubrovin's connection $\nabla$ has a unique formal asymptotic solution near
$z=0$ of the form 
\ben
\Psi(u)(1+R_1(u) z+R_2(u) z^2+\cdots ) e^{U/z},\quad
U=\diag(u_1,\dots,u_N). 
\een
Substituting this formula in the differential equation
$\nabla_{\partial_z}J=0$ and recalling Lemma \ref{le:psi} we get
\ben
z\partial_z R(u,z) + z^{-1} [U,R(u,z)] = V(u) R(u,z),
\een
where $V(u) := \Psi(u)^{-1} \theta \Psi(u)$. Comparing the
coefficients in front of $z^k$ we get 
\ben
k R_k + [U,R_{k+1}] = V R_k,\quad k\geq 0.
\een
Since we work with a Frobenius structure defined on an open subset of
$T$, the diagonal entries of $U$ are pairwise distinct, so the above
recursion has a unique solution. In particular, for the $(i,j)$-entry
of $R_1$ we get 
\ben
R_1^{ij}(u) = \frac{V_{ij}(u)}{u_i-u_j},\quad \mbox{if $i\neq j$} 
\een
and 
\ben
R_1^{ii}(u) = \sum_{j:j\neq i} V_{ij}(u) R_1^{ji}(u) = 
-\sum_{j:j\neq  i}  \frac{V_{ij}(u) V_{ji}(u) }{u_i-u_j}. 
\een
By definition $V$ is the matrix of the Hodge grading operator in the
orthonormal basis $e_i:=\eta_i^{-1/2}\partial_{u_i}$, i.e.,
\ben
\theta(e_j) = \sum_{i=1}^N V_{ij}(u) e_i. 
\een
On the other hand we have already computed the matrix of $\theta$ in
the canonical basis (see formula \eqref{theta-can}). Therefore
\ben
V_{ij} (u)= (u_i-u_j) \frac{\eta_{ij}(u)}{2\eta_i(u)^{1/2}\eta_j(u)^{1/2}}.
\een
Finally, for the 1-form $\sum_{i=1}^N R_1^{ii}(u) du_i $ we get 
\ben
\sum_{i=1}^N \sum_{j:j\neq i}
(u_i-u_j)\frac{\eta_{ij}(u)^2}{4\eta_i(u) \eta_j(u)} du_i = 
\frac{1}{8}\, \sum_{i=1}^N \sum_{j:j\neq i}
(u_i-u_j)\frac{\eta_{ij}(u)^2}{\eta_i(u) \eta_j(u)} (du_i-du_j).
\een
Note that the above form is independent of the choice of a square root
used in the definition of $R_1$. 
Let us compare this form with the 1-form $\omega$ defining the
isomonodromic $\tau$-function of the second structure connection
$\nabla^{(n)}$. We have 
\ben
\operatorname{tr}(A_i^{(n)}(u)A_j^{(n)}(u)) =
\operatorname{tr}(P_i\theta P_j \theta) -(n+1/2)
\operatorname{tr}(P_i\theta+\theta P_j). 
\een
The first trace on the RHS is 
\ben
\eta_i^{-1}(P_i\theta P_j \theta P_ie, e) = -(u_i-u_j)^2
\frac{\eta_{ij}^2}{4\eta_i \eta_j},
\een
where we used formula \eqref{trace-1}. The second trace is 0 because
\ben
\operatorname{tr}(P_i\theta) = \eta_i^{-1}(P_i \theta P_i e, e) = 
\eta_i^{-1}(\theta P_i e, P_ie) = 0, 
\een
where the last equality uses the fact that $\theta$ is skew-symmetric
with respect to the Frobenius pairing. We get
\ben
\omega = \frac{1}{2}\sum_{i=1}^N\sum_{j:j\neq i} 
\frac{\operatorname{tr}(A_i^{(n)}(u)A_j^{(n)}(u)) }{u_i-u_j} \,
(du_i-du_j) = -\sum_{i=1}^N R_1^{ii}(u) du_i.
\een
Finally we get the following relation 
\ben
e^{-48 F^{(1)}(u) }= \tau(u)^{24} \, \eta_1(u)\cdots \eta_N(u). 
\een

\medskip

\noindent
{\bf Acknowledgement.} I would like to thank S. Galkin and
A. Getmanenko for organizing a seminar at Kavli IPMU where we studied
Bolibvruch's counter example to the Riemann-Hilbert problem. This work
is partially supported by JSPS Grant-In-Aid 26800003  
and by the World Premier International Research Center Initiative (WPI
Initiative),  MEXT, Japan.

\bibliographystyle{amsalpha}

\end{document}